\documentclass[11pt]{article}

\usepackage{latexsym}
\usepackage{amsmath}
\usepackage{amssymb}
\usepackage{amsfonts}
\usepackage{graphicx,epsfig}
\usepackage{color}

\newcommand{\bitem}{\begin{itemize}}
\newcommand{\eitem}{\end{itemize}}
\newcommand{\beq}{\begin{equation}}
\newcommand{\eeq}{\end{equation}}
\newcommand{\goto}{\rightarrow}
\newcommand{\cC}{{\cal C}}

\newcommand{\cA}{{\cal A}}
\newcommand{\cP}{{\cal P}}
\newcommand{\cS}{{\cal S}}
\newcommand{\cB}{{\cal B}}
\newcommand{\bR}{{\bf R}}
\newcommand{\bP}{{\bf P}}
\newcommand{\bZ}{{\bf Z}}

\newcommand{\cR}{{\cal R}}
\newcommand{\cT}{{\cal T}}
\newcommand{\cH}{{\cal H}}
\newcommand{\cN}{{\mathcal N}}
\newcommand{\eps}{{\varepsilon}}
\newcommand{\argmin}{\mbox{argmin}}
\newcommand{\SgSp}{\mbox{sing supp}}
\newcommand{\Sp}{\mbox{supp}}
\newcommand{\cK}{{\cal K}}
\newcommand{\cL}{{\cal L}}

\def\t{\tilde}
\newcommand{\ip}[2]{\langle#1,#2\rangle}

\newcommand{\absip}[2]{| \langle#1,#2\rangle |}
\newcommand{\norm}[1]{\|#1\|}
\newcommand{\qed}{$\Box$}

\definecolor{cardinal}{rgb}{.64,0.,.11}

\numberwithin{equation}{section}
\newtheorem{theorem}{Theorem}[section]
\newtheorem{corollary}{Corollary}[section]
\newtheorem{lemma}{Lemma}[section]
\newtheorem{proposition}{Proposition}[section]
\newtheorem{definition}{Definition}[section]
\newcommand{\cF}{{\cal F}}

\title{Microlocal Analysis of the \\
Geometric Separation Problem}
\author{David L. Donoho and Gitta Kutyniok}

\begin{document}
\maketitle

\begin{abstract}
Image data are often composed of two or
more geometrically distinct constituents;
in galaxy catalogs, for instance,
one sees a mixture of pointlike structures (galaxy superclusters)
and curvelike structures (filaments). It would be ideal
to process a single image and extract two geometrically `pure' images,
 each one containing features from only one of the two  geometric constituents.
This seems to be a seriously underdetermined problem, but
recent empirical work achieved
highly persuasive  separations.

We present  a theoretical analysis showing that
accurate geometric separation of point and curve singularities
can be achieved by minimizing the $\ell_1$
norm of the representing coefficients in
two geometrically complementary frames:
wavelets and curvelets.
Driving our analysis is a specific property
of the ideal (but unachievable) representation
where each content type is expanded in the frame best adapted to it.
This ideal representation has
the property that important coefficients
are {\it clustered geometrically} in phase space, and that at fine
scales, there is very little coherence between
a cluster of elements  in one frame expansion and
individual elements in the complementary frame.
We formally introduce notions of {\em cluster coherence} and
{\it clustered sparsity} and use this machinery to show that the
underdetermined systems of linear equations can be stably solved by
$\ell_1$minimization; microlocal phase space helps
organize the calculations that cluster coherence
requires.

\end{abstract}

\vspace{.1in}
Mutual Coherence. Cluster Coherence. Tight Frames. Curvelets, Shearlets,
Radial Wavelets.
\vspace{.1in}

\vspace{.1in}
{\bf Acknowledgements.}
The authors would like to thank Inam ur Rahman,
Apple Computer, for graphics help, and Emmanuel Cand\`es,
Michael Elad, and Jean-Luc Starck, for numerous discussions on related topics.
The second author would like to thank the Statistics Department  at Stanford
and the Mathematics Department at Yale for hospitality
and support during her visits. Thanks also to the Isaac Newton  Institute of Mathematical  Sciences in Cambridge, UK for
an inspiring research environment which led to the completion of a
significant part of this work. This work was partially supported
by NSF DMS 05-05303 and DMS 01-40698 (FRG), and by Deutsche Forschungsgemeinschaft
(DFG) Heisenberg fellowship KU 1446/8-1.

\vspace{.1in}

\pagebreak

\section{Introduction}

Cosmological data analysts face tasks of {\it geometric separation} \cite{SCD03,SED04}.
Gravitation, acting over time,
drives an initially quasi-uniform distribution of matter
in 3D to concentrate near lower-dimensional structures:  points,  filaments,
and sheets.
It would be desirable to process single `maps' of matter density and somehow
extract three `pure'  maps containing just the points,  just the
filaments, and  just the sheets around which matter is concentrating.

In seemingly unrelated fields, such as medical imaging and materials
science, related questions arise frequently and naturally.  For example,
a technologist with a single confused image of an aggregate
might wish to create two images, one
containing just the fibrous and the other just the granular structures, respectively.

Such `desires', when voiced by a working scientist or engineer,
really amount to a request
for existing information technology to be put to work here and now
on data available today.  No doubt there is a wide spectrum
of image processing `hacks' and `improvisations' that might be useful,
on a case-by-case basis.
The mathematician's interest would only be piqued when
an intellectually coherent approach shows promise of success,
especially if the reasons for success are subtle and instructive.

Recently,  astronomer Jean-Luc Starck
and collaborators have been empirically successful in
numerical experiments with component
separation;  their approach
used tools from  modern harmonic analysis
in a provocative way.  They used two or more
overcomplete frames, each one specially adapted
to particular geometric structures, and
were able to obtain separation despite the
fact that the underlying system of equations is
highly underdetermined.
Here we analyze such approaches in a mathematical framework
where we can show that success stems from an
interplay between geometric properties of the objects to be separated,
and the harmonic analysis for singularities of various geometric types.
We eventually point to a much wider range of seemingly very different
`imaging' problems where our analysis techniques can provide insight.

\subsection{Singularities and Sparsity}

As a mathematical idealization
of  `image', consider a
Schwartz distribution  $f$ with domain $\bR^2$.
The distribution $f$ will be given
singularities with specified geometry: points and curves.

We plan to represent such an `image' using tools of harmonic
analysis; in particular, bases and frames.   While many such
representations are conceivable, we are interested here just in
those
 bases or frames which  can sparsely represent
$f$ -- i.e., can represent $f$ using relatively few large coefficients.

The type of basis which best sparsifies  $f$ depends on
the geometry of its singularities. If the singularities
occur at a finite number of (variable)
points, then {\it wavelets} give what is, roughly speaking, an
optimally sparse representation -- one with
the fewest significantly nonzero coefficients.
If the singularities occur at a finite number of
smooth curves, then  one of the recently
studied directional multiscale representations
({\it curvelets} or {\it shearlets})
will do the best job of sparsification.
(For careful quantitative discussions
of sparsification see,  e.g., \cite{CD04} etc.).

In fact, real-world signals
are, generally speaking, a mixture of content types and, correspondingly,
a model where singularities are of only one geometric type
is overly narrow.  If $f$ is actually a nontrivial superposition
$\cP + \cC$ where  $\cP$ has only point singularities and
$\cC$ has only curvilinear singularities, then two things happen:
\bitem
\item Neither wavelets {\it alone} nor curvelets {\it alone} will be very good
for representing $\cP + \cC$. The sparsity
either achieves alone is much less satisfactory than the {\it ideal} sparsity level -- that
which could be  achieved by using wavelets for representing $\cP$ and by curvelets in representing $\cC$. (This ideal
representation is purely notional; it assumes one can first perfectly separate
the two objects and then separately analyze the separated layers.)
\item In fact, no single basis or traditional
linear representation is very good at sparsifying $\cP+\cC$ compared
to the ideal representation.
\eitem

This immediately suggests the need to use {\it both} systems
to represent $f$ sparsely; however, since each system
is itself complete (or even overcomplete) there is no
obvious traditional way to do this.

In this paper, we consider the problem of developing sparse
representations by combining both wavelets and curvelets
and using a nonlinear representation based on $\ell_1$ minimization.
The problem we solve is a continuum variant of a problem in
image and signal processing
with considerable practical
interest, and extensive work for almost two decades.
For  references, see Subsection \ref{sec-DiscreteHistory} below.
That work, while suggestive and inspiring, concerns discretely indexed
signal/image processing, obscuring the {\it continuum} elements of geometry and
microlocal analysis which are essential to this paper.

\subsection{A Geometric Separation Problem}

Consider the following simple but
clear model problem of geometric separation.
Consider a `pointlike' object  $\cP$
made of point singularities:
\begin{equation} \label{pointdef}
   \cP = \sum_{i=1}^P |x - x_i|^{-3/2} .
\end{equation}
This object is smooth away from the $P$ given points $(x_i : 1 \leq i \leq P)$.
Consider as well a `curvelike' object
$\cC$, a singularity along a closed curve $\tau : [0,1] \mapsto \bR^2$:
\begin{equation} \label{curvedef}
    \cC = \int \delta_{\tau(t)}(\cdot)dt,
\end{equation}
where $ \delta_x$ is the usual Dirac delta function located  at $x$.
The singularities underlying these two distributions are geometrically quite different,
but the exponent $\-3/2$ is chosen so the  energy distribution across scales is similar;
if $\cA_r$ denotes the annular region $r < |\xi| < 2r$,
\begin{equation}\label{energyMatch}
      \int_{\cA_r} |\hat{\cP}|^2(\xi) \asymp r  ,
\qquad \int_{\cA_r} |\hat{\cC}|^2(\xi) \asymp r, \qquad r \goto \infty .
\end{equation}
This choice makes the components comparable
as we go to  finer scales;
the ratio of energies is more or less independent of scale.
 Separation is challenging at {\it every} scale.

Now assume that we observe
the `Signal'
\begin{equation}\label{superpose}
     f = \cP + \cC,
\end{equation}
however, the component distributions $\cP$ and $\cC$
are unknown to us.

\begin{definition}
The {\em Geometric Separation Problem}
requires to recover $\cP$ and $\cC$ from know\-ledge only of $f$;
here $\cP$ and $\cC$ are unknown to us, but obey (\ref{pointdef}),
(\ref{curvedef}) and certain regularity conditions on the curve $\tau$.
\end{definition}

As there are two unknowns ($\cP$ and $\cC$) and only
one observation ($f$), the problem seems improperly posed.
We develop a principled, rational approach which provably
solves the problem according to clearly stated standards.




\subsection{Two Geometric Frames }
\label{sec:twogeometricframes}

We now focus on two overcomplete systems  for representing the object $f$:
\bitem
\item {\it Radial Wavelets} -- a  tight frame with perfectly isotropic
generating elements.
\item {\it Curvelets} -- a highly directional tight frame with increasingly
anisotropic elements at fine scales.
\eitem
We pick these because, as is well known, point singularities
are coherent in the wavelet frame and curvilinear singularities
are coherent in the curvelet frame.
In Section \ref{sec:extensions} we discuss other system pairs. For readers
not familiar with frame theory, we refer to \cite{Chr03}, where terms
like `tight frame' --  a Parseval-like property --  are carefully discussed.

The point- and curvelike objects we defined in the previous subsection are real-valued distributions.
Hence, for deriving sparse expansions of those, we will consider radial wavelets and curvelets
consisting of real-valued functions. So  only angles associated with
radians $\theta \in [0,\pi)$ will be considered, which later on we will, as is customary, identify
with $\bP^1$, the real projective line.

We now construct the two selected tight frames as follows. Let $W(r)$ be an `appropriate' window function, where
in the following we assume that $W$ belongs to $C^\infty(\bR)$ and is compactly supported on
$[-2,-1/2] \cup [1/2,2]$ while being the Fourier transform of a wavelet. For instance, suitably scaled
Lemari\`{e}-Meyer wavelets possess these properties. We define {\em continuous radial wavelets}
at scale $a > 0$ and spatial position $b \in \bR^2$ by their Fourier transforms
\[
      \hat{\psi}_{a,b}(\xi) =  a \cdot W(a|\xi|) \cdot \exp\{ i b'\xi\}.
\]
The {\em wavelet tight frame} is then defined as a sampling of $b$ on a series of
regular lattices $\{ a_j  \bZ^2 \}$, $j \geq j_0$,
where $a_j = 2^{-j}$, i.e.,  the radial wavelets at scale $j$ and spatial position
$k = (k_1,k_2)'$ are given by the Fourier transform
\[
      \hat{\psi}_{\lambda}(\xi) =  2^{-j} \cdot W(|\xi|/2^{j}) \cdot \exp\{ i k'\xi /2^j\},
\]
where we let $\lambda = (j,k)$ index position and scale.

For the {\it same} window function $W$ and a `bump function' $V$, we define {\em continuous
curvelets} at scale $a>0$, orientation $\theta \in [0,\pi)$, and spatial position $b \in \bR^2$
by their Fourier transforms
\[
      \hat{\gamma}_{a,b,\theta}(\xi) =  a^{\frac{3}{4}}  \cdot W(a|\xi|) V(a^{-1/2}(\omega-\theta))
         \cdot \exp\{ i b' \xi \}.
\]
See \cite{CD05a} for more details.
The {\em curvelet tight frame} is then (essentially) defined as a sampling of $b$ on a
series of regular lattices
\beq \label{eq:defilattice}
   \{ R_{\theta_{j,\ell}} D_{a_j}  \bZ^2 \}, \qquad j \geq j_0, \quad \ell = 0, \dots, 2^{\lfloor j/2 \rfloor} -1 ,
\eeq
where $R_{\theta}$ is planar rotation by $\theta$ radians, $a_j = 2^{-j}$,
$\theta_{j,\ell} = \pi \ell / 2^{j/2}$, $\ell = 0, \dots, 2^{j/2}-1$,
and  $D_a$ is anisotropic dilation by $diag(a,\sqrt{a})$, i.e., the curvelets
at scale $j$, orientation $\ell$, and spatial position $k = (k_1,k_2)$ are
given by the Fourier transform
\[
      \hat{\gamma}_{\eta}(\xi) =  2^{-j\frac{3}{4}}  \cdot W(|\xi|/2^{j}) V((\omega-\theta_{j,\ell})2^{j/2})
         \cdot \exp\{ i (R_{\theta_{j,\ell}}D_{2^{-j}}k)' \xi \},
\]
where let $\eta = (j,k,\ell)$ index scale, orientation, and scale.
(For a precise statement, see \cite[Section 4.3, pp. 210-211]{CD05b}).

Roughly speaking, the radial wavelets are `radial bumps' with position $k/2^{j}$
and scale $2^{-j}$, while the curvelets live on anisotropic regions of width
$2^{-j}$ and length $2^{-j/2}$.  The wavelets are good at representing point singularities
while the curvelets are good at representing curvilinear singularities.

Using the {\it same} window $W$, we can construct a family of
filters $F_j$ with transfer functions
\[
    \hat{F}_j(\xi) = W(|\xi|/2^{j}), \qquad \xi \in \bR^2 .
\]
These filters allow us to decompose a function $f$
into pieces $f_j$ with different scales,
the piece $f_j$ at subband $j$
arises from filtering $f$ using $F_j$:
\[
   f_j = F_j \star f;
\]
the Fourier transform $\hat{f}_j$
is supported in the annulus with inner radius $2^{j-1}$ and outer radius
$2^{j+1}$. Because of our assumption on $W$,
we can reconstruct the original function from these
pieces using the formula
\[
   f = \sum_j F_j \star f_j, \qquad f \in L^2(\bR^2) .
\]

The tight frames of curvelets and radial wavelets
discussed above interact in a very local way with the
filtering $F_j$.

{\bf Lemma.}
{\sl Let $\cF_j $ denote the range of the operator of convolution
with $F_j$.  Then  curvelets at level $j'$ are orthogonal to
$\cF_j$ unless $|j' -j | \leq 1$.  Similarly, radial wavelets
at level $j'$ are orthogonal to $\cF_j$ unless $|j' -j | \leq 1$.
}

{\bf Proof.} Indeed, $\cF_j$ is the collection of all functions $f_j$
whose Fourier transform is representable as
$\hat{f}_j(\xi) = W(|\xi|/2^{j}) \hat{f} (\xi)$  where $f \in L^2(\bR^2)$.
The support in frequency space of elements of $\cF_j$
is thus an annulus $\cA_j$ (say).  The annuli
have disjoint interiors if $|j-j'| > 1$.
Hence $\cF_{j} \perp \cF_{j'}$ if $|j'-j| >1$.

However, both the radial wavelet frame elements and the
curvelet frame elements  at level $j'$  belong to $\cF_{j'}$.
 \qed

For future use,  let $\Lambda_j$ denote the collection of indices $(j,k)$ of wavelets
at level $j$, and
\[
\Lambda_j^{\pm n} = \bigcup_{j'=j-n}^{j+n} \Lambda_{j'}.
\]
Similarly,
let $\Delta_j$ denote the indices $\eta = (j,k,\ell)$ of curvelets at level $j$,
and let
\[
\Delta_j^{\pm n} = \bigcup_{j'=j-n}^{j+n} \Delta_{j'}.
\]
We conclude that elements of $\cF_j$ can be represented using either
only radial wavelets $\{\psi_\lambda : \lambda \in \Lambda_j^{\pm 1}\}$
or only curvelets $\{\gamma_\eta: \eta \in \Delta_j^{\pm 1} \}$.

\subsection{Separation via $\ell_1$ Minimization}

\subsubsection{Sparse Multiple Frame Expansions}
\label{sec:sparsemultipleframeexpansions}

We now have two complete representations
for $\cF_j$, yielding two ways of representing
the subband component $f_j$:
in terms of its wavelet expansion:
\[
f_j = \sum_{\lambda \in \Lambda_j^{\pm 1}} w_\lambda \psi_\lambda ;
\]
or in terms of its curvelet expansion:
\[
f_j = \sum_{\eta \in \Delta_j^{\pm 1}} c_\eta \gamma_\eta .
\]

Each frame exhibits a single geometric tendency -- either
highly nondirectional or highly directional -- in representing
$f_j$. However, $f_j$ may have both isotropic and directional
features.  We therefore seek a combined representation
\[
f_j  = \sum_{\lambda \in \Lambda_j^{\pm 1}} w_\lambda \psi_\lambda
+ \sum_{\eta\in \Delta_j^{\pm 1}}  c_\eta \gamma_\eta  .
\]
Because the combined frame formed by concatenating the two
frames is  overcomplete, there are many possible ways this
decomposition can be done.  Some of them may be geometrically
motivated, many are not.

Consider the following dual-frame
{\bf Component Separation}  problem based on $\ell_1$ minimization:

\begin{eqnarray*}
(\mbox{\sc CSep})  \qquad   (W_j,C_j) &=& \mbox{ argmin } \|w\|_1 + \|c\|_1 \\
              && \quad \mbox{subject to }  f_j = W_j + C_j\\
              && \qquad \mbox{    and    }  w_\lambda = \langle W_j , \psi_\lambda \rangle, \quad \lambda \in \Lambda_j^{\pm 1} \\
              && \qquad \mbox{    and    }  c_\eta = \langle C_j , \gamma_\eta \rangle,  \quad \eta \in \Delta_j^{\pm 1}  .
\end{eqnarray*}

In words, we take a given scale subband $f_j$ and
decompose it into a wavelet component $W_j$ and
a curvelet component $C_j$. The components
are chosen by the principle of $\ell_1$ minimization
on the frame coefficients: the $\ell_1$ norm of the wavelet coefficients
of the wavelet component should be small, and the $\ell_1$ norm of the
curvelet coefficients of the curvelet component should be small.

Here is our reason for the `component separation'  label:
Armed with the optimization result at each scale subband, we
define the purported {\it pointlike} component as the superposition
of all the wavelet terms:
\[
    \t{P} = \sum_j F_j \star W_j ;
\]
and the purported {\it curvelike} component as the superposition
of all the curvelet terms:
\[
    \t{C} = \sum_j F_j \star C_j .
\]

We obtain the decomposition
\[
   f = \t{P} + \t{C}.
\]

\subsubsection{Main Result}


At this stage, we have two decompositions: one by
the truly geometric pair $(\cP,\cC)$ of pointlike and curvelike
objects and one by the {\it purported} geometric pair $(\t{P},\t{C})$.
The following result justifies our interest in
the second pair.
To state it, define the scale subbands of the truly geometric components by:
\[
  \cP_j = F_j \star \cP ; \qquad \cC_j = F_j \star \cC.
\]


\begin{theorem} {\sc Asymptotic Separation.} \label{theo:maintheorem}
\begin{equation} \label{eq:mainresult}
\frac{ \| W_j - \cP_j \|_2 + \| C_j - \cC_j \|_2 }{\| \cP_j\|_2 + \|\cC_j\|_2 } \goto 0, \qquad j
\goto \infty .
\end{equation}
\end{theorem}

At fine scales, the truly pointlike component is
almost all captured by the wavelet component
and the truly curvelike component is almost all captured by the curvelet
component.  In short, the purported
pointlike and curvelike  components deserve the labelling they have been
given.

\subsection{Extensions}

Theorem 1.1  is amenable
to generalizations and extensions. Previewing
Section \ref{sec:extensions}, we mention a few examples.

\bitem
\item {\it More General Classes of Objects.}
Theorem \ref{theo:maintheorem} can be generalized to other
situations.  First, we could consider singularities of different orders.
This would allow $\cC$ to model `cartoon' images, where
the curvilinear singularities are now the boundaries of the pieces
for piecewise $C^2$ functions.
Second, we can allow smooth perurbations, i.e.,
$f = (\cP + \cC + g) \cdot h $ where $g, h$ are smooth functions of rapid decay at $\infty$.
In this situation, we let the denominator in (\ref{eq:mainresult}) be simply $\| f_j \|_2$.

\item {\it Other Frame Pairs.}
Theorem \ref{theo:maintheorem} holds without change for many other
pairs of frames and bases, such as, e.g., orthonormal separable Meyer wavelets
and shearlets.

\item {\it Noisy Data.}
Theorem \ref{theo:maintheorem} is resilient to noise impact; an image
composed of $\cP$ and $\cC$ with additive `sufficiently small' noise exhibits the same
asymptotic separation.

\item {\it Rate of Convergence.}
Theorem \ref{theo:maintheorem} can be accompanied by explicit decay
estimates.

\item {\it Other Algorithms and Other Notions of Separation.}
In the companion paper \cite{DK08a} we study thresholding as an alternative approach
to separation; it is less computationally demanding than the $\ell^1$ minimization studied here,
but also somewhat less elegant. Building on the estimates proved in this paper,
\cite{DK08a} shows that properly-tuned thresholding can also achieve asymptotic separation.
\eitem

\subsection{The Multiple-Basis Representation Problem}
\label{sec-DiscreteHistory}
Theorem 1.1 should be placed in context of a great deal of ongoing work
concerning sparsity and overcomplete representations.
Already in the early 1990's, R.R. Coifman
became interested in the problem of representing
discrete-time signals using more than one basis.  In a conversation,
he told one of us about a problem which,
in retrospect and using modern formulations,
can be posed as follows:

\bitem
\item An observed signal  $S \in \bR^n$ is thought to be a superposition of subsignals $S_i$, $i=1,2$.
\item Each subsignal $S_i$ is thought to be `coherent' in an `appropriate' basis $\Phi_i$, $i=1,2$.
\item Each subsignal 'looks incoherent' in an `inappropriate' basis. Here $\Phi_2$ is
inappropriate for $S_1$, and $\Phi_1$ is inappropriate for $S_2$.
\eitem

Coifman, Wickerhauser and co-workers at the time made a sort of heuristic exploration
motivated intuitively by these slogans.
As a published example of their work at the time, please see \cite[Fig. 26(a-h)]{CW93}.
The different `coherent parts' displayed in those figures were obtained by the following recipe:

\begin{enumerate}
\item Transform signal $S$ into basis $\Phi_1$.
\item Threshold the coefficients, yielding sparse coefficients $\tilde{\alpha}_1$.
\item Form residual $R = S - \Phi_1 \tilde{\alpha}_1$.
\item Transform $R$ into basis $\Phi_2$.
\item Threshold the coefficients, yielding sparse coefficients $\tilde{\alpha}_2$.
\item Write $\tilde{S}_i = \Phi_i \tilde{\alpha}_i$; then
\begin{equation}
\label{Recipe}
   S = \tilde{S}_1 + \tilde{S}_2 + residual.
\end{equation}
\end{enumerate}

At about the same time, St\'{e}phane Mallat and Zhifeng Zhang  became
interested in the problem of representing signals using a highly overcomplete
dictionary of time-frequency atoms; \cite{MZ93} (their
dictionary had  $\approx \log(N)$ different frames,
where $N$ is the signal length). Their approach, called {\it Matching Pursuit},
built up an approximation one-term-at-a-time iteratively,
at each stage finding the best single atom in any of the
several bases which was not yet already
forming part of the approximation and adding that term to
the approximation.

Lurking in these early numerical experiments  were two larger questions.
If there truly is a simple representation of the signal using
more than one system, can it ever be found?  Can it be found by such
a simple approach?

Formally: can one accurately recover `coherent' pieces $S_1$ and $S_2$
given knowledge of $S = S_1 + S_2$ only? For example, can we expect
that the outputs $\tilde{S}_{1}$, $\tilde{S}_{2}$  in (\ref{Recipe}) obey
 $\tilde{S}_1 \approx S_1$ and $\tilde{S}_2 \approx S_2$?
Researchers at the time said in conversation, that, when put this starkly,
the answer was simply `no', since there are twice as many unknowns as
knowns. Nevertheless, some of the empirical results at the
time were suggestive and inspiring.

\subsection{Minimum $\ell_1$ Decomposition and Perfect Separation}

A few years later, one of us worked with Scott Shaobing Chen  to
develop a formal, opti\-mi\-za\-tion-based approach
to the multiple-basis representation problem.
Given bases $\Phi_i$, $i=1,2$, one solves the following
problem
\[
(\mbox{\sc BP}) \qquad     \min \|\alpha_1\|_1 + \|\alpha_2\|_1
    \mbox{  subject to  }  S = \Phi_1 \alpha_1 + \Phi_2 \alpha_2.
\]
Here $\| \cdot \|_1$ denotes the usual $\ell_1$ norm.
Note that here there are $2n$ unknowns in $\alpha_1$
and $\alpha_2$ and only $n$ knowns in $S$, but that
an optimization principle is being used to select a
particular element from the $n$-dimensional
space of all possible solutions. (Terminological note:
the name 'Basis Pursuit' is meant to remind the reader
that (BP) actually selects a basis for the solution
out of the many conceivable bases which can be extracted
from the union of the two overcomplete systems).

Based on earlier experience of our first-named author and his collaborators,
see \cite{DJHS92,DL92,DS89},
it was known that the $\ell_1$ norm had a tendency to find sparse solutions
when they exist. And indeed, Chen's thesis showed that in some simple special
cases that this was so. Letting $\Phi_1$ be the standard basis
of $\bR^n$ (i.e., Kronecker sequences or
`spikes') and $\Phi_2$ be the Fourier basis, Chen
considered signals $S$ which were superpositions
of two spikes and two sinusoids. He showed that (BP)
recovered exactly the indices and coefficients of the terms involved in the synthesis;
and that this was true across a wide range of amplitude
ratios between the sinusoid and spike components.
In short, there was perfect separation of sinusoids
from spikes, and the true underlying simplicity of the
signal was revealed -- even though there were more unknowns than
equations.

In the years since that work, two streams of research emerged.
\bitem
 \item {\it Theoretical work}, showing that, indeed,
 one could in certain settings obtain the sparsest possible representations
  to an underdetermined problem by $\ell_1$ optimization; see, e.g.,
  \cite{CRT06b,Don06c,Don06b,Don06a,DET06,Tro04}
  for a selection of general work concerning $\ell_1$ minimization, and
  \cite{BGN08,DE03,DH01,EB02,GN03} for work somewhat relevant to component separation.
 \item {\it Empirical work}, showing that combined representations
  such as wavelets with curvelets or wavelets with sinusoids often
  gave very compelling separations of real signals and images, see, for instance,
  \cite{BSFMD07,CDS01,ESQD05,GB03,MAC02,SMBED05,SED04,SED05,SNM03,Tes07,KT09,ZP01}.
\eitem
We have already mentioned the empirical successes of Starck and collaborators.
For an overview of much recent work on sparse decompositions, see \cite{BDE09}.
Note that geometric separation is somewhat different
from the task of separation of {\it texture from smooth structure};
in that problem, sparsity in frame expansions does not play an
explicit role, nor do the geometric considerations which are
so important here; very interesting early work in such non-geometric separation was published by
Yves Meyer \cite{Meyer}, and Vese and Osher \cite{Osher}.

\subsection{Theorem 1.1 in Context, and Outline of Paper}

We can now place our result in context,
via several comparisons and contrasts,
looking ahead to themes developed below.
\bitem
\item {\it Microlocal Viewpoint.}
In Theorem 1.1  the objects of interest are
collections of point and curve singularities. The
viewpoint derives from microlocal analysis (see Section 3),
which says that points and curves are very different
objects in their joint space/orientation structure, so that
even if they happen to overlap spatially, they are microlocally distinct.
In contrast, other work on sparsity and $\ell_1$ minimization
typically has a discrete flavor, making hypotheses about the number of
nonzeros in an expansion and assuming the dictionary elements
interact randomly.
\item {\it Microlocal Asymptotics.}
Asymptotics are important for Theorem 1.1;
 the sharp separation between curves
and points in microlocal phase space
exists only as  a  limit phenomenon,
as the scale tends to zero. Asymptotic statements are
important in other literature on sparsity-driven
decompositions, but they are asymptotic in the number of random
elements in the underlying matrix, and exploit law-of-large-numbers
and concentration-of-measure effects.
For Theorem 1.1, such principles
play no role.
\item {\it Clustered Sparsity.}
In other work on sparsity-driven decompositions,
sparsity of the coefficients
plays a role primarily through the number of nonzeros.
In this work sparsity plays a role also through the arrangement of nonzeros;
Section \ref{sec:l1minimization} introduces a notion of clustered sparsity
and Sections 4-7 develop estimates bounding the locations of
significant nonzeros in the wavelet expansion of a point singularity or
in the curvelet expansion of a curve singularity; the estimates will be
organized using micolocal phase space ideas described in Section 3.
\item {\it Cluster Coherence.} In other work on sparsity-driven decompositions,
coherence or restricted isometry principles play a role;
these don't depend on the arrangement of nonzeros in an expansion.
Moreover, these are often applied to random-dictionary situations where the
interaction between frame elements is random and quasi arbitrary.
Section \ref{sec:l1minimization}   develops the notion of {\it cluster coherence} which
specifically depends on the arrangement of nonzeros.
We apply this notion to a dictionary  (wavelets + curvelets)
where interactions are geometrically driven,
and we develop estimates motivated by microlocal analysis which
provide the needed geometrical information.
\eitem

Our paper begins with Sections 2 and 3
introducing the driving ideas of
clustered sparsity,  cluster coherence,
and microlocal separation,
and describing a plan to prove Theorem 1.1
by establishing several needed estimates.
Later Sections 4, 5, 6, and 7 then develop these estimates,
by developing results about wavelet and curvelet expansions
of point and line singularities.
Section 8 then mentions a number of possible
extensions.

\section{Component Separation by $\ell_1$ Minimization}
\label{sec:compsep}

We now study  the behavior of $\ell^1$ minimization
in the two-frame case. Our analysis centers on
the use of cluster coherence to control joint concentration.

\subsection{$\ell_1$ Minimization for Separation of Two Tight Frames}
\label{sec:l1minimization}

Suppose we have two tight frames $\Phi_1$, $\Phi_2$
in a Hilbert space $\cH$, and a signal vector $S \in \cH$.
We know {\it a priori} that there exists a decomposition
\[
        S = S_1^0 + S_2^0,
\]
where $S_1^0$ is sparse in Frame 1, and
$S_2^0$ is sparsely represented in Frame 2.

Consider the following optimization problem
\begin{equation} \label{PSep}
    (\mbox{Sep}) \qquad    (S_1^\star, {S}_2^\star)
        = \argmin_{S_1,S_2}  \| \Phi_1^T S_1 \|_1 + \| \Phi_2^T S_2\|_1
          \mbox{ subject to } S = S_1 + S_2.
\end{equation}

The {\it optimization problem}
$(\mbox{\sc Sep})$ is visibly similar to,
but subtly different from, $(\mbox{\sc BP})$.
Here the $\ell_1$ norm
is being applied on the {\it analysis} coefficients of the two different
`components' rather than on the individual {\it synthesis} coefficients.
The hope in $(\mbox{\sc BP})$ is to get exactly the right nonzero coefficients in the sense of
those providing the sparsest representation.  However, this can become numerically
unstable for certain tight frames.
The hope in $(\mbox{\sc Sep})$ is merely to separate components
rather than the more ambitious goal of identifying
the true nonzero coefficients within each component's representation.
Starck and Elad
have  found this distinction to be important in their
own empirical work on separation.


To analyze this we need the following notion.
\begin{definition}
Let $\Phi_1$ and $\Phi_2$ be two tight frames. Given two sets of
coefficients $\cS_1$ and $\cS_2$, define the {\em joint
concentration} $\kappa=\kappa(\cS_1,\cS_2)$ by
\[ \kappa(\cS_1,\cS_2) =
\sup_f \frac{\norm{1_{\cS_1} \Phi_1^T f}_1 + \norm{1_{\cS_2}
\Phi_2^T f}_1}{\norm{\Phi_1^T f}_1 + \norm{\Phi_2^T f}_1}.\]
\end{definition}
In words, we consider the maximal fraction of total $\ell_1$ norm
which can be concentrated to the combined
index set $\cS_{1} \cup \cS_2$. Concepts of this kind go back
to \cite{DS89}.
Adequate control of joint concentration ensures
that the principle (\ref{PSep}) gives a successful
approximate separation.

\begin{proposition} \label{prop:mainestimate}
Suppose that $S$ can
be decomposed as $S=S_1^0+S_2^0$ so that
each component $S_i^0$ is relatively sparse in $\Phi_i$, $i=1,2$, i.e.,
\[
\norm{1_{\cS_1^c} \Phi_1^T S_1^0}_1 + \norm{1_{\cS_2^c} \Phi_2^T S_2^0}_1
\le \delta.
\]
Let $(S_1^\star,S_2^\star)$ solve (\ref{PSep}).
Then
\[
  \norm{S_1^\star-S_1^0}_2 + \norm{S_2^\star-S_2^0}_2
\le \frac{2\delta}{1-2\kappa}.
\]
\end{proposition}


\begin{definition}
Given tight frames
$\Phi=(\phi_i)_i$ and $\Psi=(\psi_j)_j$
and an index subset $\cS$ associated with
expansions in frame $\Phi$,
we define the \emph{cluster coherence}
\[
        \mu_c (\cS, \Phi; \Psi) = \max_{j}  \sum_{i \in \cS} | \langle \phi_{i}, \psi_{j} \rangle|.
\]
\end{definition}

In many studies of $\ell_1$ optimization,
one utilizes instead the mutual coherence
\begin{equation} \label{singletonCoh}
   \mu (\Phi, \Psi) = \max_{j}  \max_{i} | \langle \phi_{i}, \psi_{j} \rangle|,
\end{equation}
whose importance was shown by \cite{DH01}.
This may be called the {\it singleton coherence.}
In contrast, cluster coherence
bounds coherence between a single member of frame $\Psi$
and a cluster of members of frame $\Phi$, clustered at $\cS$.

A related notion called `cumulative coherence' was introduced in \cite{Tro04};
that notion maximizes over subsets $\cS$ of a given size, whereas here we fix a
specific set $\cS$ of coefficients.
In applying our concept,  the index subsets we will consider are  not abstract,
but will instead have a specific geometric interpretation, associated to proximity
to certain curves in phase space.
Maximizing over all subsets of a given size
 would give very loose bounds,
and would not be suitable for our purposes.
Several other coherence measures involving subsets  appear in the
literature, e.g., \cite{BTW07} and \cite{YL06}, but we do not see a strong relation to
cluster coherence.

\begin{lemma} \label{lemm:kappamuc}
We have
\[ \kappa(\cS_1,\cS_2) \le \max\{ \mu_c(\cS_1,\Phi_1;\Phi_2), \mu_c(\cS_2,\Phi_2;\Phi_1)\}.\]
\end{lemma}

The proofs for Proposition \ref{prop:mainestimate} and Lemma \ref{lemm:kappamuc} are presented
in Section \ref{subsec:proofs_1}.

\subsection{Intended Application}
\label{sec-IntendApplic}

The concepts  of this section will now be applied
 to  (\mbox{\sc CSep}),  at scale $j$ only.
With $f  = \cP + \cC$  our distribution of interest,
and $F_j$ our bandpass filter, and
we set $f_j  = F_j \star f$.
Throughout this section,
the  object $S \equiv f_j$ and
the tight frames are
 $\Phi_1$,  the full radial wavelet frame,
 and $\Phi_2$, the full curvelet tight frame.
 We apply the optimization problem (\mbox{\sc Sep}),
 getting subsignal components $S_1^\star$ and $S_2^\star$,
 which we then relabel as
 the wavelet component $W_j$ and curvelet component $C_j$;
 {\sl one should check that with this
 sequence of substitutions, problem (\mbox{\sc Sep}) of this section
 becomes  (\mbox{\sc CSep}) of the introduction}.

The key problem in the application of Proposition \ref{prop:mainestimate}
to the aforementioned setting is the correct choice of the clusters of significant
coefficients at each scale. If
those clusters are chosen `too small asymptotically', the relative sparsity
will blow up, and if chosen `too large asymptotically', we lose control
of the cluster coherence.  We define those clusters on the ideal decomposition,
where wavelets are used to analyze the point singularity and curvelets
are used to analyze the curve singularity.  (Such
clusters are theoretical, non-observable  entities.)

For a series of wavelet-coefficient thresholds $\eps_{j,1}$  to be specified,
the cluster of significant wavelet coefficients can be provisionally defined as
\begin{equation} \label{threshClust1}
  \cS_{1,j}  =  \{ \lambda \in \Lambda_j^{\pm 1} :  |\langle \psi_\lambda , \cP_j \rangle| > \eps_{j,1} \cdot  \| (\langle \psi_{\lambda'} , \cP_j \rangle)_{\lambda'} \|_{\ell_\infty(\Lambda_j^{\pm 1})} \}.
\end{equation}
For a  series of curvelet-coefficient thresholds $\eps_{j,2}$ to be specified,
the cluster of significant curvelet coefficients can be provisionally taken as
\begin{equation} \label{threshClust2}
  \cS_{2,j}  =  \{ \eta \in  \Delta_j^{\pm 1} :  |\langle \gamma_\eta, \cC_j \rangle| > \eps_{j,2} \cdot  \| (\langle \gamma_{\eta'}, \cC_j \rangle )_{\eta'}\|_{\ell_\infty(\Delta_j^{\pm 1})} \} .
\end{equation}
Each threshold choice picks a specific point on the tradeoff
between relative sparsity and cluster coherence.

Then $\delta_j$ will denote the degree of approximation by significant coefficients,
the sum $\delta_j = \delta_{j,1} + \delta_{j,2}$
of the wavelet approximation error to the point singularity:
\[
   \delta_{j,1} = \sum_{\lambda \in \cS_{1,j}^c} |\langle  \psi_\lambda , \cP_j \rangle| ;
\]
and the curvelet approximation error to the curvilinear singularity:
\[
   \delta_{j,2} = \sum_{\eta \in \cS_{2,j}^c} |\langle \gamma_\eta, \cC_j \rangle|  .
\]

Finally,  $\kappa(\cS_{1,j},\cS_{2,j})$, the degree of joint wavelet-curvelet concentration
at the significant subsets, will be controlled by two cluster coherences:
$ \mu_c(\cS_{1,j}  , \Phi_1 ; \Phi_2)$
the maximal coherence of a curvelet to a cluster of significant
wavelet coefficients; and
$ \mu_c(\cS_{2,j}  , \Phi_2 ; \Phi_1)$
the maximal coherence of a wavelet to a cluster of significant
curvelet coefficients.  We have

\begin{corollary} \label{coro:maincorollary}
Suppose that the sequence of transform-space
clusters $(\cS_{1,j})$, and $(\cS_{2,j})$
has {\bf both}  of the following two properties:
(i) asymptotically negligible cluster coherences:
\[
   \mu_c(\cS_{1,j}  , \Phi_1 ; \Phi_2),    \mu_c(\cS_{2,j}  , \Phi_2 ; \Phi_1) \goto 0, \qquad j \goto \infty,
\]
and (ii) asymptotically negligible cluster approximation errors:
\[
  \delta_j = \delta_{1,j} + \delta_{2,j} = o(\|f_j\|_2) , \qquad j \goto \infty.
\]
Then we have asymptotically near-perfect separation:
\[
    \frac{ \|W_j - \cP_j \|_2 +  \|C_j - \cC_j \|_2 }{\|f_j\|_2} \goto 0, \qquad j \goto \infty.
\]
\end{corollary}
The main result \---Theorem 1.1 \--- follows from this
lemma,  but this will require
sufficiently good estimates for cluster coherence
for clusters defined as sufficiently good approximations
to the objects of interest.

We remark that although the threshold in the {\it provisional} definition of the
clusters in (\ref{threshClust1})-(\ref{threshClust2}) provides a
means to balance between relative sparsity and cluster coherence,
there is no {\it a priori} guarantee  that there exists a threshold for which
conditions (i) and (ii) of Corollary \ref{coro:maincorollary} are true.
The main achievement of this paper is to show that this is indeed
possible with a revised definition.
In fact we do  not finally define the clusters by  (\ref{threshClust1})-(\ref{threshClust2});
note that the clusters must simply exhibit properties (i) and (ii) of the Corollary. We intend to
make use of our freedom of definition in order to get the needed properties.

\section{Microlocal Analysis Viewpoint}
\label{sec:microlocal}

The proof of the main result boils down to
defining clusters and bounding cluster coherences
and cluster approximation errors; our approach
is inspired by microlocal analysis.  In effect, we consider
the case where the cluster is either a string
of curvelet coefficients in the cone of influence
of a curvilinear singularity, or a block of wavelet
coefficients in the cone of influence of a point singularity,
and we must bound interactions between clusters in one
frame and elements in the other frame.  Microlocal analysis
provides a simple organizational framework that
immediately suggests which interactions 'ought' to be
small and large, based on
the geometry of the overlaps between phase portraits
in the underlying microlocal phase space.

\subsection{Microlocal Analysis Concepts}

The singular support of a distribution $f$, $\SgSp(f)$,
is the set of points where $f$ is not locally $C^\infty$.
In the geometric separation setting, we have
\[
  \SgSp(f) = \SgSp(\cP+\cC) = \SgSp(\cP) \cup \SgSp(\cC) = \{ x_i \}  \cup image(\tau)
\]
because we have constructed  the distributions
$\cP$ and $\cC$ so their singularities have this
form.

\begin{figure}
\centering
\includegraphics[height=1.75in]{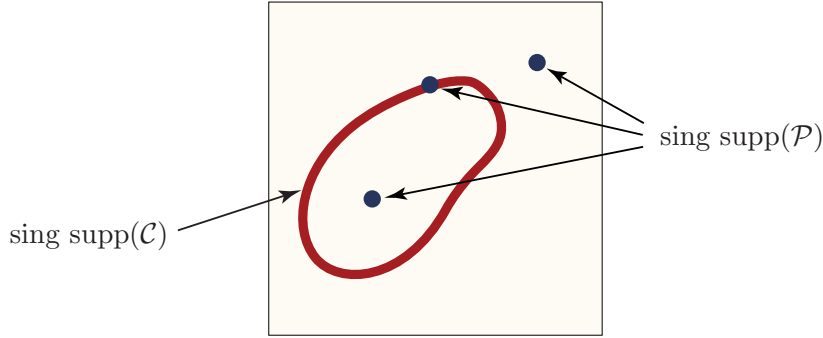}
\put(-240,35){$\SgSp(\cC)$}
\put(6,72){$\SgSp(\cP)$}
\caption{Singular supports of  the point singularity $\cP$ and
the curvilinear singularity $\cC$. The two supports overlap
in one point.}
\label{fig:PandC}
\end{figure}

Note that the points $x_i$ can intersect the image of the curve $\tau$ --
we make no separation hypothesis asking the point singularities
to `stay away'  from the curvilinear singularities. Figure \ref{fig:PandC} displays
the singular support of $f$ and the contributions
from $\cP$ and $\cC$.

To properly separate between pointlike and curvelike singularities
we need to consider a {\it phase space} for microlocal analysis
 indexed by  position-orientation pairs $(b,\theta)$;
such pairs can describe the locations and orientations
where $f$ has singular behavior. The orientational component $\theta$
will be regarded as an element in $\bP^1$, the real projective space
in $\bR^2$. (Here we identify $\bP^1$ with $[0,\pi)$ and freely
write one or the other in what follows.
It may at first seem more natural to think of  directions $[0,2\pi)$ rather than
orientations $[0,\pi)$, note however that in this paper
we consider {\it real-valued} distributions  $\cP +\cC$ measured by real-valued
curvelets $\gamma_\eta$ so directions are not resolvable, only orientations.
We also frequently abuse notation as follows: we will write $|\theta - \theta'|$
when what is actually meant is geodesic distance between two points on $\bP^1$.)

Living in this phase space
is the wavefront set $WF(f)$; roughly, this is the set of position-orientation
pairs at which $f$ is nonsmooth; for more
details, see: \cite{Hoe03,CD05a,KL07}.

Under the geometric separation model of Section 1, we
have
\[
      WF(\cP) = \Sp(\cP) \times \bP^1,
\]
since a point singularity is singular in all directions
on its singular support, and is singular nowhere else; while
\[
      WF(\cC) =  \{ (\tau(t), \theta(t)): t \in [0,L(\tau)] \}
\]
where $\tau(t)$ is a unit-speed parametrization of $\cC$
and $\theta(t)$ is the normal direction to $\cC$ at $\tau(t)$ regarded in $\bP^1$.

It is convenient to think of the parameter space
for microlocal analysis as a
plane $\bR^2$ of positions $b$ lying beneath a third dimension
of orientations $\theta$.  Then the wavefront set of a point singularity
$\cP$ concentrated on a single point $x_1$
is a vertical line segment $\{x_1\} \times \bP^1$, corresponding to
singular behavior in every direction at a given point, while
the wavefront set of $\cC$ is a more general curve in phase space.
Even if $x_1$ meets $image(\tau)$,
so the singular support of
the point singularity and a curvilinear singularity
overlap at $x_1$, they behave quite differently
as wavefront sets in the full 3D parameter space, which gives us hope for separation.

\begin{figure}
\centering
\includegraphics[height=2.3in]{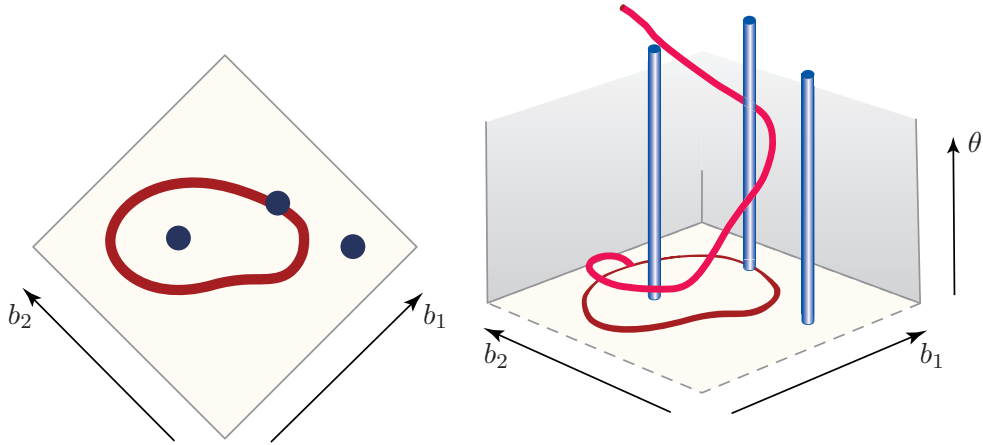}
\put(-360,45){$b_2$}
\put(-203,45){$b_1$}
\put(-180,30){$b_2$}
\put(-15,30){$b_1$}
\put(3,110){$\theta$}
\caption{Left panel: singular supports of  $\cP$ and $\cC$ (compare Figure \ref{fig:PandC}).
Right panel: wavefront sets of  $\cP$ and $\cC$  are indicated by blue and red tubes
in $(b_1,b_2,\theta)$ space.
The singular support  of $\cC$ is also indicated;
it is the $b_1$-$b_2$ projection of the wavefront set.
Due to the additional dimension,  the  wavefront sets
are more distinctly separated than  the corresponding singular supports.}
\label{fig:PandCinPhaseSpace}
\end{figure}

Figure \ref{fig:PandCinPhaseSpace} illustrates the 3D phase space with
pointlike and curvelike singularities superposed.

\newcommand{\cD}{{\cal D}}
\subsection{Support of Frame Elements }

A radial wavelet $\psi_\lambda$, $\lambda = (a_0,b_0)$
is  `morally' supported  in $x$-space  in a spatial ball $B(a_0,b_0)$,
defined by
\[
  B(a_0,b_0) = \{ x:  | x - b_0 |/a \leq 1 \} ;
\]
this statement should not however be taken too literally, as the
radial wavelets we study are not
of compact support. The more precise statement is that
the wavelet decays rapidly in the variable  $ | x - b_0 |/a $.
The next statement uses the notation
\[
   \langle a \rangle = (1 + a^2)^{1/2} .
\]

\begin{lemma}
 For each $N = 1,2, \dots$
there is a constant $c_N$ so that
\[
   |\psi_{a,b}(x)| \leq c_N \cdot  a^{-1} \cdot  \langle | x - b |/a \rangle^{-N}, \qquad \forall a \in \bR^{+} \;
   \forall   b,x \in \bR^2.
\]
\end{lemma}

An individual curvelet $\gamma_\eta$,
$\eta = (a_0,b_0,\theta_0)$ is  `morally' supported in
$x$-space inside an anisotropic spatial ellipse.
To make this precise, let $D_{1/a}$ be the diagonal
matrix $diag(1/a,1/\sqrt{a})$ and $R_{-\theta}$ denote
planar rotation by $-\theta$ radians.
Let
\[
  P_{a,\theta} = D_{1/a} R_{-\theta}
\]
denote the parabolic directional dilation operator
which dilates much more strongly in the $\theta$ direction
than in the orthogonal direction. For a vector $v \in \bR^2$
define the norm
\[
   |v|_{a,\theta} = | P_{a,\theta}(v)| ;
\]
the unit ball in this norm is ellipsoidal, with minor axis pointing in direction
$\theta$.
A curvelet is morally supported in the ellipse
 $E(a_0,b_0,\theta_0)$
defined by
\[
   E(a_0,b_0,\theta_0) = \{ x:  | x - b_0 |_{a_0,\theta_0}  \leq 1 \} ,
\]
Again the correct formal statement is that there
is rapid decay in the variable $| x - b |_{a,\theta} $.
The following is proved in \cite[Lemma 1, page 168]{CD05a}.

\begin{lemma}
\label{lemma:estimate_gamma}
For each $N = 1,2, \dots$
there is a constant $c_N$ so that
\begin{equation} \label{curveletMod}
   |\gamma_{a,b,\theta}(x)| \leq c_N \cdot  a^{-3/4} \cdot  \langle | x - b |_{a,\theta} \rangle^{-N}, \qquad \forall a \in \bR^{+} \;  \forall \theta \in [0,\pi) \;
   \forall   b,x \in \bR^2.
\end{equation}
\end{lemma}

\begin{figure}
\centering
\includegraphics[height=1.75in]{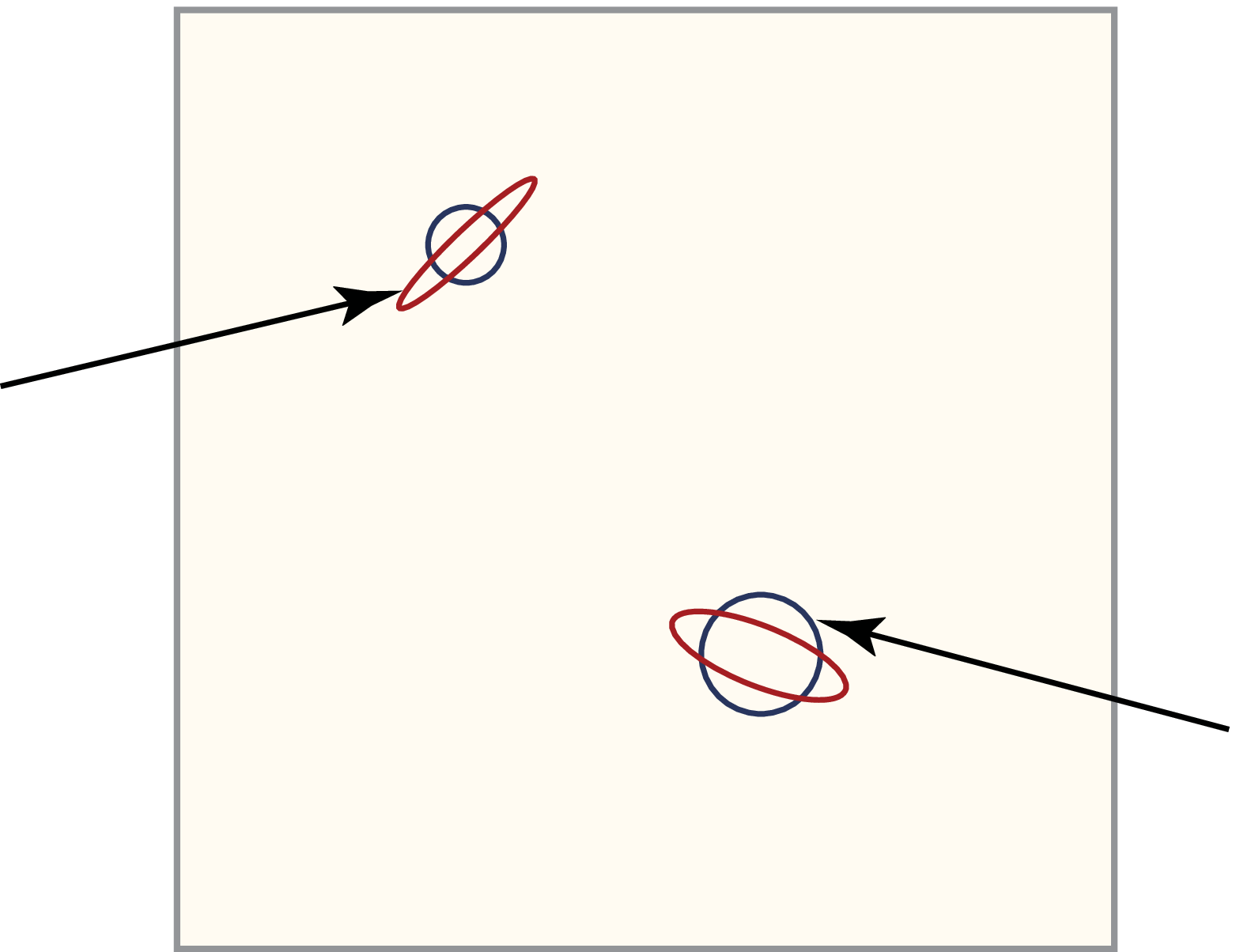}
\put(5,25){$\Sp(\psi_{a',b'})$}
\put(-222,72){$\Sp(\gamma_{a,b,\theta})$}
\caption{Effective supports
of wavelets and curvelets at different scales. The shapes of those supports
become increasingly distinct with increasingly fine scales.}
\label{fig:WaveletCurvelet}
\end{figure}

In Figure \ref{fig:WaveletCurvelet}, we visualize  these support relationships.

\subsection{Phase-Space Support of Frame Elements}

We also study the location-orientation
behavior of frame elements, i.e., the attribution
of regions in an orientation/location domain as
regions of significant activity in a distribution $f$.
The wavefront set gives a qualitative way to do this;
we use the continuous curvelet transform (CCT)
to do so quantitatively. This transform is defined by
\[
  \Gamma_f(a,b,\theta) = \langle \gamma_{a,b,\theta} , f \rangle
\]
and is indexed by triples $(a,b,\theta)$, where
$a > 0$, $b \in \bR^2$, and $\theta \in [0,\pi)$; this
associates a function $f$ to a scale/location/direction
domain.
There is a natural measure on this domain:
\[
    m(da,db, d\theta ) = a^{-3} \cdot da \cdot db \cdot  d\theta .
\]
Consider a high-pass function $f$, i.e.,
one whose Fourier transform of $f$ vanishes near the origin;  $\hat{f}(\xi) = 0$, $|\xi| \leq \pi/a_0$,
see \cite{CD05a}.
The CCT  offers  a Parseval relationship for high-pass functions:
\[
    \int\int\int |\Gamma_f(a,b,\theta)|^2 m(da,db, d\theta ) = \| f \|_2^2,
\]
see \cite{CD05b}.
Hence  the energy in $f$ is distributed through
the scale-location-direction domain by the curvelet
transform offering a  portrait of the function's
significant activity.

Consider the transform of a radial wavelet:
\[
  \Gamma_{\psi_{a_0,b_0}}(a,b,\theta) = \langle \gamma_{a,b,\theta} , \psi_{a_0,b_0} \rangle , \qquad a>0, \; b \in \bR^2,   \; \theta \in [0,\pi) .
\]
Because $\psi_{a_0,b_0}$ is radial, $\Gamma_{\psi_{a_0,b_0}}(a,b,\theta)$
is constant, independent of $\theta$, and decays rapidly
in variables $|\log_2(a/a_0)|$ and $|b - b_0|/a_0$.
It is morally localized to a cell of the form
\[
   W(a_0,b_0) =  \{ (a,b,\theta):  |\log_2(a/a_0)| < 1, \theta  \in [0,\pi),
   |b - b_0|/a_0 \leq 1 \};
\]
we have the following  formal statement, proved in Subsection \ref{pf-Lemm-3.3} below:

\begin{lemma} \label{lemm:gamma_psi_estimate} For each $N=1,2,\dots$,
there is  a constant $c_N$ so that
\[
  |\langle \gamma_{a,b,\theta} , \psi_{a_0,b_0} \rangle |
  \leq c_N \cdot a^{1/4} \cdot 1_{\{|\log_2(a/a_0)| < 3\}}  \cdot
    \langle |b - b_0|_{a,\theta} \rangle^{-N} .
\]
\end{lemma}

It implies the following for a scale-conditional phase portrait
of a wavelet (i.e., we freeze the analysis scale $a$
at a specific value, and inspect $\Gamma_{\psi_{a_0,b_0}}(a,b,\theta) $
as a function of variables $b$ and $\theta$): when freezing $a =a_0$,
we see that $\Gamma$ is `morally' supported in a vertical tube above the
point $b_0$; each horizontal cross-section is a ball of width $a_0$, i.e., $B(a_0,b_0)$.

 Consider now the transform of a curvelet:
\[
  \Gamma_{\gamma_{a_0,b_0,\theta_0}}(a,b,\theta) = \langle \gamma_{a,b,\theta} , \gamma_{a_0,b_0,\theta_0} \rangle , \qquad a>0, \; b \in \bR^2,   \; \theta \in [0,\pi) .
\]
This is `morally' localized to a cell
of the following form:
\[
   Q(a_0,b_0,\theta_0) =  \{ (a,b,\theta):  |\log_2(a/a_0)| < 1, |\theta - \theta_0| < \sqrt{a_0},
   |b - b_0|_{a_0,\theta_0} \leq 1 \};
\]
the correct formal statement being [see (21) in \cite{CD05b}]
\begin{lemma} \label{lemm:gamma_gamma_estimate}
For each $N=1,2,\dots$,
there is a constant $c_N$ so that
\[
  |\langle \gamma_{a,b,\theta} , \gamma_{a_0,b_0,\theta_0} \rangle |
  \leq c_N \cdot 1_{\{|\log_2(a/a_0)| < 3\}} \cdot 1_{\{ |\theta - \theta_0| < 10 \sqrt{a_0}  \}} \cdot
    \langle |b - b_0|_{a_0,\theta_0} \rangle^{-N} .
\]
\end{lemma}
(We remind the reader of our convention that, for two points $\theta,\theta^\prime \in [0,\pi)$,
$|\theta-\theta'|$ really means geodesic distance in $\bP^1$.)

Thus each curvelet is supported in scale, location, direction
in a set which effectively has a product structure, and is compactly supported in
both scale and orientation. In a scale-conditional
phase portrait, freezing  $a = a_0$, we see a vaguely ellipsoidal structure
with slice $\theta = \theta_0$ exhibiting an
anisotropic footprint, like $E(a_0,b_0,\theta_0)$.

\begin{figure}
\centering
\includegraphics[height=2.3in]{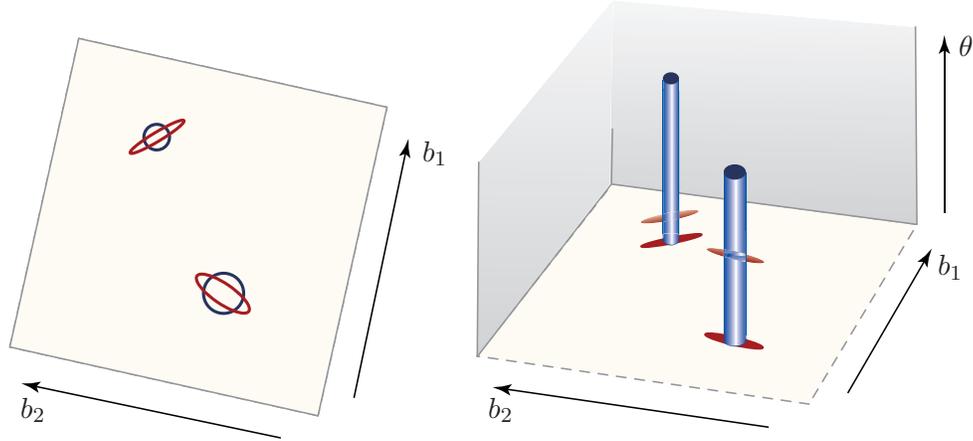}
\put(-352,8){$b_2$}
\put(-200,105){$b_1$}
\put(-175,8){$b_2$}
\put(-5,62){$b_1$}
\put(3,145){$\theta$}
\caption{Left panel: effective supports of wavelets (blue) and curvelets (red);
compare Figure \ref{fig:WaveletCurvelet}.
Right Panel:  Phase space portrait -- obtained by CCT -- of the same wavelets and curvelets.
Effective supports of curvelets are indicated in the $(b_1,b_2)$ plane.
The portrait conveys the intuition that while the
 effective support of wavelets and curvelets may overlap substantially in the spatial domain,
 in the full phase space, the relative overlap is significantly reduced.}
\label{fig:WaveletCurveletPhaseSpace}
\end{figure}

In Figure \ref{fig:WaveletCurveletPhaseSpace} we visualize  the scale-conditional portraits
of a wavelet and a curvelet.
The intuition to be fostered from these figures is that
curvelets don't interact very heavily with
wavelets, because they have such different
support in $(b,\theta)$, even when they have the
same scale and location parameters, $a_0$ and $b_0$.

This support structure  of curvelets is notable
for other reasons.
The support structure is implicit in the natural measure:
\begin{eqnarray*}
    m(da,db, d\theta ) &=& a^{-3} \cdot da \cdot db \cdot  d\theta \\
                       &=& \quad \frac{da}{a} \cdot \frac{db}{a^{3/2}} \cdot \frac{d\theta}{a^{1/2}},
                        \qquad 0 < a < a_0, b \in \bR^2 , \; \theta \in \bP^1 .
\end{eqnarray*}
The second expression has the following interpretation.
We assign roughly unit measure to curvelet cells
\[
    m(Q(a,b,\theta)) \approx 1 , \qquad 0 < a < a_0, \; b \in \bR^2, \theta \in \bP^1 ;
\]
 there are about $a^{-3/2}$ locations per unit volume
 at a fixed scale and orientation
and about  $a^{-1/2}$ orientations at a fixed scale and location.

The support structure is mimicked by the discretization
of the curvelet tight frame, which morally
breaks the $(a,b,\theta)$ domain into disjoint cells
$Q(a_j,b_{j,k,\ell},\theta_{j,\ell})$ and samples
the continuous transform once per cell; for details
see \cite{CD05b}.


%

\subsection{CCT  of Singularities}

The CCT can be used to analyze the singularities
$\cP$ and $\cC$ directly.  Let's impose regularity conditions
on $\cC$.


Let the Hausdorff pseudo-distance $d$ between $b$ and
a curve $\tau$ be defined by
\[
d(a,b,\theta;\tau) = \min \{  |x - b|_{a,\theta} : x \in image(\tau) \}.
\]

\begin{definition}
A finite-length planar curve $\tau$ will be called {\em parabolically regular},
if, for $N=1,2,\dots$,  there is a constant $c_N$
so that  for $a \in (0,1)$ and all $b,\theta$,
\begin{equation} \label{paraReg}
\int \langle |\tau(t) - b|_{a,\theta} \rangle^{-N} dt    \leq c_N \cdot  a^{1/2}  \cdot \langle d(a,b,\theta;\tau) \rangle ^{-N}.
\end{equation}
\end{definition}

Traditional nice curves, such as line segments, circles, etc. are parabolically regular,
though we skip the demonstration.

\begin{lemma} \label{lem-tube}
Let the singularity $\cC$ be defined as in (\ref{curvedef})
by a parabolically regular curve $\tau$.
Then, for each $N = 1,2,\dots$,
\[
  |\Gamma_{\cC}(a,b,\theta)| \leq c_N \cdot a^{-1/4} \cdot \langle d((a,b,\theta); \tau) \rangle^{-N} .
\]
\end{lemma}

In words, the CCT may be large at phase space points close
to $\tau$, but elsewhere it is very small.

\noindent
{\bf Proof.} The definition of $\cC$ as a linear functional on the Schwartz space gives
\[
     \Gamma_{\cC}(a,b,\theta) = \int \gamma_{a,b,\theta}(\tau(t)) dt.
\]
The estimates (\ref{curveletMod}) and (\ref{paraReg}) give
\begin{eqnarray*}
     |\Gamma_{\cC}(a,b,\theta)| & \leq& \int c_N \cdot a^{-3/4} \cdot \langle |\tau(t) - b|_{a,\theta} \rangle^{-N} dt   \\
   & \leq&  c_N \cdot a^{-3/4}  \cdot C'_N \cdot  a^{1/2}  \cdot \langle d(a,b,\theta;\tau) \rangle ^{-N} \\[1ex]
   &=& C''_N \cdot  a^{-1/4}  \cdot \langle d(a,b,\theta;\tau) \rangle ^{-N} .  \mbox{\qed}
\end{eqnarray*}

\subsection{Heuristics}

We are now in a position to give a heuristic explanation why the strategy
announced in Section \ref{sec-IntendApplic} is likely to work. In effect,
there are very instructive analogies between the calculations
needed to implement that strategy
and the behavior of certain `tubes' in phase space.
The reader will have noticed that the curvelet parameter
 $(a,b,\theta)$ and the wavefront set parameter
$(b,\theta)$ differ only by the latter's provision of a scale.
Hence there is some analogy between the scale-conditional
portrait by CCT and the wavefront set --
both provide measures of the activity of an object,
indexed by location and orientation.

In effect, the scale-conditional portrait by CCT is a ``thickened-out''
version of the wavefront set.
A point singularity has a wavefront
set which is a vertical line in phase space, and a scale-conditional
portrait which is localized near a thin vertical tube.
A curvilinear singularity
has a wavefront set which is a curve in phase space,
while Lemma \ref{lem-tube} says that, morally,
the curvelet transform of the object $\cC$ `lives'
near a tube.
The tube in question has thickness $\approx a$.
Look at the scale-conditional
portrait, and define the  tube
\[
   T(a_0) = \bigcup_{b,\theta}   (E(a_0,b,\theta)  \times \{\theta\}),
\]
where the union is over $b$, $\theta$ satisfying
\[
    d((a_0,b,\theta); \tau) \leq 1.
\]
As $a_0 \goto 0$, this tube shrinks
down to a curve $\tilde{\tau}$ in phase space $\bR^2 \times \bP^1$
defined by
\[
     \tilde{\tau}(t) = (\tau(t),\theta(t)),
\]
where $\theta(t)$ is the orientation of the normal to $\tau(t)$.
In short, in the sense of set convergence
\begin{equation*}
   T(a_0) \goto WF(\cC), \qquad a_0 \goto 0.
\end{equation*}
Thus, the wavefront set and the curvelet transform
both signal that the activity in location-orientation space
is concentrated near $image(\tilde{\tau})$.

More is true.  A wavelet has a scale-conditional
portrait which is a thin vertical tube -- similar to the
phase portrait of a point singularity --
while a curvelet has a scale-conditional portrait which
is a tube surrounding a little `piece of a curve' in phase space, i.e., it morally
has a position and orientation.

This visual analogy suggests
that curvelets are incoherent to wavelets -- because of the low overlap
in phase space. Indeed, from Parseval,
\[
    \int\int\int \Gamma^\star_{\psi_{a_0,b_0}}(a,b,\theta)   \Gamma_{\gamma_{a_1,b_1,\theta_1}}(a,b,\theta)  m(da,db, d\theta ) =  \langle \psi_{a_0,b_0}, \gamma_{a_1,b_1,\theta_1} \rangle,
\]
and so the low overlap between the two phase portraits indeed
will cause relatively low singleton coherence (\ref{singletonCoh});
indeed the tubelet associated to a given curvelet and the tube
associated to a given wavelet visibly have relatively small overlap in
the scale-conditional phase portrait; for example, if we compare
the overlap of effective supports in phase space to the overlap of effective supports
in the spatial domain, we see that the fractional overlap is dramatically smaller at fine scales
in the phase space portrait than it is in the spatial domain portrait.

\begin{center}
\begin{figure}[h]
\hspace*{2cm}
\includegraphics[height=2.2in]{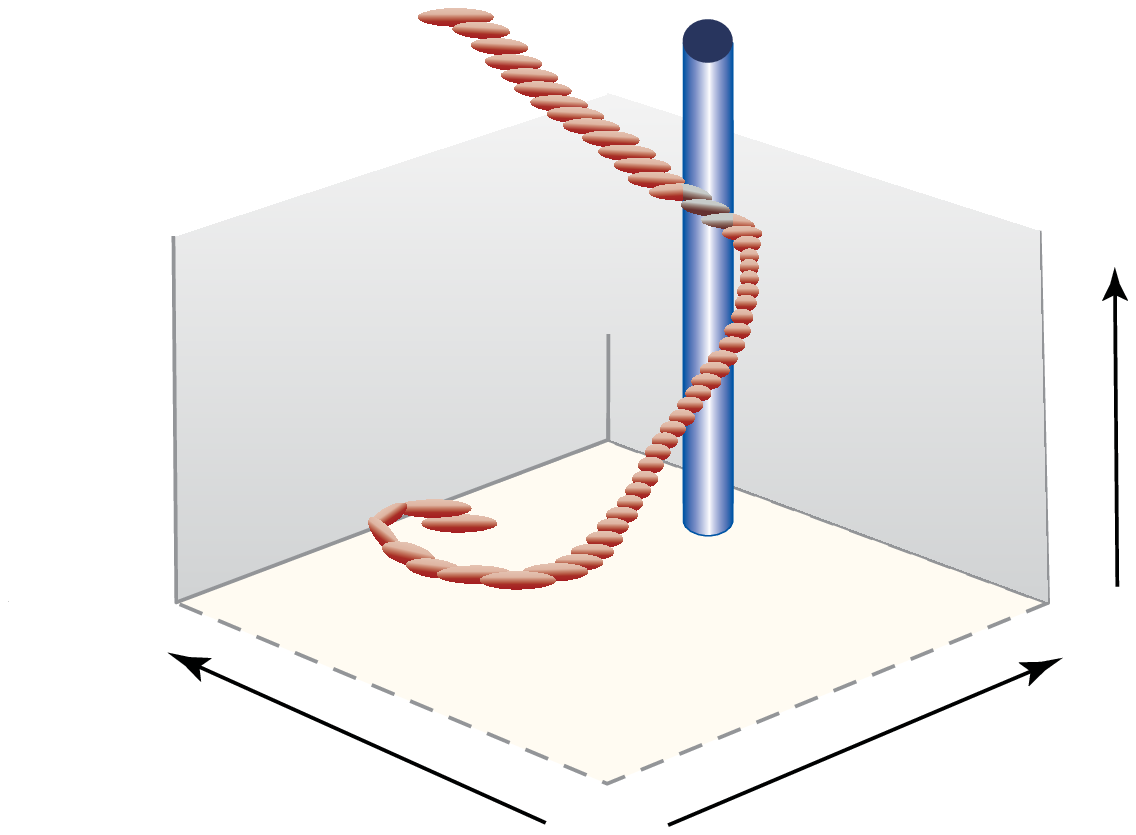}
\put(-180,30){$b_2$}
\put(-15,30){$b_1$}
\put(3,110){$\theta$}
\caption{Phase space portrait of a cluster of curvelets and one single wavelet.
Note the visually small overlap of these two geometrical objects. This suggests by analogy
that each single wavelet is not coherent with anything built from such a cluster of
curvelets.}
\label{fig:baseballs}
\end{figure}
\end{center}

However, the singleton coherence is not sufficiently small to
be powerful in the present setting.
Instead, this paper develops cluster coherence. The visual
analogy presented in Figure  \ref{fig:baseballs} suggests how to bound the cluster coherence
and suggests that the proof
strategy of Section \ref{sec-IntendApplic} will succeed.
To understand that analogy, let's study Figure \ref{fig:baseballs}.
If we let $\cS_1$ denote the set of significant
 wavelet coefficients in the radial wavelet transform
 of $\cP$ at scale $a_0$, and $\cS_2$ denote the set of significant
 curvelet coefficients in the curvelet transform
 of $\cC$ at scale $a_0$,
 we believe the reader will be easily able
 to motivate the following assertions on the basis
 of Figures \ref{fig:WaveletCurvelet}-\ref{fig:baseballs}:
 \bitem
 \item Wavelets in $\cS_1$ are associated to vertical tubes clustering
 around the point singularities in $\cP$;
 \item Curvelets in $\cS_2$ are associated with tubes
 clustering around the curvilinear phase portrait
 of $\cC$;
 \item No single wavelet's phase portrait overlaps much
 with the cluster of curvelet phase portraits of $\cS_2$;
 \item No single curvelet phase portrait overlaps with the
cluster of wavelets in $\cS_1$.
 \eitem

Let's outline a pseudo-calculation inspired by these
visual observations.  First, we consider a pseudo-calculation
of the cluster coherence,
\[
\mu_c( \cS_{1,j} , \mbox{\sc Wavelet scale $j$} ; \mbox{\sc Curvelet scale $j$})
= \sup_\eta \sum_{\lambda \in \cS_{1,j}} |\langle \gamma_\eta , \psi_\lambda \rangle |    .
\]
With $(\psi_i)_i$ an enumeration of the significant wavelets in the expansion at scale $a_j = 2^{-j}$,
\begin{eqnarray*}
\sum_i |\langle \gamma_\eta , \psi_i \rangle |  &=& \sum_i   \left| \int\int\int  \Gamma_{\gamma_\eta}^* \Gamma_{\psi_i} dm \right|  \\
   &\leq& \sum_i \int\int\int  |\Gamma_{\gamma_\eta}|(a,b,\theta)  |\Gamma_{\psi_i}|(a,b,\theta) dm \\
   &\leq&\int\int\int  |\Gamma_{\gamma_\eta}|(a,b,\theta)    \sum_i  |\Gamma_{\psi_i}|(a,b,\theta) dm  .
\end{eqnarray*}
Now use the bounds on $ |\Gamma_{\psi_i}|(a,b,\theta)$ given above in Lemma \ref{lemm:gamma_psi_estimate},
and deploy the slogan that {\it only a  bounded number of significant wavelets at any
given scale interact strongly  with any specific phase space point}; we have
\[
    \sum_i   |\Gamma_{\psi_i}|(a,b,\theta) \lesssim C_1 \cdot a_j^{1/4} ,
\]
where the sum is over the significant wavelets
at scale $a_j$ and $C_1$ does not depend on $j$,
and the symbol $\lesssim$ indicates an inequality motivated heuristically -- in this case
by the preceding italicized slogan.
We also observe that Lemma \ref{lemm:gamma_gamma_estimate} implies that the integral
over phase space of a curvelet phase portrait obeys
$\int\int\int  |\Gamma_{\gamma_\eta}|(a,b,\theta)  dm  \leq C_2$.
Combining this with the previous displays, we pseudo-conclude that
\[
  \mu_c( \cS_1 , \mbox{\sc Wavelet scale $j$} ; \mbox{\sc Curvelet scale $j$}) \goto 0,  \qquad j \goto \infty.
\]

Next, we consider a pseudo-calculation
of the cluster coherence
\[
\mu_c( \cS_2 , \mbox{\sc Curvelet scale $j$} ; \mbox{\sc Wavelet scale $j$}) = \sup_\lambda \sum_{\eta \in \cS_{2,j}} |\langle \gamma_i , \psi_\lambda \rangle |  .
\]
With $(\gamma_i)_i$ an enumeration of the significant curvelets in the expansion at scale $a_j = 2^{-j}$,
then for a fixed wavelet index $\lambda$ we have
\begin{eqnarray*}
 \sum_i |\langle \gamma_i , \psi_\lambda \rangle |  &=& \sum_i   \left| \int\int\int  \Gamma_{\gamma_i}^* \Gamma_{\psi_\lambda} dm \right|  \\
   &\leq& \sum_i \int\int\int  |\Gamma_{\gamma_i}|(a,b,\theta)  |\Gamma_{\psi_\lambda}|(a,b,\theta) dm \\
   &\leq& \int\int\int   |\Gamma_{\psi_\lambda}|(a,b,\theta)  \sum_i  |\Gamma_{\gamma_i}|(a,b,\theta)  dm.
\end{eqnarray*}
Now use the bounds on $ |\Gamma_{\gamma_i}|(a,b,\theta)$ given above in Lemma \ref{lemm:gamma_gamma_estimate},
and apply the slogan that {\it only a  bounded number of significant curvelets at any
given scale interact strongly with a given phase space point}.
We have
\[
    \sum_i   |\Gamma_{\gamma_i}|(a,b,\theta) \lesssim C_1 ,
\]
where the sum is over the significant curvelets
at scale $a_j$ and $C_1$ is a constant.
We also observe that Lemma \ref{lemm:gamma_psi_estimate} implies that the integral
over phase space of a wavelet phase portrait obeys
$\int\int\int  |\Gamma_{\psi_\lambda}|(a,b,\theta)  dm \leq C_2 \cdot a_0^{1/4}$,
where $a_0 = a_0(\lambda)$.
Combining this with the previous displays, we pseudo-conclude that
\[
  \mu_c( \cS_2 , \mbox{\sc Curvelet scale $j$} ; \mbox{\sc Wavelet scale $j$}) \goto 0,  \qquad j \goto \infty.
\]
We now turn to the approximation tasks posed by the strategy
announced in Section \ref{sec-IntendApplic}.
To pseudo-bound $\delta_{j,1}$, we fix $\eps > 0$
and define the tube in phase space
$\cT_{j,1}$ consisting of all scale/location pairs $(a,b)$ where the bound
provided in  Lemma \ref{lemma:pointwavelet} permits  coefficients larger
than $a_j^{1-\eps}$. Also, we let $i_j$ denote the index in the wavelet enumeration
beyond which such potentially significant coefficients
can no longer arise. We can heuristically approximate a sum of wavelet coefficients with
an integral over the phase space region covered by the union
of their phase space supports; then we have
\begin{eqnarray*}
    \delta_{j,1} &=& \sum_{i > i_j} | \langle \psi_i , \cP_j \rangle |  \nonumber \\
                        &\lesssim& c_N \cdot a_j^{-1/2} \cdot \int\int\int_{\cT_{j,1}^c} \langle b/a_j\rangle^{-N}  dm.
\end{eqnarray*}
For sufficiently large $N$, the integrand has powerful decay;  for large $j$, the width of the tube  $\cT_{j,1}$ is
significantly wider than the decay scale $a_j$; and so the integral becomes negligible
for large $j$, i.e., we pseudo-conclude  that $ \delta_{j,1}  \goto 0$.
To pseudo-bound $\delta_{j,2}$ we fix $\eps > 0$
and define the tube in phase space
$\cT_{j,2}$ consisting of all triples $(a,b,\theta)$ where the bound
provided in  Lemma \ref{lem-tube} permits  coefficients larger
than $a_j^{1-\eps}$. Also we let $i_j$ denote the index in the curvelet enumeration
beyond which such potentially significant coefficients
can no longer arise; again heuristically identifying a sum with a phase-space integral we have
\begin{eqnarray*}
    \delta_{j,2} &=& \sum_{i > i_j} | \langle \gamma_i , \cC_j \rangle | \nonumber \\
                        &\lesssim& c_N \cdot \int\int\int_{\cT_{j,2}^c} a^{-1/4} \cdot \langle d(a_j,b,\theta; \tau) \rangle^{-N}  dm .
\end{eqnarray*}
For large $N$ the integrand decays strongly; for large $j$ the width of the tube is
significantly wider than the decay scale $a_j$; and so again the integral becomes negligible
for large $j$, i.e., we pseudo-conclude  that $ \delta_{j,2}  \goto 0$.

In short, phase space diagrams, some elementary estimates motivated by tube overlaps,
and some cardinality `slogans'
combine to show plausibility of the strategy announced in Section \ref{sec-IntendApplic}.
In the sections to come, we rigorously carry out that strategy. While the details
are much more delicate than this plausibility argument would suggest,
the architecture of our full demonstration remains faithful to the
geometric viewpoint.

\section{The Cluster $\cS_{1,j}$ and its estimates}
\label{sec-Sec1j}

In this section, we define the cluster of wavelet coefficients $\cS_{1,j}$
of the filtered point singularity $F_j \star \cP$,
and estimate relative sparsity as well as cluster coherence using this cluster.
We intend to show that with this definition of cluster set,
\begin{equation} \label{wavClust1}
   \delta_{1,j} = o(\|f_j \|_2 ) , \qquad j \goto \infty ,
\end{equation}
and
\begin{equation} \label{wavClust2}
   \mu_c(\cS_{1,j}, \{ \psi_\lambda\};\{ \gamma_\eta \} ) \goto 0 , \qquad j \goto \infty.
\end{equation}
As explained in Section \ref{sec-IntendApplic},
this gives the needed part of Theorem 1.1
having to do with $\cS_{1,j}$.

WLOG we can assume that
\[
\cP=|x|^{-3/2} .
\]
The result for the more general $\cP$ of \eqref{pointdef} follows
easily by combining translation invariance with finitely many uses of the
triangle inequality.
We also, from now on, fix some
\[
\eps \in (0,1/32).
\]

Our first lemma is used frequently in what follows,  and is crucial
for our definition of the cluster of wavelet coefficients.  For the proof, see
Subsection \ref{subsec:proofs_3_1}.

\begin{lemma} \label{lemma:pointwavelet}
For each $N = 1,2, \dots$, there is a constant $c_N$ so that
\[
   \absip{\psi_{a_{j'},b}}{\cP_j} \le c_N \cdot 2^{j/2} \cdot 1_{\{|j-j'| \le 1\}} \cdot \langle |b/a_j| \rangle^{-N},
   \qquad \forall j, j' \in \bZ_+, \; \forall b \in \bR^{2}.
\]
\end{lemma}

In line with the heuristics of the previous section,
we think of our estimate as describing relative overlaps of
tubes in phase space. However, in the
particular case of wavelets, there is no directional selectivity,
so all that matters is the projection of phase space onto the spatial domain.
We measure spatial distances with
\[
d_2(x,A) = \min_{a \in A}\norm{x-a}_2, \quad x \in \bR^2,\,A \subseteq \bR^2,
\]
the Euclidean distance between a point $x$ and a set $A$.

Of course, since we are dealing with frames,
ultimately we have to consider discrete indices.
To support geometric intuition,  most of our arguments will be in the continuum
setting, restrictions to discrete sampling grids being delayed as
late in each argument as possible.

Morally, the points in phase space associated with
significant wavelet coefficients are contained
in a tube around $WF(\cP)$ in phase space.
This neighborhood of $WF(\cP)$ can be explicitly defined by
\[
\cN_{1}^{PS}(a) = \{b \in \bR^2 : d_2(b,\{0\}) \le D_1(a)\} \times [0,\pi),
\]
where
\[
D_1(a) = a^{(1-\eps)}.
\]
The shape of the tube reflects the isotropic behavior of $WF(\cP)$.
For an illustration of $\cN_{1}^{PS}(a)$, we refer to Figure \ref{fig:tubes}.

We define the cluster of wavelet coefficients around the point-singularity by intersecting
the tube $\cN_{1}^{PS}(a)$ with the wavelet lattice, i.e.,
\[
  \cS_{1,j} = \{ (j,k)\in\Delta_j^{\pm 1} : \{b_{j,k}\} \times [0,\pi) \in \cN_{1}^{PS}(a_j)\} .
\]
The remainder of the section establishes (\ref{wavClust1})-(\ref{wavClust2}) for this
definition of cluster.

\subsection{Size of $f_j$}

 \begin{lemma}
 \label{lemm:estimateforfj}  For some $c > 0$,
\[
\| f_j\|_2 \geq  c 2^{j/2} , \qquad j \goto \infty.
\]
\end{lemma}
{\bf Proof.}
Apply (\ref{energyMatch}) and (\ref{superpose}):
\begin{eqnarray*}
   \| f_j \|_2^2 &=&  \int_{\bR^2}  W^2(|\xi|/2^j) |\hat{f}(\xi)|^2 d\xi
         \geq  c \cdot \int_{\cA_j} |\xi|^{-1} d\xi
               \geq  c \cdot 2^j. \qquad \mbox{\qed}
\end{eqnarray*}

\subsection{$\cS_{1,j}$ offers low approximation error}

%


Now we are ready to state and prove the approximation error of $\cS_{1,j}$.

 \begin{lemma} \label{estimatefordelta1j}
 \[
\delta_{1,j} = \sum_{\lambda \in \cS_{1,j}^c} | \langle
\psi_\lambda, \cP_j \rangle | = o(\|f_j \|_2 ), \qquad  j \goto  \infty .
\]
\end{lemma}
{\bf Proof.}
Due to the specific filtering we use,
\[
\delta_{1,j} = \sum_{\lambda \in \cS_{1,j}^c} | \langle
\psi_\lambda, \cP_j \rangle | = \sum_{j'=j-1}^{j+1} \sum_{\{k :
\{b_{j',k}\} \times [0,\pi) \not\in \cN_1^{PS}(a_{j'})\}} | \langle
\psi_{(j',k)}, \cP_j \rangle |.
\]
Applying Lemma \ref{lemma:pointwavelet},  and picking $N$ so large that
$(N-1)\eps > 1/2$
\[
\delta_{1,j}
\le c \cdot \sum_{|k|  > 2^{j\eps} } c_N \cdot 2^{j/2} \cdot  \langle |k| \rangle^{-N}
= o(1), \qquad  j > j_0 .
\]
So the lemma is proved. \qed


\subsection{$\cS_{1,j}$ offers low cluster coherence}

\begin{lemma} \label{decayofmuc1}
\[
\mu_{c}(\cS_{1,j}, \{ \psi_\lambda\} ; \{ \gamma_\eta\}  ) \goto 0, \qquad  j \goto \infty .
\]
\end{lemma}

{\bf Proof.}
By Lemma \ref{lemm:gamma_psi_estimate}, if $\lambda  \in \Lambda_j^{\pm 1}$,
\[
     |\langle \psi_\lambda,\gamma_\eta \rangle|   \leq c \cdot 2^{-j/4}, \qquad \forall \eta.
\]
The definition of $\cS_{1,j}$ via $\cN_{1}^{PS}(a_j)$ implies that
\[
    \#(\cS_{1,j}) \leq c \cdot 2^{2j \eps}, \qquad \mbox{for sufficiently large }j.
\]
We conclude from $\eps < 1/8$ that
\[
    \sup_\eta \sum_{\lambda \in \cS_{1,j}} |\langle \psi_\lambda,\gamma_\eta \rangle |  \leq  c
\cdot 2^{-j(1/4 - 2\eps) }  \goto 0, \qquad j \goto \infty. \qquad \mbox{\qed}
\]

\section{Sparse Expansion of a Linear Singularity}
\label{sec-LinearSing}

After the relative ease with which we obtained
concentration estimates (\ref{wavClust1})-(\ref{wavClust2})
for the cluster of significant wavelet coefficients,  we must now
brace ourselves for the considerably harder challenge posed by the analogous estimates
for the cluster of significant curvelet coefficients. This extra work seems, at least
to us, much more rewarding, as it involves a full-blown use of phase space geometry.

In this section, we develop essential infrastructure for  the analysis to come,
documenting the sparsity of curvelet coefficients of a special linear singularity.
Let $w_2: \bR \mapsto [0,1]$ be a smooth function to be specified later (cf. Subsection \ref{subsec:cutting}),
supported in $[-1,1]$, and define the very special distribution $w\cL$ supported on a line segment
$\{0\} \times [-\rho,\rho]$ by
\[
w\cL = w_2(x_2/\rho) \cdot \delta_0(x_1).
\]
Then we can write
\[
\widehat{w\cL} = \hat{w} \star \hat{\cL},
\]
where
\[
\hat{w} = \hat{w}_2(\rho \xi_2) \cdot \rho \cdot \delta_0(\xi_1)
\quad \mbox{and} \quad
\hat{\cL} = \delta_0(\xi_2).
\]
Thus the action of $w\cL$ on a continuous function $f$ is given by
\beq \label{eq:actionofwL}
2\pi \langle w\cL , f \rangle  = \langle \cL , \hat{w} \star \hat{f} \rangle
= \int (\hat{w} \star \hat{f})(\xi_1,0) d\xi_1.
\eeq
Conceptually, $w\cL$ is a straight curve fragment; our analysis of $\cC$
in Section \ref{sec-CurvilinearSing} will reduce to the study of this case.

\begin{figure}
\centering
\includegraphics[height=2.3in]{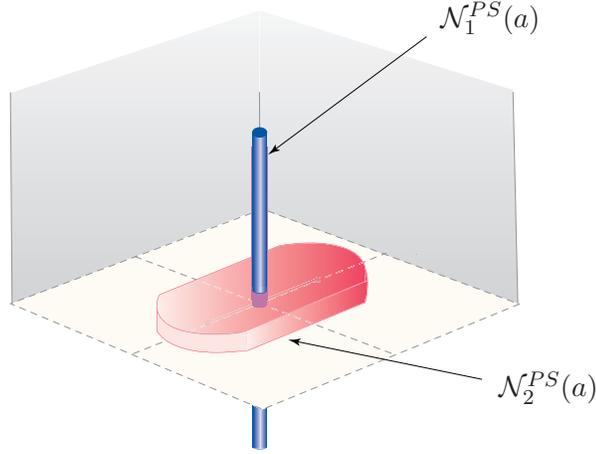}
\put(-30,160){$\cN_{1}^{PS}(a)$}
\put(-8,20){$\cN_{2}^{PS}(a)$}
\caption{The tubes $\cN_{1}^{PS}(a)$ and $\cN_{2}^{PS}(a)$ in phase space.
The clusters of significant coefficients correspond roughly to phase-space support regions
overlapping these tubes.}
\label{fig:tubes}
\end{figure}

Define a
tube in phase space, in which the significant curvelet coefficients will be located.
This will now be a neighborhood of $WF(w\cL)$, defined by
\[
\cN_{2}^{PS}(a) = \{b \in \bR^2 : d_2(b,\{0\} \times [-2\rho,2\rho]) \le D_2(a)\} \times [0,\sqrt{a}],
\]
where
\[
D_2(a) = a^{(1-\eps)}.
\]
For an illustration of $\cN_{2}^{PS}(a)$, and its relation to $\cN_{1}^{PS}(a)$,
we refer to Figure \ref{fig:tubes}. The actual definition of the cluster of curvelet coefficients is much more involved.
In Lemma \ref{lemm:cluster-wS}, we will introduce a first set which helps to determine
its location.

Several bounds will control the curvelet coefficients of a linear singularity.
Lemma \ref{lem-tube} gives
\beq \label{trivialestimate}
|\langle w\cL ,  \gamma_{a,b,\theta}  \rangle | \leq c \cdot 2^{j/4} , \qquad \forall a,b,\theta ;
\eeq
in fact that lemma even gives a decay estimate, as the microlocation $(b,\theta)$ moves
away from $(\{0\} \times [-\rho,\rho]) \times \{0\}$.  In the situations where we would
use that decay estimate, the next lemma is more convenient.

\begin{lemma} \label{lemma:linecurvelet_estimate_close0}
Suppose that $\theta \in [0,\sqrt{a}]$, and set
\[
\tau := \cos \theta \sin \theta (a^{-1}-a^{-2}),\quad d_1^2 := b_1^2(\sigma_2^2-\sigma_1^{-2}\tau),
\]
and
\[
d_2^2
:= \left\{ \begin{array}{ccl}
\min_\pm  \left( (\pm \rho-b_2)\sigma_1 + \sigma_1^{-1}b_1\tau\right)^2 & : & b_2-\sigma_1^{-2}b_1\tau \not\in[-\rho,\rho],\\
0 & : & b_2-\sigma_1^{-2}b_1\tau \in[-\rho,\rho],
\end{array} \right.
\]
where
\[
\sigma_1 = (a^{-2} \sin^2 \theta + a^{-1} \cos^2 \theta)^{1/2} \quad \mbox{and} \quad
\sigma_2 = (a^{-1} \sin^2 \theta + a^{-2} \cos^2 \theta)^{1/2}.
\]
Then, for $N=1, 2, \ldots$,
\[
|\langle w\cL ,  \gamma_{a,b,\theta}  \rangle |
\le c_N \cdot a^{-3/4} \cdot \sigma_1^{-1} \cdot \langle d_1 \rangle^{-1}
\cdot \langle|(d_1,\sigma_1 d_2)|\rangle^{2-N}.
\]
\end{lemma}

In some cases, spatial decay alone is insufficient and we also need to exploit
directional localization; for such cases we employ the following lemma.

\begin{lemma} \label{lemma:linecurvelet_estimate_other}
Suppose that $\theta \in (\sqrt{a},\pi)$. Then, for $L,M=0,1, 2, \ldots$,
\begin{eqnarray*}
\absip{w\cL}{\gamma_{a,b,\theta}}
& \le & c_{L,M} \cdot a^{-1/4} \cdot |\cos \theta| \cdot e^{-\rho \frac{|\sin \theta|}{2a}} \cdot \langle |b_1| \rangle^{-L} \cdot
(a^{1/2} |\sin \theta| + a |\cos \theta|)^L\\
& &  \cdot \langle |b_2| \rangle^{-M} \cdot (\rho + a^{1/2} |\cos \theta| + a |\sin \theta|)^M.
\end{eqnarray*}
\end{lemma}

Both previous lemmas will be proved in Subsections \ref{subsec:proof_of_linecurvelet_estimate_close0}
and \ref{subsec:proof_of_linecurvelet_estimate_other}, respectively.
Together, they imply that the curvelet frame coefficients of  $w\cL$ are
sparse. Indeed, in the directional panels where $\theta$ is close to $0$,
we have about $2^{j/2} \rho $ significantly nonzero coefficients, which
are bounded by  $c 2^{j/4}$, while in the directional panels
where $\theta$ is far from zero, we have few significantly nonzero coefficients.
Formally,

\begin{lemma} \label{Lemma-sparse-sum}
Let $\alpha_j =  (\langle w\cL_j ,  \gamma_\eta  \rangle )_\eta$
denote the curvelet frame coefficients of $w\cL_j = F_j \star w\cL$.
For each $p > 0$, there is $c_p > 0$ so that,
\[
    \| \alpha_j \|_p \leq c_p  \cdot 2^{j \left( (1/2+\eps)/p + 1/4 \right)}, \qquad j > j_0.
\]
\end{lemma}

This will be proved in Section \ref{proof-sparse-sum}.
The  next  result, making
 precise the location of the significant coefficients, is proved in Section  \ref{proof-cluster}.

\begin{lemma} \label{lemm:cluster-wS}
Put
\[
  \t{\cS}_{j} = \{ (j,k,\ell) \in \Delta_j^{\pm 1} : (b_{j,k,\ell},\theta_{j,\ell}) \in \cN_{2}^{PS}(a_j)\}.
\]
Then
\[
     \| \alpha_j 1_{ (\t{\cS}_j)^c } \|_1 = O(1), \qquad j \to \infty.
\]
\end{lemma}


\subsection{Proof of Lemma \ref{Lemma-sparse-sum}}
\label{proof-sparse-sum}

We first observe that the full curvelet coefficient vector $\alpha_j$
is simply the extension of $(\langle w\cL_j ,  \gamma_\eta  \rangle )_{\eta \in \Delta^{\pm 1}_j}$
to scales away from
 $\{j-1,j,j+1\}$ by zero filling.
Also WLOG we can assume that
$\eta \in \Delta_j$, since the terms related to the scales $j-1$ and $j+1$
only change the constant factor of the final estimate independent on $j$.

Define the following four regions in phase space:
\begin{eqnarray*}
\cN_{2}^{PS}(a) & = & \{b \in \bR^2 : d_2(b,\{0\} \times [-2\rho,2\rho]) \le D_2(a)\} \times [0,\sqrt{a}],\\
\cN_{3}^{PS}(a) & = & (\{b \in \bR^2 : d_2(b,\{0\} \times \bR) \le D_2(a)\} \times [0,\sqrt{a}]) \setminus \cN_{2}^{PS}(a),\\
\cN_{4}^{PS}(a) & = & (\bR^2 \times [0,\sqrt{a}]) \setminus (\cN_{2}^{PS}(a) \cup \cN_{3}^{PS}(a)),\\
\cN_{5}^{PS}(a) & = & \bR^2 \times (\sqrt{a},\pi).
\end{eqnarray*}
We now split the norm $\| \alpha_j \|_p$
according to each curvelet coefficient's microlocation.
For a phase space set $\cN$ write $\eta \sim \cN$,
meaning the set $\{ \eta: \eta \in \Delta_j \mbox{ and }  (b_{j,k,\ell},\theta_{j,\ell}) \in \cN\}$.
We have the following decomposition:
\begin{eqnarray}\nonumber
\| \alpha_j \|_p^p & = & \sum_{\eta \sim  \cN_{2}^{PS}(a_j)} \absip{w\cL_j}{\gamma_{\eta}}^p +
\sum_{\eta \sim \cN_{3}^{PS}(a_j)} \absip{w\cL_j}{\gamma_{\eta}}^p\\ \label{eq:Lemma-sparse-sum_0}
& &  + \sum_{\eta \sim \cN_{4}^{PS}(a_j)} \absip{w\cL_j}{\gamma_{\eta}}^p
+\sum_{\eta \sim \cN_{5}^{PS}(a_j)} \absip{w\cL_j}{\gamma_{\eta}}^p.
\end{eqnarray}

\begin{figure}
\centering
\includegraphics[height=2.5in]{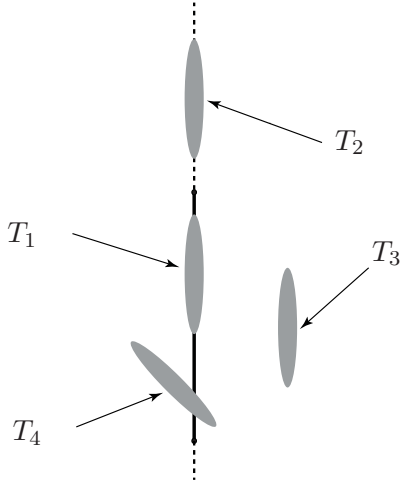}
\put(-137,90){$T_1$}
\put(-13,125){$T_2$}
\put(1,83){$T_3$}
\put(-135,15){$T_4$}
\caption{Curvelets associated with the cases $T_1$--$T_4$. The line
segment is the support of $w\cL_j$. The dotted line is the affine extension of
that segment.}
\label{fig:Lemma63}
\end{figure}

We have the following approximate equivalences:
\begin{eqnarray*}
\eta \sim \cN_{2}^{PS}(a_j) & \approx &
\{(j,(k_1,k_2),0) : |k_1| \leq 2^{j\eps}, k_2 \in [-2\rho/\sqrt{a},2\rho/\sqrt{a}]\},\\
\eta \sim  \cN_{3}^{PS}(a_j)& \approx &
\{(j,(k_1,k_2),0) : |k_1| \leq 2^{j\eps}, k_2 \not\in [-2\rho/\sqrt{a},2\rho/\sqrt{a}]\},\\
\eta \sim \cN_{4}^{PS}(a_j) & \approx &
\{(j,k,0) : k \in \bZ^2, |k_1| \geq 2^{j\eps}\},\\
\eta \sim \cN_{5}^{PS}(a_j) & \approx &
\{(j,k,\ell) : \ell=1, ..., a_j^{-1/2}-1\}.
\end{eqnarray*}

Thus, in place of the original continuum-domain splitting \eqref{eq:Lemma-sparse-sum_0}, we consider
instead the `discrete-domain' splitting
\begin{eqnarray}\nonumber
\| \alpha_j \|_p^p & = & \sum_{|k_1| \leq 2^{j\eps}} \sum_{|k_2| \leq 2\rho/\sqrt{a}}\absip{w\cL_j}{\gamma_{j,(k_1,k_2),0}}^p +
\sum_{|k_1| \leq 2^{j\eps}} \sum_{|k_2| > 2\rho/\sqrt{a}} \absip{w\cL_j}{\gamma_{j,(k_1,k_2),0}}^p\\ \nonumber
& &  + \sum_{|k_1| > 2^{j\eps} } \absip{w\cL_j}{\gamma_{j,k,0}}^p +
\sum_{\ell=1}^{a_j^{-1/2}-1} \sum_{k \in \bZ^2} \absip{w\cL_j}{\gamma_{j,k,\ell}}^p\\[2ex] \label{eq:Lemma-sparse-sum_1}
& = & T_1 + T_2 + T_3 + T_4.
\end{eqnarray}
These terms correspond, respectively, to nearly vertical curvelets lying on the line segment singularity ($T_1$),
nearly vertical curvelets centered elsewhere on the line containing the line segment ($T_2$),
nearly vertical curvelets centered elsewhere ($T_3$),
and all other curvelets ($T_4$), as illustrated in Figure \ref{fig:Lemma63}.


To estimate $T_1$, use (\ref{trivialestimate})
\beq \label{eq:Lemma-sparse-sum_2}
\sum_{k_2 =-2\rho/\sqrt{a}_j}^{2\rho/\sqrt{a}_j} \absip{w\cL_j}{\gamma_{j,(0,k_2),0}}^p
\le c \cdot 2^{j \frac p4 } \cdot \sum_{k_2 =-\rho/\sqrt{a}_j}^{\rho/\sqrt{a}_j} 1
\le c \cdot a_j^{-1/2 + p/4},
\eeq
so $T_1 \leq c \cdot 2^{j\eps} a_j^{-1/2 -p/4} = c 2^{j(1/2 + p/4 + \eps)}$.
Here and below, when we write a sum taken over integers with non-integer bounds, we implicitly mean
that the sum extends over all integers between the bounds.

To derive estimates for $T_2$--$T_3$, we first transfer the estimates for $\theta=0$
derived in Lemma \ref{lemma:linecurvelet_estimate_close0} (in terms of continuum parameters)
to statements about
$|\langle w\cL_j ,  \gamma_{j,k,\ell}  \rangle |$, in terms of the discrete lattice parameters.

\begin{lemma} \label{lemma:linecurvelet_estimate_close0_lattice}
For  $N=1,2,\ldots$,
\[
|\langle w\cL_j ,  \gamma_{j,k,0}  \rangle | \leq
c_N \cdot a_j^{-1/4} \cdot \langle |k_1| \rangle^{-1} \cdot
\langle a_j^{-1} [|a_jk_1|^2 + \min_\pm |a_j^{1/2}k_2\pm\rho|^2]^{1/2} \rangle^{2-N}.
\]
\end{lemma}

\noindent
{\bf Proof.}
This follows directly from the `in particular'-part of Lemma \ref{lemma:linecurvelet_estimate_close0} and the relation
between continuous coefficients and lattice parameters given by \eqref{eq:defilattice}.
\qed


To estimate $T_2$, let $\cK_2 = \{k_2  \in \bZ : k_2 > 2 \rho/\sqrt{a_j}\}$. Then
\begin{eqnarray*}
\lefteqn{\sum_{|k_2|  > 2\rho/\sqrt{a}_j} \absip{w\cL_j}{\gamma_{j,(0,k_2),0}}^p}\\
& \le & c_{N,p} \cdot a_j^{-p/4} \sum_{\cK_2}
\langle a_j^{-1}\min_\pm |a_j^{1/2}k_2 \pm \rho| \rangle^{(2-N)p}\\
& \le & c_{N,p} \cdot a_j^{-p/4} \sum_{\cK_2} \langle | a_j^{-1/2}k_2-a_j^{-1}\rho| \rangle^{(2-N)p}.
\end{eqnarray*}
For $(N-2)p > 1$, we have
\[
\int_{2\rho/\sqrt{a}_j}^\infty \langle | a_j^{-1/2}x-a_j^{-1}\rho| \rangle^{(2-N)p}dx
= a_j^{1/2} \int_{\rho/{a_j}}^\infty \langle | y| \rangle^{(2-N)p}dy
\le c_{N,p} \cdot a_j^{1/2 + ((N-2)p-1) }.
\]
This estimate concerns $k_1 = 0$. For the other cases with $|k_1| < 2^{j\eps}$ we use this
same estimate, getting
\beq \label{eq:Lemma-sparse-sum_3}
T_2
\le c_{N,p} \cdot 2^{j(-1/2 + \eps)}.
\eeq

To estimate $T_3$ let $\cK_3 := \{ k \in \bZ^2 , |k_1| > 2^{j\eps} \}$
and choose $N$ so that  $\eps \cdot (N-2)p > 3 \frac 14$. Then
\[
\sum_{k \in \cK_3} \absip{w\cL_j}{\gamma_{j,k,0}}^p
 \le  c_N \cdot a_j^{-p/4} \sum_{k \in \cK_3}
\langle [|k_1|^2 + \min_\pm |a_j^{-1/2}k_2\pm a_j^{-1} \rho|^2]^{1/2} \rangle^{(2-N)p}.
\]
Partition the set $\cK_3 =  \cK_3^0 \cup \cK_3^1$, where $\cK_3^0 = \{ |k_2| < 2^{j}\rho\} \cap \cK_3$.
The sum over $\cK_3^0$ involves sites where $k_2$ may as well be zero;
it is not asymptotically larger than the LHS in this display:
\begin{eqnarray*}
   4 \cdot  \int_0^{2^{j}\rho} \int_{2^{j\eps}}^\infty \langle |x_1| \rangle^{(2-N)p} dx_1
    &= &  2^{j+2} \rho \cdot \int_{2^{j\eps}}^\infty \langle |x_1| \rangle^{(2-N)p} dx_1  \leq c \cdot 2^{j(1 - \eps\cdot ((N-2)p - 1))}.
\end{eqnarray*}
Since $\eps \cdot (N-2)p > 5/4 + \eps$, this last term is $O(2^{-j/4})$.
The sum over $\cK_3^1$ is not asymptotically larger than the LHS of the next display;
the RHS uses Lemma \ref{lemm:rayintegral}:
\[
     \int_{2^{j} \rho}^\infty \int_{2^{j\eps}}^\infty \langle |(x_1,x_2)| \rangle^{(2-N)p} dx_1  dx_2  \leq
         c \langle | ( 2^{j\eps}, 2^{j}\rho) |\rangle ^{(2-N)p + 2} \leq 2^{j(1 - \eps\cdot ((N-2)p - 2))}
\]
Since $\eps \cdot (N-2)p > 3 \frac 14$, this last term is $O(2^{-j/4})$.
%
We conclude that
\beq \label{eq:Lemma-sparse-sum_4}
\sum_{k \in \cK_3} \absip{w\cL_j}{\gamma_{j,k,0}}^p
\le c_{N,p} \cdot a_j^{(1-p)/4}.
\eeq

Before estimating $T_4$, we translate Lemma \ref{lemma:linecurvelet_estimate_other} into
a simple form involving discrete
 curvelet parameters.
\begin{lemma} \label{lemma:linecurvelet_estimate_other_lattice}
Let \begin{eqnarray*}
T_{j,k,\ell} & = & \langle |a_j k_1\cos \theta_{j,\ell} -a_j^{1/2} k_2\sin \theta_{j,\ell} | \rangle^{-1} \cdot
 (a_j^{1/2} |\sin \theta_{j,\ell}| + a_j |\cos \theta_{j,\ell}|)\\
& & \cdot \langle |a_j k_1\sin \theta_{j,\ell} +a_j^{1/2} k_2\cos \theta_{j,\ell} | \rangle^{-1} .
\end{eqnarray*}
There exist constants $c_N$ so that, for $j,k$, $N=1,2,\ldots$,
\[
\absip{w\cL_j}{\gamma_{j,k,\ell}}
\le c_{N} \cdot a_j^{-1/4} \cdot e^{-\rho \frac{|\sin \theta_{j,\ell}|}{2a_j}} \cdot T_{j,k,\ell}^{N}.
\]
\end{lemma}

\noindent
{\bf Proof.}
This follows directly from Lemma \ref{lemma:linecurvelet_estimate_other},
from $\langle v/u \rangle < \langle v \rangle /u$ for $0 < u \leq 1$ and the relation \eqref{eq:defilattice}
between continuous coefficients and lattice parameters.
\qed

By Lemma \ref{lemma:linecurvelet_estimate_other_lattice}, the term $T_4$ can now be estimated by
\begin{eqnarray}\nonumber
T_4 & \le & c_{N} \cdot a_j^{-p/4} \cdot \sum_{\ell=1}^{a_j^{-1/2}-1} e^{-\rho p \frac{|\sin \theta_{j,\ell}|}{2a_j}}
\Big[\sum_{k \in \bZ^2} T_{j,k,\ell}^{Np} \Big]. \ \label{eq:Lemma-sparse-sum_5}
\end{eqnarray}

Let $\cB_{j,\ell}$ be the rotated anisotropic cartesian grid of curvelet coefficient
locations at scale $j$ and orientation $\theta_{j,\ell}$.  Note that for $N$ large,
\[
    a_j^{3/2} \sum_{b \in \cB_{j,\ell}} \langle | b_1 | \rangle^{-N} \cdot \langle | b_2 |\rangle^{-N}
     \goto \int  \int \langle | b_1 | \rangle^{-N} \cdot \langle | b_2 |\rangle^{-N}  db_1 db_2 , \qquad j \goto \infty .
\]
Indeed, the function $F(b) =  \langle | b_1 | \rangle^{-N} \cdot \langle | b_2 |\rangle^{-N}$ is
smooth and the above display just expresses the fact that Riemann sums of $F$ converge to the integral of $F$.
In fact it is quite evident that the convergence is uniform in  $\theta$. We conclude that
\[
   \max_\ell \sum_{k \in \bZ^2} T_{j,k,\ell}^{Np} \leq  c \cdot a_j^{-3/2} \cdot  a_j^{-Np/2} .
\]

We obtain
\[
T_4 \le c_{N,p} \cdot a_j^{-p/4} \cdot a_j^{-Np/2} \cdot
a_j^{-3/2} \cdot \sum_{\ell=1}^{a_j^{-1/2}-1} e^{-\rho p \frac{|\sin \theta_{j,\ell} |}{2a_j}}.
\]
On the interval $0 < \omega \leq \pi/2$, $\sin(\omega)/\omega \geq 2/\pi$.
We have $|\sin (\pi \ell \sqrt{a_j})| \ge 2 \ell \sqrt{a_j}$, $0 \leq \ell < a_j^{1/2}/2$.
Hence, using $|\sin(\pi/2 + \omega)| = |\sin(\pi/2 - \omega)|$,
\[
\sum_{\ell=1}^{a_j^{-1/2}-1} e^{-\rho p \frac{|\sin \theta_{j,\ell} |}{2a_j}} \leq 2 \cdot \sum_{\ell=1}^{a_j^{-1/2}/2} e^{-\rho p 2 \ell a_j^{-1/2}}.
\]
Summing the geometric series  $\sum_{\ell=1}^{\infty} z^\ell$ with $z = e^{-\rho p 2  a_j^{-1/2}}$,
we finally obtain
that for all $N$ with $Np$ sufficiently large:
\beq \label{eq:Lemma-sparse-sum_7}
\sum_{\ell=1}^{a_j^{-1/2}-1} \sum_{k \in \bZ^2} \absip{w\cL_j}{\gamma_{j,k,\ell}}^p
\le c_{N,p} \cdot a_j^{-p/4-Np/2-3/2} \cdot e^{-\rho p 2  a_j^{-1/2}} \le c_{N,p} \cdot a_j^{N}.
\eeq

Summarizing our bounds on  $T_1$--$T_4$ and using \eqref{eq:Lemma-sparse-sum_1},
we obtain:

\begin{lemma}
\label{lemma:finally_correct_estimates}
For $j$, $N=1, 2, \ldots,$ and $p > 0$, the following holds.
\begin{itemize}
\item[{\rm (i)}] We have
\[
\sum_{\{\eta \in \Delta_j^{\pm 1} : (b_{j,k,\ell},\theta_{j,\ell}) \in \cN_{2}^{PS}(a_j)\}}
\absip{w\cL_j}{\gamma_{\eta}}^p \le c_{N,p} \cdot a_j^{-(1/2+ \eps) -p/4}.
\]
\item[{\rm (ii)}] We have
\[
\sum_{\{\eta \in \Delta_j^{\pm 1} : (b_{j,k,\ell},\theta_{j,\ell}) \not\in \cN_{2}^{PS}(a_j)\}}
\absip{w\cL_j}{\gamma_{\eta}}^p \le c_{N,p} \cdot a_j^{(1-p)/4}.
\]
\end{itemize}
\end{lemma}

\noindent
{\bf Proof.}
Again reducing to scale $j$ and to the discrete
setting as in the proof of Lemma \ref{Lemma-sparse-sum},
 (i) follows from $T_1$, i.e., from \eqref{eq:Lemma-sparse-sum_2}. (ii) follows from $T_2$--$T_4$, i.e.,
from \eqref{eq:Lemma-sparse-sum_3}, \eqref{eq:Lemma-sparse-sum_4}, and \eqref{eq:Lemma-sparse-sum_7}.
\qed\\[1ex]
Finally, this lemma now implies that
\[
\| \alpha_j \|_p^p \le c_{N,p} \cdot a_j^{-(1/2+\eps) - p/4},
\]
which is what was claimed in Lemma \ref{Lemma-sparse-sum}. \qed

\subsection{Proof of Lemma \ref{lemm:cluster-wS}}
\label{proof-cluster}

Using the special properties of our subband filters, WLOG we can assume that $\eta \in \Delta_j$, and
can conclude that
\[
\| \alpha_j 1_{ (\t{\cS}_j)^c } \|_1
\le c \cdot \sum_{\eta \not\sim \cN_{2}^{PS}(a_j)}
\absip{w\cL_j}{\gamma_{\eta}}^p.
\]
Applying Lemma \ref{lemma:finally_correct_estimates}(ii), we obtain
\[
\| \alpha_j 1_{ (\t{\cS}_j)^c } \|_1
\le c_{N,1} \cdot a_j^{0} = O(1), \qquad j \to \infty.
\]
\qed

\section{Sparse Expansion of a Curvilinear Singularity}
\label{sec-CurvilinearSing}


Continuing our `infrastructure development', we now study properties of
curvelet coefficients of a curved singularity.  The strategy is to smoothly
partition the curve into pieces and then straighten each piece, enabling us
to apply results from the previous section.




\newcommand{\length}{\mbox{length}}

\subsection{Tubular Neighborhood}
\label{sec:tubular}

First, we develop a quantitative `tubular neighborhood
theorem'. By regularity, we note that the radius of curvature of $\tau$
is bounded below, by $r > 0$ say.
We can find $\rho$ small compared
to $r$ and an integer $m$ so that
\[
   m \cdot \rho = \length(\tau)
\]
and so that the integrated curvature of $\tau$ on each interval
$[(i-1)\rho,(i+1)\rho]$ is controlled:
\beq \label{curvBound}
     \int_{(i-1)\rho}^{(i+1)\rho} |\tau^{''} (t) | dt \leq \eps .
\eeq
Consider the following local coordinate system
in the vicinity of $\tau$. Let $t_i =  i \rho$, for $i=0 , \dots , m$, and
$T^i = [t_{i-1}, t_{i+1}]$ for $i \in 1, \dots, m-1$. If $\tau$ is a closed curve,
let $T^0 = [t_{m-1},t_1]$ and $T^m = T^0$ (as $\tau(t_0) \equiv \tau( t_m)$).
Let $n_i$ be some choice of unit normal  vector to $\tau(t_i)$.
For $y \in \bR^2$ a point near $\tau(t)$, consider
the closest point in $image(\tau)$; this has arclength parameter
\[
     x_2(y) = \mbox{argmin} \{ | \tau(t) - y | : 0 \leq t \leq \length(\tau) \}
\]
and signed distance parameter
\[
     x_1^i(y) =  \langle n_i , y - \tau(x_2(y)) \rangle  \cdot
     \mbox{min} \{ | \tau(t) - y | : 0 \leq t \leq \length(\tau)\}.
\]
Define the correspondences
\[
  \phi^i ( y) = (x_1^i(y),x_2(y)-t_i) ,  \qquad i=1,\dots , m-1 ,
\]
with similar definitions, slightly amended for the case $i \in \{0,m\}$
if $\tau$ is a closed curve. Recall the curvature bound
 $\eps$ in (\ref{curvBound}).

\begin{lemma} {\bf (Tubular Neighborhood Theorem)}
For sufficiently small $\eps > 0$,
there is some $\eps' > 0$ so that, for
$X_{\eps'} = [-\eps',\eps'] \times [-\rho,\rho]$, we have:

\bitem
\item the correspondence $\phi^i$ is
one-one on the set   $Y_{\eps'}^i \equiv (\phi^i)^{-1}[X_{\eps'}]$,
\item  the mapping $\phi^i : Y_{\eps'}^i \mapsto X_{\eps'}$ is a diffeomorphism, and
\item  the mapping $\phi^i $ extends to a diffeomorphism from $\bR^2$
to $\bR^2$ which reduces to the identity outside a compact set.
\eitem
\end{lemma}
In what follows, $\phi^i$ always denotes the extended diffeomorphism
from $\bR^2$ to $\bR^2$.

\begin{figure}
\centering
\includegraphics[height=1.75in]{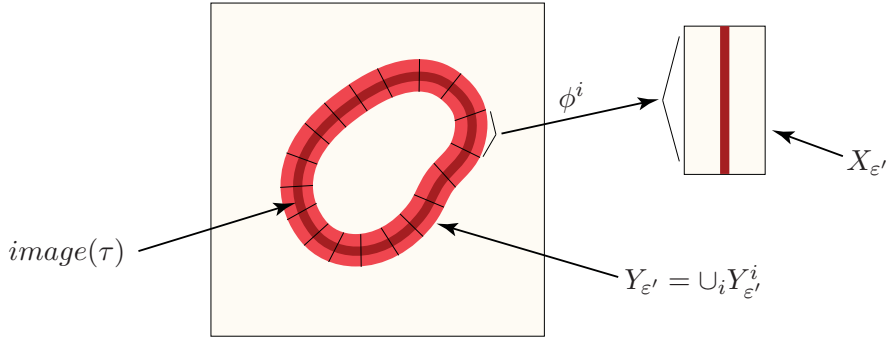}
\put(-316,29){$image(\tau)$}
\put(-83,18){$Y_{\eps'} = \cup_i Y_{\eps'}^i$}
\put(-108,87){$\phi^i$}
\put(1,63){$X_{\eps'}$}
\caption{The tubular neighborhood $Y_{\eps'} = \cup_i Y_{\eps'}^i$ of $image(\tau)$ and the
mapping $\phi^i : Y_{\eps'}^i \mapsto X_{\eps'}$.}
\label{fig:tube}
\end{figure}

The set $Y_{\eps'} = \cup_i Y_{\eps'}^i $ is a tubular neighborhood of $image(\tau)$
on which we have nice local coordinate systems, see Figure \ref{fig:tube}.
This will allow us to locally {\it bend} the curve $\tau$ into something straight.

\subsection{Cutting into pieces}
\label{subsec:cutting}

Choose a $C^\infty$ function
$w_2: \bR \mapsto [0,1]$ (cf. Section \ref{sec-LinearSing}) supported in $[-1,1]$ so that
\[
      w_2(t/\rho) + w_2((t-1)/\rho) = 1, \qquad -1/2 \leq t \leq 0 ,
\]
and
\[
      w_2(t/\rho) + w_2((t+1)/\rho) = 1, \qquad  0 \leq t \leq 1/2 .
\]
In addition,  we require $w_2$ to satisfy
\beq \label{eq:exp_w2}
|\hat{w}_2(\omega)| \le c \cdot e^{-|\omega|},\qquad \omega \in \bR.
\eeq
Define now a smooth partition of unity of $[0, 1]$ using $w_2$:
\[
   w_{2,i}(t/\rho) = w_2((t-t_i)/\rho) , \qquad 1 \leq i \le m-1,
\]
with a modification for $i \in \{0,m\}$ that depends on
whether $\tau$ is closed or not. Then
\begin{equation} \label{propUni}
   \sum_i w_{2,i}(t/\rho) = 1 \qquad \forall t \in [0,1].
\end{equation}
This will allow us to chop the curve $\tau$ into something that
can be bent.

Now define the distributions
\[
     \cC^i = \int_{t_{i-1}}^{t_{i+1}}  w_{2,i} (t/\rho) \delta_{\tau(t)} dt;
\]
the partition of unity property (\ref{propUni})  gives
$ \sum_i \cC^i = \cC$.

We note that
$ \{ \tau(t) :  t \in T^i \} \subset Y_{\eps'}^i$,
and hence $\phi^i$
diffeomorphically straightens  the piece of curve $ \{ \tau(t) :  t \in T^i \} $
into the line segment $\{0\} \times [-\rho,\rho]$.




\subsection{Bending one piece}
\label{sec:bending}

Now consider a diffeomorphism $\phi : \bR^2 \mapsto \bR^2$;
it acts on the distribution $f$ by change of variables
\[
   \phi^\star f = f \circ \phi .
\]
This action induces a linear transformation on the space of
curvelet coefficients. With $\alpha(f)$ the curvelet coefficients
of $f$ and $\beta(f)$ the curvelet coefficients of $\phi^\star f$, we obtain
a linear operator
\[
    M_\phi ( \alpha(f)) = \beta(f) .
\]

It is by now well-known that diffeomorphisms preserve
sparsity of frame coefficients when the frame
is based on parabolic scaling (as with curvelets and shearlets).
For example, the following can be derived from
Hart Smith's work  \cite{Smi98} by a simple atomic decomposition.

\begin{lemma} \label{lemm:operatorquasinorm}
\cite[Theorem 6.1, page 219]{CD05b} For $p > 0$, define the operator quasi-norm
\[
      \| M_\phi \|_{Op,p} =   \max \left\{ \sup_\eta  \| (\langle \gamma_\eta , \phi^\star  \gamma_{\eta'} \rangle )_{\eta} \|_p,
      \sup_{\eta'}  \| (\langle \gamma_\eta , \phi^\star  \gamma_{\eta'} \rangle )_{\eta'}  \|_p \right\}
\]
Let $\phi$ denote a diffeomorphism that reduces to the identity
outside of a compact set. Then for $0 < p \leq 1$,
\[
    \| M_\phi \|_{Op,p}  < c_p < \infty .
\]
\end{lemma}

Far more detailed and precise results on the invariance of curvelet
coefficients under $C^2$ changes of variables, with optimal regularity
conditions, were developed by Cand\`es and Demanet \cite{CD05f} --
as we will see in the next section.

The $p$-triangle inequality  $|a+b|^p \leq |a|^p + |b|^p$
for $p \in (0,1]$ implies the following:

\begin{lemma} For $p \in (0,1]$, a vector $\alpha=(\alpha_\eta)_\eta$,
and a linear operator $M$,
\[
     \| M  \alpha \|_p \leq  \| M \|_{Op,p}  \| \alpha \|_p.
\]
\end{lemma}



\subsection{Gluing pieces together}
\label{sec:glue}

Now define
\[
   \beta_j =  ( \langle \gamma_\eta ,  \cC_j \rangle )_\eta.
\]
From the decomposition $\cC_j = \sum_{i=1}^m \cC^i_j$ we have
\[
   \beta_j =   \sum_{i=1}^m M_{\phi^i} \alpha_j .
\]
This decomposition allows us to relate sparsity of coefficients of the linear singularity to
those of the curvilinear singularity:
\begin{eqnarray*}
     \| \beta_j \|_p & \leq& m^{1/p} \cdot \left( \max_i    \|M_{\phi^i}\|_{Op,p} \right) \cdot  \| \alpha_j  \|_p  .  
\end{eqnarray*}
This decomposition will be useful below; however, the above argument,
which implies sparsity, will not be enough for our main result, which requires
also to know the geometric arrangement of the significant coefficients.
The next section develops a much finer estimation approach.

\section{The cluster $\cS_{2,j}$ and its estimates}
\label{sec-Sec2j}

We finally turn to the definition of the cluster set $\cS_{2,j}$ and
the decisive estimates
\begin{equation} \label{curvClust1}
   \delta_{2,j} = o(\|f_j \|_2 ) , \qquad j \goto \infty,
\end{equation}
and
\begin{equation} \label{curvClust2}
   \mu_c(\cS_{2,j},\{ \gamma_\eta \} ; \{ \psi_\lambda\}) \goto 0 , \qquad j \goto \infty.
\end{equation}
As explained in Section \ref{sec-IntendApplic},
combining these results with the
results of Section \ref{sec-Sec1j} will
complete the proof of Theorem 1.1.

We define the cluster of curvelet coefficients indirectly.
We first define $\t{\cS}_j$, the cluster of significant coefficients
of our `straight' model singularity $w\cL_j$; then by cutting, bending,
and filtering, we induce a cluster for the curvilinear singularity $\cC_j$.
Set
\beq \label{eqtSjDef}
  \t{\cS}_{j} = \{ (j,k,\ell) \in \Delta_j^{\pm 1} : (b_{j,k,\ell},\theta_{j,\ell}) \in \cN_{2}^{PS}(a_j)\}.
\eeq
Lemma \ref{lemm:cluster-wS} shows that this set contains the significant coefficients of $w\cL_j$.

Let $M_{F_j} = (\ip{\gamma_{\eta}}{F_j \star \gamma_{\eta'}})_{\eta,\eta'}$ be the filtering matrix
associated with the filter $F_j$, and recall the definition of the mapping matrix $M_{(\phi^i)^{-1}}$ from
Subsection \ref{sec:bending}. Our analysis will require us to consider their product, hence for the
sake of brevity we define $M_j^i$ to be
\[
M_j^i = M_{F_j} \cdot M_{(\phi^i)^{-1}}
\]
and the entries of this matrix by $M_j^i(\eta,\eta')$.
Further, we let $t_{\eta',n}$ denote the amplitude of the $n$'th largest element of the $\eta'$'th column.
Also let $n_j = 2^{j\eps}$, where $\eps$ was fixed at the beginning of Section \ref{sec-Sec1j}.
We can think of $\eps$ being arbitrarily small, however for our analysis the condition $\eps < 1/28$ will
be sufficient.

Morally, what we would like to do is study  a cluster of curvelet coefficients built from
the cluster pieces
\[
  \cS_{j}^i = \{ \eta : \eta' \in \t{\cS}_{j} \mbox{ and }  |M_j^i(\eta,\eta')| > t_{\eta',n_j} \} .
\]
In words, $\cS_{j}^i$ consists of the `top-$n_j$' curvelet coefficients affected by some significant coefficient
in $\t{\cS}_j$. The overall cluster set would then be made by combing the pieces:
\[
   \cS_{2,j} = \bigcup_i  \cS_{j}^i .
\]

While this morally explains what we do in this
section, it turns out that the exact behavior of $\cS_j^i$ and $\cS_{2,j}$
defined in this natural manner  would be rather delicate.
In fact, this section uses a more robust
definition of cluster set that is similar in spirit;
see (\ref{defCurvClust0})-(\ref{defCurvClust1}) below.
This definition depends on
 some more sophisticated ideas,
 which we now  develop.

\subsection{Decay Estimates for the Curvelet Representation of FIO's}
\label{subsec:decay}

We first recall some results from \cite{CD05f} on
sparsity of curvelet representations of Fourier Integral Operators (FIO's) and
decay estimates of such a representation, which will later on be applied
to the matrix $M_j^i$.

In order to state decay estimates of the curvelet representation of FIO's, we first require a
notion of distance between two curvelet indices. A suitable distance has first been
introduced by Hart Smith in \cite{Smi98}. Our analysis will employ results
obtained by Cand\`{e}s and Demanet in their work on the curvelet representation
of wave propagators \cite{CD05f}, in which they use the following
variation of Hart Smith's distance:
\[
d_{HS}(\eta,\eta')=|\theta_{j,\ell}-\theta_{j',\ell'}|^2 + |b_k - b_{k'}|^2 + |\ip{e_\eta}{b_k - b_{k'}}|,
\]
where
\[
b_k = R_{\theta_{j,\ell}} D_{a_j} k \quad \mbox{and} \quad e_\eta = (\cos(\theta_{j,\ell}),\sin(\theta_{j,\ell})),
\]
and the difference $|\theta_{j,\ell}-\theta_{j',\ell'}|$ is understood to refer to geodesic distance in $\bP^1$.
In \cite{CD05f}, this distance was then extended to derive a distance adapted to discrete curvelet
indices, which means, in particular, including the scaling component. For a pair of
curvelet indices $\eta = (j,k,\ell)$ and $\eta'=(j',k',\ell)$, this so-called {\em
dyadic-parabolic pseudo-distance} is defined by
\[
\omega(\eta,\eta')=2^{|j-j'|}\left(1+\min\{2^j,2^{j'}\}d_{HS}(\eta,\eta')\right).
\]
We will require the following property of this pseudo-distance:

\begin{lemma}\cite[Prop. 2.2 (3.)]{CD05f} \label{lemm:CD1}
For sufficiently large $N>0$,
there is a constant $c_N>0$ such that
\[
\sum_{\eta''} \omega(\eta,\eta'')^{-N} \cdot \omega(\eta'',\eta')^{-N}
\le c_N \cdot \omega(\eta,\eta')^{-(N-1)}, \qquad \forall  \eta,\eta' .
\]
\end{lemma}
Another property which will come in handy is the following estimate:

\begin{lemma}\cite[Proof of Thm. 1.1]{CD05f} \label{lemm:CD2}
There exist some $N>0$
and constant $c_{N}> 0$ obeying
\[
\sum_\eta \omega(\eta,\eta')^{-N} \le c_N , \qquad \forall \eta, \eta' .
\]
\end{lemma}

Before we can state the next result we have to briefly recall some of the key notions
in microlocal analysis. Let $S^\star(\bR^2)$ denote the cosphere bundle of $\bR^2$ --
roughly speaking $\{(b_0,\theta_0) : b_0 \in \bR^2,\, \theta_0 \in \bP^1\}$ --, and
let $\phi$ be a diffeomorphism of $\bR^2$. Then the associated so-called {\em canonical
transformation} $\chi$ maps some element $(b_0,\theta_0)$ of phase space into
$\chi(b_0,\theta_0) = (\phi(b_0),\phi^\star \theta_0)$, where $\phi^\star \theta_0$ is
the codirection into which the codirection $\theta_0$ based infinitesimally at $b_0$ is
mapped under $\phi$. Phrasing it differently, we can say that each diffeomorphism
of the base space $\bR^2$ induces a diffeomorphism of phase space. Such a canonical
transformation induces a mapping of curvelet indices which -- abusing notation -- we
again denote by $\chi$. Since we will consider discrete curvelet coefficients $\eta$,
we have to be careful how to define this extension. In fact, we will define the image of
$\eta$ to be the closest point using the pseudo-distance $\omega$ to the image of $\eta$
under the canonical transformation.
As already remarked in \cite{CD05f}, choosing a different neighbor only affects the
constants in the key inequalities.

The basic insights about parabolic scaling and FIO's are already present in
\cite{Smi98}, implying sparsity of FIO's of order $0$, as explained in  \cite{CD05b}.
But utilizing the dyadic-parabolic pseudo-distance, Cand\`{e}s and Demanet derived
phase space decay estimates for the curvelet representation of FIO's of each order $m$,
which imply sparsity, but also inform about geometry.

\begin{theorem}\cite[Thm. 5.1]{CD05f} \label{theo:CD}
Let $T$ be a Fourier Integral Operator of order $m$ acting on functions of $\bR^2$. Then, for each
$N>0$, there exists some positive constant $c_N$ such that
\beq \label{eq:theoCD}
|\ip{\gamma_{\eta}}{T\gamma_{\eta'}}| \le c_N \cdot 2^{mj'} \omega(\eta,\chi(\eta'))^{-N}.
\eeq
Moreover, for each $0<p<\infty$, $(\ip{\gamma_{\eta}}{T\gamma_{\eta'}})$ is bounded from
$\ell_p$ to $\ell_p$.
\end{theorem}
In the sequel we will use the first part of the result for $m=0$. Let us now turn to the
decay estimate of the cluster approximate error $\delta_{2,j}$.

\subsection{$\cS_{2,j}$ offers low approximation error}
\label{sec:delta2j}

In this section we give two decisive lemmas
which drive our analysis,  and define $\cS_{2,j}$. From now
on $\chi^i$ denotes the extension to curvelet indices of the canonical transformation associated
with $(\phi^i)^{-1}$.

\begin{lemma} \label{lemm:estimateforM1}
For any $N>0$, there exists a positive
constant $c_N$ such that
\beq \label{eq:decayofMji}
|M_j^i(\eta,\eta')| \le c_N \cdot \omega(\eta,\chi^i(\eta'))^{-N}, \qquad \forall \eta, \eta' .
\eeq
\end{lemma}

\noindent
{\bf Proof.}
Given $N>0$, by Theorem \ref{theo:CD}, there exists some positive constant $c_N$ such that
\eqref{eq:theoCD} holds, which implies both
\[
|\ip{\gamma_{\eta}}{F_j \star \gamma_{\eta'}}| \le c_N \cdot \omega(\eta,\eta')^{-(N+1)}
\quad \mbox{and} \quad
|\ip{\gamma_{\eta}}{\gamma_{\eta'} \circ (\phi^i)^{-1}}| \le c_N \cdot \omega(\eta,\chi^i(\eta'))^{-(N+1)}.
\]
Now applying Lemma \ref{lemm:CD1} proves the claim.
\qed

\begin{lemma} \label{lemm:estimateforM2}
There is a constant $c_1> 0$ such that
for each  vector  $\alpha = (\alpha_\eta)_\eta$ ,
\[
\norm{M_j^i \alpha}_1
\le c_1 \cdot \norm{\alpha}_1.
\]
\end{lemma}

\noindent
{\bf Proof.}
This already follows from Lemma 6.3. However, it is instructive to reprove it using
 \eqref{eq:decayofMji} of Lemma \ref{lemm:estimateforM1}:
\[
\norm{M_j^i \alpha}_1 \le \sup_{\eta'} \sum_{\eta} |M_j^i(\eta,\eta')| \cdot \norm{\alpha}_1
\le c_N \cdot \sup_{\eta'} \sum_{\eta} \omega(\eta,\chi^i(\eta'))^{-N}   \cdot \norm{\alpha}_1.
\]
Realizing that the $\sup_{\eta'}$ allows us to omit $\chi^i$ and applying Lemma \ref{lemm:CD2},
\[
\sup_{\eta'} \sum_{\eta} \omega(\eta,\chi^i(\eta'))^{-N}
= \sup_{\eta'} \sum_{\eta} \omega(\eta,\eta')^{-N}
\le c_N. \qquad \mbox{\qed}
\]

\medskip

These two lemmas say that, in place of studying $M_j^i$ and its detailed properties, we
can simply study its majorant $c_N \cdot \omega(\eta,\chi^i(\eta'))^{-N}$.
So fix $N$ large and let $\tilde{M}_j^i$ denote the `model'
\begin{equation} \label{defCurvClust0}
\tilde{M}_j^i(\eta,\eta') = c_N \cdot \omega(\eta,\chi^i(\eta'))^{-N} .
\end{equation}
We define our cluster set $\cS_{2,j}$
 in terms of the model rather than in terms of $M_j^i$, via
 \begin{equation} \label{defCurvClust1}
   \cS_{2,j} = \bigcup_i  \cS_{j}^i,
\end{equation}
where
\[
\cS_{j}^i = \{ \eta : \eta' \in \t{\cS}_{j} \mbox{ and }  |\tilde{M}_j^i(\eta,\eta')| > t_{\eta',n_j} \} .
\]
In this definition, $\cS_j^i$ is not truly the set
of significant coefficients, but rather a set of sites
where significant coefficients could {\it potentially} occur,
given the geometry of the problem; so it is a bit larger.
We still speak of $\cS_{2,j}$ as if it were
exactly the set of significant coefficients.

The set $\t{\cS}_{j}$ is explicitly defined by a tube in phase space; the tube becomes narrower at finer scales and `converges' to
$WF(w\cL)$. The set of potentially significant coefficients $\cS_{j}^i$ is a much thicker set and gets progressively thicker
relative to $\t{\cS}_{j}$ with increasing $j$, however, geometrically the corresponding `tube' is still
becoming very narrow as $j$ increases. This device already appeared in the Heuristics section;
it allows to conveniently bound all the insignificant interactions $M_j^i(\eta,\eta')$;
in particular, see the estimate of $T_2$ in the proof of Lemma \ref{lemma:delta2j}.

We can now prove the estimate  (\ref{curvClust1})
for the cluster approximate error $\delta_{2,j}$.

\begin{lemma}
\label{lemma:delta2j}
\[
 \delta_{2,j} = o(\|f_j \|_2 ) , \qquad j \goto \infty.
\]
\end{lemma}

\noindent
{\bf Proof.}
As in Section \ref{sec:glue}, let $\beta_j = (\ip{\gamma_\eta}{\cC_j})_\eta$ as well as
$\beta_j^i=(\ip{\gamma_\eta}{\cC^i_j})_\eta$, $i\in \{0,\ldots,m-1 \}$. The
decomposition $\cC_j = \sum_i \cC^i_j$ implies
$\beta_j =   \sum_i \beta_j^i$.
Now
\begin{eqnarray*}
\sum_{\eta \in \Delta \setminus \cS_{2,j}} | \beta_j (\eta) |
     &\leq& \sum_i  \sum_{\eta \in \Delta_j^{\pm 1} \setminus \cS_{j}^i} | \beta_j^i (\eta) |   \\
     &\leq&  m \cdot \max_i  \sum_{\eta \in  \Delta_j^{\pm 1} \setminus \cS_{j}^i} | \beta_j^i (\eta) |
\leq m\cdot \max_i  \norm{\beta_j^i \cdot 1_{\Delta_j^{\pm 1} \setminus \cS_{j}^i}}_1   .
\end{eqnarray*}
We now decompose $\norm{\beta_j^i \cdot 1_{\Delta_j^{\pm 1} \setminus \cS_{j}^i}}_1$ into three
components and estimate each separately.  Let $u_j := \lceil \eps \cdot j \rceil$
which will be the `radius' of the scale-neighborhood about scale $j$ that we distinguish from the remaining
(unimportant) scales. In what follows, remember the definition of $\t{\cS}_j$ in (\ref{eqtSjDef});
and let $\t{\cS}_j^{\pm u_j} = \bigcup_{j'} \t{\cS}_{j'}$ . Then
\begin{eqnarray} \nonumber
\lefteqn{\norm{\beta_j^i \cdot 1_{\Delta_j^{\pm 1} \setminus \cS_{j}^i}}_1} \nonumber \\ \nonumber
& \le &
\norm{\t{M}_j^i(\alpha \cdot 1_{\Delta_j^{\pm u_j} \setminus \t{\cS}_j^{\pm u_j}})
\cdot 1_{\Delta_j^{\pm 1} \setminus \cS_{j}^i}}_1
+ \norm{\t{M}_j^i(\alpha \cdot 1_{\t{\cS}_j^{\pm u_j}})
\cdot 1_{\Delta_j^{\pm 1} \setminus \cS_{j}^i}}_1\\ \nonumber
& & + \norm{\t{M}_j^i(\alpha \cdot 1_{\Delta \setminus \Delta_j^{\pm u_j}})
\cdot 1_{\Delta_j^{\pm 1} \setminus \cS_{j}^i}}_1\\ \label{eq:d2j}
& = & T_1 + T_2 + T_3.
\end{eqnarray}

First, let's estimate $T_1$. From Lemma \ref{lemm:estimateforM2}, we obtain
\[
T_1
\le \sum_{j'=j-u_j}^{j+u_j} \norm{\t{M}_j^i(\alpha \cdot 1_{\Delta_{j'} \setminus \t{\cS}_{j'}})}_1
\le c_1 \cdot \sum_{j'=j-u_j}^{j+u_j} \norm{\alpha \cdot 1_{\Delta_{j'} \setminus \t{\cS}_{j'}}}_1.
\]
By Lemma \ref{lemm:cluster-wS},
\[
  \|\alpha \cdot 1_{\Delta_{j'} \setminus \t{\cS}_{j'}}\|_1 = O(1), \qquad j' \to \infty ;
\]
hence,
\beq \label{eq:C1}
T_1
\le c_1 \cdot (2u_j+1)
= O(2^{j\eps}), \qquad j \to \infty.
\eeq

Next, turn to $T_2$. Observe that
\beq \label{eq:C2first}
T_2
\le \sum_{j'=j-u_j}^{j+u_j}
\Big[\sup_{\eta' \in \t{\cS}_{j'}} \sum_{\eta \in \Delta_j^{\pm 1} \setminus \cS_j^i} |\t{M}_j^i(\eta,\eta')|\Big]
\cdot \norm{\alpha \cdot 1_{\t{\cS}_{j'}}}_1 .
\eeq
We now need the following standard lemma about $n$-term approximations.

\begin{lemma} Let $x = (x_i)_i$ denote a sequence of numbers
and let $|x|_{(n)}$ denote the $n$th-largest element in the decreasing rearrangement.
For $0 < p < 1$ we  have the inequality:
\[
   \sum_i | x_i | \cdot 1_{\{ |x_i| \leq |x|_{(n)} \}} \leq c_p \cdot \| x \|_p  \cdot n ^{-(1-p)/p} , \quad n = 1,2, \dots .
\]
\end{lemma}

Recalling that $ \cS_{j}^i $ consists of elements $\eta$ such that $|\t{M}_j^i(\eta,\eta')| > t_{\eta',n_j}$,
we conclude
\[
\sup_{\eta' \in \t{\cS}_{j'}} \sum_{\eta \in \Delta_j^{\pm 1} \setminus \cS_j^i} |\t{M}_j^i(\eta,\eta')|
  \leq c_p \cdot  \left( \sup_{\eta'} \sum_{\eta} |\t{M}_j^i(\eta,\eta')|^p \right)^{1/p} \cdot n_j^{-(1-p)/p}
   \leq c_1 \cdot c_p  \cdot n_j^{-(1-p)/p} .
\]
Also, from Lemma \ref{Lemma-sparse-sum}, we obtain
$\| \alpha \cdot 1_{\t{\cS}_{j'}}\|_1 \leq \| \alpha \cdot 1_{\Delta_{j'}} \|_1 \leq \t{c}_1  2^{j'(3/4+\epsilon)}$.
Choose $p$ so that $p < 4\eps/(3+7\eps)$; returning to \eqref{eq:C2first},
\beq \label{eq:C2}
T_2
\le c  \cdot n_j^{-(1-p)/p} \sum_{j'=j-u_j}^{j+u_j} 2^{j'(3/4+\epsilon)}
\le c \cdot 2^{-(1-p)\eps j/p} \cdot  (2 u_j+1) \cdot 2^{3(j+u_j)/4+j\epsilon}
= O(2^{j \eps}), \; j \to \infty,
\eeq

At last, we consider $T_3$. Notice that
\beq \label{eq:C3first}
T_3
\le \sum_{|j-j'| \ge u_j} \Big[ \sup_{\eta' \in \Delta_{j'}} \sum_{\eta \in \Delta_j^{\pm 1} \setminus \cS_j^i}
|\t{M}_j^i(\eta,\eta')| \Big] \cdot \norm{\alpha \cdot 1_{\Delta_{j'}}}_1.
\eeq
Using the definition of $\t{M}_j$, we proceed as in the proof of
Lemma \ref{lemm:estimateforM2}, and employ the definition of the pseudo-distance $\omega$; we'll obtain
\begin{eqnarray*}
\sup_{\eta' \in \Delta_{j'}} \sum_{\eta \in \Delta_j^{\pm 1} \setminus \cS_j^i} |\t{M}_j^i(\eta,\eta')|
& \le & \t{c}_N \cdot \sup_{\eta' \in \Delta_{j'}} \sum_{\eta \in \Delta_j^{\pm 1} \setminus \cS_j^i}
\omega(\eta,\eta')^{-N}\\
& \le &  c \cdot 2^{-N|j-j'|}.
\end{eqnarray*}
By Lemma \ref{Lemma-sparse-sum}, the second term in \eqref{eq:C3first} can be estimated by
$\norm{\alpha \cdot 1_{\Delta_{j'}}}_1\leq c \cdot 2^{j'(3/4+\epsilon)}$.
We conclude that for $N$ sufficiently large,
\beq \label{eq:C3}
T_3
\le c_N \sum_{|j-j'| \ge u_j} 2^{j'(3/4+\epsilon)-N|j-j'|}
= O(2^{j\eps}), \; j \to \infty.
\eeq
Combining \eqref{eq:C1},  \eqref{eq:C2}, and \eqref{eq:C3} with
\eqref{eq:d2j} yields
\[
   \delta_{2,j} = \sum_{ \eta \in (\cS_j^i)^c } |\beta_j^i (\eta)| =
   O(2^{j\eps})  = o(\|f_j \|_2 ) , \qquad j \goto \infty. \quad  \mbox{\qed}
\]

\subsection{$\cS_{2,j}$ offers low cluster coherence}

This section proves (\ref{curvClust2}), the asymptotically negligible cluster coherence of $\cS_{2,j}$.

\begin{lemma} \label{decayofmuc2}
We have
\[
\mu_{c}(\cS_{2,j}, \{ \gamma_\eta\}  ; \{ \psi_\lambda\}  ) \goto 0, \qquad  j \goto \infty .
\]
\end{lemma}

Before giving the proof, we state two useful lemmas.
Both use the variables $n_j$ introduced earlier.
The first lemma implies that a given significant
curvelet coefficient in the analysis of $w\cL_j$
pushes forward to produce significant coefficients
at roughly the same scale, and  near a certain fixed orientation
and location.  Thus the pushforward acts roughly like a rigid motion.
That first lemma is proved in Section \ref{subsec:proofs_5_1}.

\newcommand{\Fwd}{\mbox{\sc fwd}}
\begin{definition}
Let the canonical transformation $\chi$ be given. For a specific curvelet index $\eta$
the  {\em forward set} of radius $n$ is:
\[
   \Fwd(\eta ;  \chi, n) = \{ \eta' : \omega(\eta', \chi(\eta)) \leq n \} .
\]
\end{definition}
In words, $\Fwd$ is the set of curvelet indices close to the
pushforward of $\eta$ by $\chi$. Note that the forward set covers
the set of significant interactions with $\eta$:
\[
\{\eta' : |\t{M}_j^i(\eta',\eta)| > t_{\eta,n_j}\} \subset \Fwd( \eta; \chi^i , n_j).
\]
Consequently
\[
   \cS_{2,j} \subset \bigcup_{\eta \in \t{\cS}_j}  \bigcup_{i \in \{0, \dots, m\}}   \Fwd( \eta; \chi^i , n_j).
\]

\begin{lemma} \label{lemm:t_estimate}
Let $\eta = (j,k,\ell)$ be a curvelet index, with its image
under the canonical transformation $\chi^i$ denoted by $\chi^i(\eta) = (\t{j},\t{k},\t{\ell})$
for a fixed $i \in \{0,\ldots,m\}$.
There exists some positive constant $c>0$ such that,
for $j$ sufficiently large,
\begin{eqnarray*}
\Fwd(\eta; \chi^i, n_j )
& \subseteq & \{(j',k',\ell') : |j'-\t{j}| \le c \cdot \log n_j,\,
|2^{(j'-\t{j})/2} \t{\ell}- \ell'|^2 \le c \cdot n_j,\,\\
& & 2^{j'}\cdot \max\{|b_{\t{k}}-b_{k'}|^2,|b_{\t{k}}-b_{k'}|\} \le c\cdot  n_j\}
\end{eqnarray*}
\end{lemma}

We conclude that for sufficiently large $j_0$, there exists $c > 0$ so that
\[
        \# \Fwd (\eta ; \chi^i , n_j ) \leq c \cdot n_j^4 , \qquad j > j_0, \quad \eta \in \t{\cS}_j .
\]
\newcommand{\kmin}{d_{min}}
Let $k(\eta)$ denote the $k$-component of  $\eta = (j,k,\ell)$.
Define
\[
          \kmin(\eta,i,n_j)  = \min \{ |k(\eta')| : \eta' \in \Fwd(\eta; \chi^i,n_j)  \}.
\]
We also need the fact that points in the forward set of $\eta$ have spatial components
almost as far from the origin as the spatial component of $\eta$ itself.
\begin{lemma} \label{lemm:lboundk} There are $j_0, c_{1,0},c_{2,0}$ so that
we have
\[
       \kmin(\eta; i , n_j) \geq c_1(j) | k(\eta) | - c_2(j) ,  \qquad j > j_0,
\]
where
\[
         c_1(j) = c_{1,0} / n_j \quad
\mbox{ and } \qquad
       c_2(j)  = c_{2,0} \cdot n_j .
\]
\end{lemma}
The proof is given in the appendix, as is the proof of
\begin{lemma} \label{lemm:sumlemma}
For $N > 2$, there are constants $c_3$,$c_4$,$c_5$ such that
\[
\sum_{k \in \bZ^2} \langle  (a |k| - b)_+ \rangle ^{-N}
\leq (b/a)^2  \cdot ( c_3 + c_4 b^{-N})  + c_5  .
\]
\end{lemma}
With the last three lemmas we can now prove the main result of this section.

\medskip

\noindent
{\bf Proof of Lemma  \ref{decayofmuc2}.}
The definition of $\cS_{2,j}$ implies
\begin{eqnarray*}
\mu_{c}(\cS_{2,j}, \{ \gamma_\eta\}  ; \{ \psi_\lambda\}  )
& = & \max_{j',k'} \sum_{\eta \in \cS_{2,j}} \absip{\gamma_\eta}{\psi_{j',k'}} 
=\max_{j',k'} \sum_{i=0}^{m-1} \sum_{\eta \in \cS_{j}^i} \absip{\gamma_\eta}{\psi_{j',k'}}\\
& \le & m \cdot \max_i \max_{j',k'} \sum_{\eta \in \cS_{j}^i} \absip{\gamma_\eta}{\psi_{j',k'}}
\end{eqnarray*}
WLOG assume that $j' = j$ and $k'=0$, reducing the task
to proving that $\sum_{\eta \in \cS_{j}^i} \absip{\gamma_\eta}{\psi_{j,0}} \to 0$
as $j \to \infty$ for all $i$. We have the estimate
\beq \label{eq:Aimmuc}
\mu_{c}(\cS_{2,j}, \{ \gamma_\eta\}  ; \{ \psi_\lambda\}  )
\le  m \cdot \max_i \sum_{\eta \in \cS_{j}^i} \absip{\gamma_\eta}{\psi_{j,0}}.
\eeq

We can now continue \eqref{eq:Aimmuc} by applying
Lemma \ref{lemm:gamma_psi_estimate} and the three lemmas immediately above:
\begin{eqnarray*}
\sum_{\eta \in \cS_{j}^i} \absip{\gamma_\eta}{\psi_{j,0}}
& \le & \sum_{\eta \in \t{\cS}_{j}} \sum_{\eta' \in \Fwd(\eta)}  \absip{\gamma_{\eta'}}{\psi_{j,0}}\\
& \le & \sum_{\eta \in \t{\cS}_{j}} \sum_{\eta' \in \Fwd(\eta)}   c_N \cdot 2^{-j/4} \cdot \langle | k(\eta') | \rangle^{-N} \\
& \le & c \cdot 2^{-j/4} \cdot \sum_{\eta \in \t{\cS}_{j}}   \#\Fwd(\eta) \cdot  \langle \kmin(\eta) \rangle^{-N} \\
& \le & c \cdot 2^{-j/4} \cdot  n_j^4 \cdot   \sum_{\eta \in \t{\cS}_{j}}\langle (c_1(j) |k(\eta)| - c_2(j))_+ \rangle^{-N} \\
& \le & c \cdot 2^{-j/4} \cdot  n_j^4  \cdot   \left( \frac{c_2(j)^2}{c_1(j)^2}\cdot( c_3 +  c_4 \cdot c_1(j)^2 (c_2(j))^{-N})  + c_5  \right)  \\
& \leq & c \cdot 2^{-j/4} \cdot n_j^8,  \qquad j > j_0.
\end{eqnarray*}

We have assumed that $\eps < 1/32$, so $2^{-j/4} n_j^8 = 2^{-j/4} \cdot  2^{8\eps j}\goto 0$ as $j \goto \infty$;
substituting this into  \eqref{eq:Aimmuc} proves the lemma.
\qed

\section{Discussion}

\subsection{Extensions}
\label{sec:extensions}

So far we focused entirely on a very special separation problem using
very specific tools of harmonic analysis.  Our goal was to show that
a certain set of  questions and results make sense and provide insight.
This is the `tip of the iceberg':  the main results
are susceptible of very extensive generalizations and extensions.

\bitem
\item {\it More General Classes of Objects.}
We may vary the problem,  taking point and curve singularities
whose `strength' is different than the ones
we chose in (\ref{pointdef})-(\ref{curvedef}); however,
always matching the strength of the point singularity
to that of the curve singularity.
For example, consider a `cartoon'
image model, where $\cC$ is a function smooth away from discontinuities,
and the components of the continuity set are bounded by a complex
of smooth curves.  Such cartoons still exhibit curvilinear singularities,
but the singularities are of order zero rather than order $-1$.
For a separation problem with nontrivial asymptotics,
we  replace the point singularity
$|x-x_i|^{-3/2}$ in (\ref{pointdef}) by  $|x-x_i|^{-1/2}$, preserving
an energy-matching condition like (\ref{energyMatch}), with $r^{-1}$ replacing $r^{1}$.
Recall that, without energy matching, the whole problem is trivial.
With such changes  the proof of Theorem \ref{theo:maintheorem}
will run very closely in parallel.  As a general rule,
if
$\langle \cC , f \rangle = \int  (\triangle^\gamma f)(\tau(t)) dt$, where $\gamma$
is a fractional power of the Laplacian $\triangle$, then matching point singularities have strength
$\alpha = (-3 + 4 \gamma)/2$. The case we studied in this paper was $\gamma=0$ and
hence $\alpha=-3/2$.

\item {\it Other Frame Pairs.}
Theorem \ref{theo:maintheorem}
holds without change for many other
pairs of frames and bases. Consider this  pair:
\bitem
\item[$\diamond$] {\it Orthonormal Separable Meyer Wavelets} -- an orthonormal basis of perfectly iso\-tropic
generating elements.
\item[$\diamond$] {\it Shearlets} -- a highly directional tight frame with increasingly
anisotropic elements at fine scales.
\eitem
In this pair, the wavelets are actually orthonormal, and both wavelets
and shearlets correspond very closely to discrete transforms
used in digital image processing.
In digital image processing, the notions of  `radial', `directional', `rotation'
and so on are problematic; both orthonormal wavelets and shearlets
avoid such concepts. At the same time this pair offers the same ability
to sparsify point and curve singularities as the counterparts pair we introduced above.
This allows to
provide a complete methodology for the continuous and discrete
setting (see, e.g., \cite{GKL06,KL07}) as well as for algorithmic realizations (see, e.g., \cite{DKS08,KS07}).

While the proof arguments explicitly cover the
one frame pair we have taken pains to define so far, those arguments extend
immediately to other `compatible' pairs -- where the cross-frame matrices
are almost diagonal in a suitable sense.   This grants us the freedom
to prove results in one system which is convenient, but apply those
to another compatible system.  The arguments showing that
shearlets and curvelets are compatible are supplied
in \cite{DK08b}.
In this paper we discussed the pair radial wavelets/curvelets.
However, all results hold true in a similar way for the pair orthonormal
wavelets/shear\-lets.

\item {\it Noisy Data.}
Are the results studied here robust against small modelling errors?
In fact they are. Consider an image
composed of $\cP$ and $\cC$ with additive noise $\cN$, hence we
 measure $\t{f}=\cP + \cC + \cN$ instead of $f$. We then -- as in the noiseless case --
filter to obtain subband components $\t{f}_j= \cP_j + \cC_j + \cN_j$ and apply $(\mbox{CSep})$
to $\t{f}_j$ to obtain a pair $(\t{W}_j,\t{C}_j)$. Provided that the noise
component $\cN$ has `sufficiently' small curvelet coefficients in the sense
that at each scale $j$ the $\ell_1$ norm of the analysis coefficients satisfies $o(2^{j/2})$ as
$j \to \infty$, we again obtain asymptotically perfect separation:
\[
\frac{ \| \t{W}_j - \cP_j \|_2 + \| \t{C}_j - \cC_j \|_2 }{\| \cP_j\|_2 +
\|\cC_j\|_2 } \goto 0, \qquad j \goto \infty
\]
This can be proved along the lines of the proof of Theorem \ref{theo:maintheorem}.
Indeed, consider a composed signal $\t{S}=S_1^0+S_2^0
+n$ with components $S_1^0$ and $S_2^0$ relatively sparse as
in Proposition \ref{prop:mainestimate}, and noise term $n$
satisfying $\norm{\Phi_1^T n}_1 < \eps$ or $\norm{\Phi_2^T n}_1 < \eps$.
Let $(\t{S}_1^\star,\t{S}_2^\star)$ solve (\ref{PSep}) with $S$ substituted by $\t{S}$.
Then following the proof  of  Proposition \ref{prop:mainestimate}
line by line and adapting the arguments accordingly
shows:
\beq \label{eq:noise}
  \norm{\t{S}_1^\star-S_1^0}_2 + \norm{\t{S}_2^\star-S_2^0}_2
\le \frac{2\delta + 5 \eps}{1-2\kappa}.
\eeq

Substituting Proposition \ref{prop:mainestimate} by  estimate \eqref{eq:noise} in the proof
of Theorem \ref{theo:maintheorem} implies the result on geometric separation
of noisy data stated in Subsection \ref{sec:extensions}.

We conclude that our analysis is indeed stable.
\item {\it Rate of Convergence.}
One might wonder about the rate of separation. The lemmas proven in
Sections \ref{sec-Sec1j}-\ref{sec-Sec2j}
imply the following upper bound on the rate of convergence for $\ell_1$ minimization:
\[
\frac{ \| W_j - \cP_j \|_2 + \| C_j - \cC_j \|_2 }{\| \cP_j\|_2 + \|\cC_j\|_2 } =
O(2^{-j(1/2-\epsilon)}), \qquad j \goto \infty.
\]
Such information might be the key to getting even stronger separation conclusions.

\item {\it Other Algorithms and Other Notions of Separation.}
In the companion paper \cite{DK08a} we show that one pass of alternating
hard thresholding, properly tuned, can achieve asymptotic separation.
Surprisingly, we can even show clean separation at the level of wavefront sets.
\eitem

\subsection{Interpretation as an Uncertainty Principle}

Separation results such as the `birth problem' of $\ell_1$ component separation,
a combination of sinusoids and spikes \cite{DH01}, have been interpreted at that
time as uncertainty principles. As a reminder to the reader, the classical
uncertainty principle states that a signal cannot be highly concentrated in both
time and frequency; and a lower bound is placed on the product of the concentration
in time and in frequency. The core property which allows the separation of sinusoids and
spikes by using a dictionary consisting of the unit basis and the Fourier basis,
is the non-existence of a sparse representation of a signal both in time and in
frequency.

Considering the present separation problem, these core ideas need to be extended,
thereby providing us with yet another interpretation than the one already presented
in the previous sections. The two representation `domains' are now the isotropic
system of wavelets and the anisotropic system of curvelets. Hence, we might regard
the separation result of Theorem \ref{theo:maintheorem} as a statement that a 2D Schwartz distribution cannot be
sparsely represented via analysis coefficients both in the `isotropic world' and
in the `anisotropic world'. In particular, if a 2D Schwartz distribution has a sparse
representation in wavelets, it is not sparse in curvelets and vice versa.
Phrasing it in more general terms, a 2D Schwartz distribution having only isotropic
features cannot be sparsely represented using an anisotropic system, and if it
has only exhibits anisotropic phenomena, it does not possess a sparse representation
in terms of an isotropic system.

Summarizing, comparison with the classical uncertainty principle
shows that we here derive an uncertainty principle for the isotropy-anisotropy relation
instead of the classical time-frequency relation.

\section{Proofs}

\subsection{Proofs of Results from Section \ref{sec:compsep}}
\label{subsec:proofs_1}

\subsubsection{Proof of Proposition \ref{prop:mainestimate}}

\noindent {\bf Proof.}
Since $\Phi_1$ and $\Phi_2$ are tight frames,
\begin{eqnarray*}
\norm{S_1^\star-S_1^0}_2 + \norm{S_2^\star-S_2^0}_2
& = & \norm{\Phi_1^T(S_1^\star-S_1^0)}_2 +
\norm{\Phi_2^T(S_2^\star-S_2^0)}_2\\
& \le & \norm{\Phi_1^T(S_1^\star-S_1^0)}_1 +
\norm{\Phi_2^T(S_2^\star-S_2^0)}_1.
\end{eqnarray*}
Now invoke exact decomposition:
$S_1^0+S_2^0=S=S_1^\star+S_2^\star$.
Rewrite the last display:
\[ \norm{S_1^\star-S_1^0}_2 + \norm{S_2^\star-S_2^0}_2 \le
\norm{\Phi_1^T(S_1^\star-S_1^0)}_1 +
\norm{\Phi_2^T(S_1^\star-S_1^0)}_1.\]
By definition of $\kappa$,
\begin{eqnarray*}
\lefteqn{\norm{\Phi_1^T(S_1^\star-S_1^0)}_1 + \norm{\Phi_2^T(S_1^\star-S_1^0)}_1}\\
& = & \hspace*{-0.2cm} \norm{1_{\cS_1}\Phi_1^T(S_1^\star-S_1^0)}_1 +
\norm{1_{\cS_2}\Phi_2^T(S_1^\star-S_1^0)}_1  +
\norm{1_{\cS_1^c}\Phi_1^T(S_1^\star-S_1^0)}_1 +
\norm{1_{\cS_2^c}\Phi_2^T(S_2^\star-S_2^0)}_1 \\
& \le & \kappa \cdot \left( \norm{\Phi_1^T(S_1^\star-S_1^0)}_1 +
\norm{\Phi_2^T(S_1^\star-S_1^0)}_1 \right) +
\norm{1_{\cS_1^c}\Phi_1^T(S_1^\star-S_1^0)}_1 +
\norm{1_{\cS_2^c}\Phi_2^T(S_2^\star-S_2^0)}_1;
\end{eqnarray*}
use relative sparsity of the subsignals $S_i^0$, $i=1,2$,
\begin{eqnarray} \nonumber
\lefteqn{\norm{\Phi_1^T(S_1^\star-S_1^0)}_1 + \norm{\Phi_2^T(S_1^\star-S_1^0)}_1}\\ \nonumber
& \le & \frac{1}{1-\kappa} (\norm{1_{\cS_1^c}\Phi_1^T(S_1^\star-S_1^0)}_1 + \norm{1_{\cS_2^c}\Phi_2^T(S_2^\star-S_2^0)}_1)\\ \nonumber
& \le & \frac{1}{1-\kappa} (\norm{1_{\cS_1^c}\Phi_1^TS_1^\star}_1 +
\norm{1_{\cS_1^c}\Phi_1^TS_1^0}_1 +
\norm{1_{\cS_2^c}\Phi_2^TS_2^\star}_1+\norm{1_{\cS_2^c}\Phi_2^TS_2^0}_1)\\ \label{eq:l1min_1}
& \le & \frac{1}{1-\kappa} \left(\norm{1_{\cS_1^c}\Phi_1^TS_1^\star}_1+
\norm{1_{\cS_2^c}\Phi_2^TS_2^\star}_1+
\delta \right) .
\end{eqnarray}
Apply  minimality of $S_1^\star$ and $S_2^\star$,
\begin{eqnarray*}
\norm{1_{\cS_1^c}\Phi_1^TS_1^\star}_1+\norm{1_{\cS_1}\Phi_1^TS_1^\star}_1+
\norm{1_{\cS_2^c}\Phi_2^TS_2^\star}_1+
\norm{1_{\cS_2}\Phi_2^TS_2^\star}_1
& =& \norm{\Phi_1^TS_1^\star}_1+\norm{\Phi_2^TS_2^\star}_1\\
& \le & \norm{\Phi_1^T S_1^0}_1+\norm{\Phi_2^T S_2^0}_1.
\end{eqnarray*}
Again use sparsity of the subsignals $S_i^0$,
$i=1,2$,
\begin{eqnarray*}
\lefteqn{\norm{1_{\cS_1^c}\Phi_1^TS_1^\star}_1+ \norm{1_{\cS_2^c}\Phi_2^TS_2^\star}_1}\\
& \le & \norm{\Phi_1^T S_1^0}_1+\norm{\Phi_2^T S_2^0}_1 - \norm{1_{\cS_1}\Phi_1^TS_1^\star}_1 - \norm{1_{\cS_2}\Phi_2^TS_2^\star}_1\\
& \le & \norm{\Phi_1^T S_1^0}_1+\norm{\Phi_2^T S_2^0}_1
+ \norm{1_{\cS_1}\Phi_1^T(S_1^\star-S_1^0)}_1 - \norm{1_{\cS_1}\Phi_1^T S_1^0}_1\\
& & + \norm{1_{\cS_2}\Phi_2^T(S_2^\star-S_2^0)}_1 - \norm{1_{\cS_2}\Phi_2^T S_2^0}_1\\
& \le & \norm{1_{\cS_1}\Phi_1^T(S_1^\star-S_1^0)}_1 +
\norm{1_{\cS_2}\Phi_2^T(S_2^\star-S_2^0)}_1 +
\delta.
\end{eqnarray*}
Using  \eqref{eq:l1min_1}, this leads to
\begin{eqnarray*}
\lefteqn{\norm{\Phi_1^T(S_1^\star-S_1^0)}_1 + \norm{\Phi_2^T(S_1^\star-S_1^0)}_1}\\
& \le & \frac{1}{1-\kappa} \left [\norm{1_{\cS_1}\Phi_1^T(S_1^\star-S_1^0)}_1
+ \norm{1_{\cS_2}\Phi_2^T(S_1^\star-S_1^0)}_1+
2\delta \right]  \\
& \le & \frac{1}{1-\kappa} \left[\kappa \cdot (\norm{\Phi_1^T(S_1^\star-S_1^0)}_1 +
\norm{\Phi_2^T(S_1^\star-S_1^0)}_1)+
2\delta \right].
\end{eqnarray*}
Thus, finally we obtain
\[\norm{S_1^\star-S_1^0}_2 + \norm{S_2^\star-S_2^0}_2 \le \left(1- \frac{\kappa}{1-\kappa}\right)^{-1} \cdot \frac{2\delta}{1-\kappa} =
\frac{2\delta}{1-2\kappa}. \Box\]

\subsubsection{Proof of Lemma \ref{lemm:kappamuc}}
\noindent {\bf Proof.}
For each $f$, we choose coefficient sequences $\alpha_1$ and $\alpha_2$
such that $f=\Phi_1 \alpha_1 = \Phi_2 \alpha_2$ and $\norm{\alpha_i}_1 \le
\norm{\beta_i}_1$ for all $\beta_i$ satisfying $f=\Phi_i \beta_i$, $i=1,2$.
Then, employing the fact that, because $\Phi_1$ and $\Phi_2$ are
tight frames, also $f = \Phi_i \Phi_i^T \Phi_i \alpha_i$, $i=1,2$, we
obtain {\allowdisplaybreaks
\begin{eqnarray*}
\lefteqn{\norm{1_{\cS_1} \Phi_1^T f}_1 + \norm{1_{\cS_2} \Phi_2^T f}_1}\\[1ex]
& = & \norm{1_{\cS_1} \Phi_1^T \Phi_2 \alpha_2}_1 + \norm{1_{\cS_2} \Phi_2^T \Phi_1 \alpha_1}_1\\
& \le & \sum_{i \in \cS_1}\left(\sum_j
\absip{\phi_{1,i}}{\phi_{2,j}} |\alpha_{2,j}|\right)
+ \sum_{j \in \cS_2}\left(\sum_i \absip{\phi_{1,i}}{\phi_{2,j}} |\alpha_{1,i}|\right)\\
& = & \sum_{j}\left(\sum_{i \in \cS_1}
\absip{\phi_{1,i}}{\phi_{2,j}}\right) |\alpha_{2,j}|
+ \sum_i \left(\sum_{j \in \cS_2} \absip{\phi_{1,i}}{\phi_{2,j}}\right) |\alpha_{1,i}|\\[1ex]
& \le & \mu_c(\cS_1,\Phi_1;\Phi_2)
 \norm{\alpha_2}_1 + \mu_c(\cS_2, \Phi_2; \Phi_1) \norm{\alpha_1}_1\\[1ex]
& \le & \max\{ \mu_c(\cS_1,\Phi_1;\Phi_2), \mu_c(\cS_2,\Phi_2;\Phi_1)\} (\norm{\alpha_1}_1+\norm{\alpha_2}_1)\\[1ex]
& \le &\max\{ \mu_c(\cS_1,\Phi_1;\Phi_2), \mu_c(\cS_2,\Phi_2;\Phi_1)\} (\norm{\Phi_1^T \Phi_1 \alpha_1}_1+\norm{\Phi_2^T \Phi_2 \alpha_2}_1)\\[1ex]
& = & \max\{ \mu_c(\cS_1,\Phi_1;\Phi_2), \mu_c(\cS_2,\Phi_2;\Phi_1)\} (\norm{\Phi_1^T
f}_1+\norm{\Phi_2^T f}_1).
\end{eqnarray*}
} 
\hfill \qed

\subsection{Proofs of Results from Section \ref{sec:microlocal}}
\label{subsec:proofs_2}

\subsubsection{Proof of Lemma  \ref{lemm:gamma_psi_estimate} }
\label{pf-Lemm-3.3}

 \noindent
 Using Parseval, $\langle \gamma_{a,b,\theta},\psi_{a_0,b_0} \rangle  =
     2\pi  \int    \hat{\gamma}_{a,b,\theta}(\xi) \hat{ \psi}_{a_0,b_0}(\xi) d\xi$,
we consider
\[
   \int   \hat{\gamma}_{a,b,\theta}(\xi) \hat{ \psi}_{a_0,b_0}(\xi)   d\xi \
=    \int   a_0 W(a_0 r) e^{-i b_0'\xi}   \cdot  a^{3/4}   W(a r) V((\omega-\theta)/\sqrt{a})
e^{-ib'\xi}  d\xi.
\]
Now WLOG we may consider the special case
 $\theta=0$, so that $R_\theta = I$. Recall that
 $W$ is supported on $[1/2,2]$ by construction. Then
 \[
   W(a_0 r)    W(a r) = 0 \qquad  \forall r \geq 0,  | \log_2(a/a_0)| \ge 3.
 \]
Hence we need only consider the case where $| \log_2(a/a_0)|  < 3$,
and in that circumstance we may WLOG take $a = a_0$.
We may also assume $b_0 = 0$. Apply the change of variables
$\zeta = D_{a} \xi$ and $d\zeta = a^{3/2} d\xi$,
\[
\int   \hat{\gamma}_{a,b,\theta}(\xi) \hat{ \psi}_{a,b_0}(\xi)   d\xi
= a^{1/4} \cdot \int W^2(\norm{\zeta_a}) V(\omega(\zeta_a)/\sqrt{a})e^{-i(D_{1/a}b)'\zeta} d\zeta,
\]
where $\zeta_a = (\zeta_1,\sqrt{a}\zeta_2)$ and $\omega(\zeta_a)$ denotes the angular component
of the polar coordinates of $\zeta_a$.
Applying integration by parts, for any $k=1, 2, ...$,
\begin{eqnarray}\nonumber
\absip{\gamma_{a,b,\theta}}{\psi_{a_0,b_0}}
& = &2\pi \cdot a^{1/4} \cdot |D_{1/a}b|^{-k} \left|\int \Delta^k [W^2(\norm{\zeta_a}) V(\omega(\zeta_a)/\sqrt{a})]e^{-i(D_{1/a}b)'\zeta}
d\zeta\right|\\ \nonumber
& \le & 2\pi \cdot a^{1/4} \cdot |D_{1/a}b|^{-k} \int \left|\Delta^k[W^2(\norm{\zeta_a}) V(\omega(\zeta_a)/\sqrt{a})]\right| d\zeta.
\end{eqnarray}
Hence
\begin{eqnarray}\nonumber
\lefteqn{(1+|D_{1/a}b|^{k}) \cdot \absip{\gamma_{a,b,\theta}}{\psi_{a_0,b_0}}}\\ \label{eq:gammapsi_1}
& \le & 2\pi \cdot a^{1/4} \int \left[ \left|W^2(\norm{\zeta_a})\right|\left|V(\omega(\zeta_a)/\sqrt{a})\right|
+ \left|\Delta^k[W^2(\norm{\zeta_a}) V(\omega(\zeta_a)/\sqrt{a})]\right| \right]d\zeta.
\end{eqnarray}
Next we show that, for each $k$, there exists $c_k < \infty$ such that
\beq \label{eq:gammapsi_2}
\int \left[ \left|W^2(\norm{\zeta_a})\right|\left|V(\omega(\zeta_a)/\sqrt{a})\right|
+ \left|\Delta^k[W^2(\norm{\zeta_a}) V(\omega(\zeta_a)/\sqrt{a})]\right| \right]d\zeta \le c_k,\quad \forall \;  a>0.
\eeq
We have
\[
\frac{\partial}{\partial \zeta_1}W^2(\norm{\zeta_a}) = \frac{\partial}{\partial \zeta_1}W^2(\norm{(\,\cdot\,,\sqrt{a}\zeta_2)})(\zeta_1)
\]
and
\[
\frac{\partial}{\partial \zeta_2}W^2(\norm{\zeta_a}) = \sqrt{a}\cdot\frac{\partial}{\partial \zeta_2}W^2(\norm{(\zeta_1,\sqrt{a}\,\cdot\,)})(\zeta_2).
\]
Hence, by induction, the absolute values of the derivatives of $W^2(\norm{\zeta_a})$ are upper bounded independently of $a$.
Also,
\[
\frac{\partial}{\partial \zeta_1} V(\omega(\zeta_a)/\sqrt{a})
= \frac{\partial}{\partial \zeta_1}V(\omega((\,\cdot\,,\zeta_2)_a)/\sqrt{a})(\zeta_1)\cdot g_1(\zeta,a)
\]
and
\[
\frac{\partial}{\partial \zeta_2} V(\omega(\zeta_a)/\sqrt{a})
= \frac{\partial}{\partial \zeta_2}V(\omega((\zeta_1,\,\cdot\,)_a)/\sqrt{a})(\zeta_1)\cdot g_2(\zeta,a),
\]
and tedious computations show that both $|g_1|, |g_2|$ possess an upper bound independently of $a$. Thus,
by induction, the absolute values of the derivatives
of $V(\omega(\zeta_a)/\sqrt{a})$ are upper bounded independently of $a$.
These observations imply \eqref{eq:gammapsi_2}.

Further, for each $k =1, 2, ..$,
\beq \label{eq:gammapsi_3}
\langle|D_{1/a}b|\rangle^k = (1+|D_{1/a}b|^2)^\frac{k}{2} \le \frac{k}{2}(1+|D_{1/a}b|^{k}).
\eeq

To finish, simply combine \eqref{eq:gammapsi_1}, \eqref{eq:gammapsi_2}, and \eqref{eq:gammapsi_3},
and recall that we chose coordinates so that $\theta = 0$.  Translating back to the case
of general $\theta$ gives the full conclusion. \qed

\subsection{Proofs of Results from Section \ref{sec-Sec1j}}
\label{subsec:proofs_3}

\subsubsection{Proof of Lemma \ref{lemma:pointwavelet}}
\label{subsec:proofs_3_1}

{\bf Proof.}
By Parseval,
\[
\langle \psi_{a_{j'},b}, \cP_j \rangle =  2\pi\int  \hat{
\psi}_{a_{j'},b}(\xi)  \hat{\cP}_j(\xi) d\xi =   2\pi\int  a_{j'}
W(a_{j'}r) e^{-ib'\xi} W(a_{j}r) r^{-1/2}  d\xi,
\]
where of course $r = |\xi|$. Now
\[
   W(a_{j'} r)    W(a_{j}r) = 0 \qquad  \forall r \geq 0  , \qquad  |j - j'| > 1,
 \]
 hence we may as well assume that $j'=j$.
Making the change of variables $\zeta = a_j\xi$, $d\zeta = a_{j}^2 d\xi$
and defining the annulus $\cA = \{ \zeta : 1/2 \leq |\zeta| \leq 2 \}$,
\[
\langle \psi_{a_j,b}, \cP_j \rangle = 2\pi \cdot a_{j}^{-1/2} \cdot
\int_{\cA} W^2(|\zeta|)  |\zeta|^{-1/2} e^{-i(b/a_j)'\zeta} d\zeta.
\]
Applying integration by parts, for any $k=0, 1, ...$,
\begin{eqnarray}\nonumber
\absip{\psi_{a_j,b}}{\cP_j}
& = & 2\pi\cdot a_{j}^{-1/2} \cdot |b/a_j|^{-k} \cdot \left|\int_\cA \Delta^k [W^2(|\zeta|)  |\zeta|^{-1/2}] e^{-i(b/a_j)'\zeta} d\zeta\right|\\ \nonumber
& \le & 2\pi\cdot a_{j}^{-1/2} \cdot |b/a_j|^{-k} \cdot \int_\cA \left|\Delta^k [W^2(|\zeta|)  |\zeta|^{-1/2}]\right| d\zeta.
\end{eqnarray}
Hence
\beq \label{eq:psiP_1}
(1+|b/a_j|^{k}) \cdot \absip{\psi_{a_j,b}}{\cP_j}
\le2\pi \cdot a_{j}^{-1/2} \cdot \int_\cA \left[ \left| W^2(|\zeta|)  |\zeta|^{-1/2}\right| + \left|\Delta^k [W^2(|\zeta|)  |\zeta|^{-1/2}]\right|\right] d\zeta.
\eeq
For $W$ suitably chosen,
\[
 \int_\cA \left[ \left| W^2(|\zeta|)  |\zeta|^{-1/2}\right| + \left|\Delta^k [W^2(|\zeta|)  |\zeta|^{-1/2}]\right|\right] d\zeta < \infty.
\]
Further, for each $k =1, 2, ..$,
\[
\langle|b/a_j|\rangle^{k} = (1+|b/a_j|^2)^\frac{k}{2} \le \frac{k}{2}(1+|b/a_j|^{k}).
\]
Infusing these two last observations into \eqref{eq:psiP_1}, for any $N=1, 2, ...$,
\[
\absip{\psi_{a_{j'},b}}{\cP_j} \le c_N \cdot a_{j}^{-1/2} \cdot 1_{|j-j'| \le 1}\cdot  \langle |b/a_j| \rangle^{-N}.
\]
 \qed

\subsection{Proofs of Results from Section \ref{sec-LinearSing}}
\label{subsec:proofs_4}

\subsubsection{Proof of  Lemma \ref{lemma:linecurvelet_estimate_close0}}
\label{subsec:proof_of_linecurvelet_estimate_close0}

We study the situation geometrically, and for each $(a,b,\theta)$ define the line segment
\[
Seg(a,b,\theta) = \Big\{D_{1/a} R_{-\theta} \left(\begin{array}{c}-b_1\\y-b_2\end{array}\right) : |y| \le \rho\Big\}.
\]
Two special points associated to these line
segments will play an essential role in our
estimate; they are defined in

\begin{lemma}
\label{lemma:PLPS}
Retain the definitions for $d_1, d_2, \sigma_1, \sigma_2,$ and $\tau$
from the statement of Lemma \ref{lemma:linecurvelet_estimate_close0}.
Then, for each $a,b,\theta$, the following conditions are fulfilled.
\begin{itemize}
\item[{\rm (i)}] Consider the line
\[
Line(a,b,\theta) = \Big\{D_{1/a} R_{-\theta} \left(\begin{array}{c}-b_1\\y-b_2\end{array}\right) : y \in \bR\Big\} .
\]
The closest point $P_L$ to the origin on $Line(a,b,\theta)$ satisfies
\[
\norm{P_L}_2^2 = b_1^2(\sigma_2^2-\sigma_1^{-2}\tau) = d_1^2.
\]
\item[{\rm (ii)}] Let $P_S$ be the closest point on $Seg(a,b,\theta)$ to the origin. Then
\[
\norm{P_S-P_L}_2^2 = d_2^2.
\]
\end{itemize}
\end{lemma}

\begin{figure}
\centering
\vspace*{-1.5cm}
\includegraphics[height=2.5in]{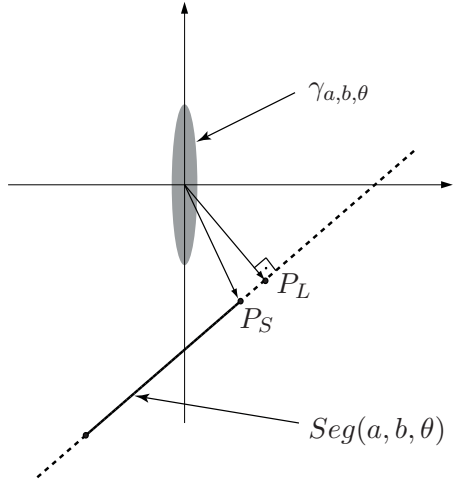}
\put(-81,57){$P_S$}
\put(-67,70){$P_L$}
\put(-55,145){$\gamma_{a,b,\theta}$}
\put(-55,15){$Seg(a,b,\theta)$}
\caption{Relation between a curvelet $\gamma_{a,b,\theta}$, a line segment
$Seg(a,b,\theta)$, its affine hull $Line(a,b,\theta)$, and the points $P_L$ and $P_S$.}
\label{fig:Lemma65}
\end{figure}

Figure \ref{fig:Lemma65} shows a general configuration featuring $P_L$ and $P_S$.

\noindent
{\bf Proof.}
Set
\[
L(y) := D_{{1}/{a}} R_{-\theta} \left(\begin{array}{c}-b_1\\y-b_2\end{array}\right).
\]
Then
\begin{eqnarray} \nonumber
\norm{L(y)}_2^2
& = & \norm{(a^{-1}(-b_1\cos \theta + (y-b_2) \sin \theta),a^{-1/2}(b_1\sin \theta + (y-b_2) \cos \theta)}_2^2\\ \label{eq:PLPS_1}
& = & b_1^2\sigma_2^2 + (y-b_2)^2\sigma_1^2 + 2 b_1 (y-b_2) \tau.
\end{eqnarray}
Since
\[
\frac{d}{dy}\norm{L(y)}_2^2 = 2(y-b_2)\sigma_1^2+2b_1\tau,
\]
it follows by definition of $P_L$ that
\beq \label{eq:PLPS_2}
P_L=L(b_2-\sigma_1^{-2}b_1\tau).
\eeq
Hence, by \eqref{eq:PLPS_1},
\[
d_1^2=\norm{P_L}_2^2 =b_1^2\sigma_2^2 + (-\sigma_1^{-2}b_1\tau)^2\sigma_1^2 + 2 b_1 (-\sigma_1^{-2}b_1\tau) \tau=  b_1^2(\sigma_2^2-\sigma_1^{-2}\tau).
\]
This proves (i).

To prove (ii), observe that, by \eqref{eq:PLPS_2},  $P_L\in Seg(a,b,\theta)$ if and only if $b_2-\sigma_1^{-2}b_1\tau \in [-\rho,\rho]$,
which are the two different cases the definition of $d_2^2$ is separated into.
Now, if $P_L\in Seg(a,b,\theta)$, then, obviously, $d_2^2=0$. Next assume that $P_L\not\in Seg(a,b,\theta)$. Then
\begin{eqnarray*}
d_2^2
& = & \min_\pm \norm{L(\pm \rho)-P_L}_2^2\\
& = & \min_\pm \; (\pm \rho-b_2)^2\sigma_1^2+\sigma_1^{-4}b_1^2\tau^2\sigma_1^2+2(\rho-b_2)\sigma_1^{-2}b_1\tau\sigma_1^{2}\\
& = & \min_\pm \; ((\pm \rho-b_2)\sigma_1-\sigma_1^{-1}b_1\tau)^2.
\end{eqnarray*}
\qed

Now define the {\it ray integral}
\[
R_N(x_0,y_0)= \int_{y_0}^\infty \langle |(x_0,t)|\rangle^{-N} dt,
\]
i.e., we integrate along the vertical ray $\cR(x_0,y_0)$ whose `lowest' point is $(x_0,y_0)$.
The geometry of the previous lemma allows to
control the curvelet coefficient of a linear singularity by a ray integral,
properly deployed. Farther below we will prove:

\begin{lemma} \label{lemm:raybound}
Let
\[
    \langle  w \cL , \gamma_{a,b,\theta} \rangle = \int_{-\rho}^\rho w(x_2/\rho) \gamma_{a,b,\theta}(0,x_2) dx
\]
where $\| w \|_\infty \leq 1$. Then
\[
    | \langle  w \cL , \gamma_{a,b,\theta} \rangle| \leq a^{-3/4} \cdot \sigma_1^{-1} \cdot  R_N(d_1, \sigma_1 d_2) .
\]
\end{lemma}

The next lemma gives a bound on the ray integral which, combined
with the last lemma, finishes the proof of Lemma \ref{lemma:linecurvelet_estimate_close0}.

\begin{lemma} \label{lemm:rayintegral}
For $y_0 \ge 0$,
\beq \label{eq:PLPS_3}
R_N(x_0,y_0) \le \pi \cdot \langle|x_0|\rangle^{-1} \cdot  \langle|(x_0,y_0)|\rangle^{2-N}.
\eeq
\end{lemma}

{\bf Proof of Lemma \ref{lemm:rayintegral}.} For $\beta \in (0,1)$,
\[
\int_0^\infty|f(t)| dt\le (\sup_{t \in (0,\infty)}|f(t)|^\beta)\cdot\int_0^\infty|f(t)|^{1-\beta}dt.
\]
Now setting $(1-\beta)N=2$ and $f(t) = \langle|(x_0,y_0+t)|\rangle^{-N}$,
we have
\[
R_N(x_0,y_0) \le (\sup_{v \in \cR(x_0,y_0)}\langle|v|\rangle^{2-N})
\cdot\int_0^\infty|\langle|(x_0,y_0+t)|\rangle^{-2}dt.
\]
Since
\[
\int_{-\infty}^\infty|\langle|(x_0,y)|\rangle^{-M} dy
=\langle|x_0|\rangle^{-M} \cdot \int_{-\infty}^\infty \langle y/\langle x_0 \rangle \rangle^{-M} dy
=\langle|x_0|\rangle^{-M+1} \cdot \int_{-\infty}^\infty \langle t \rangle^{-M} dt,
\]
setting $M=2$ and recalling $\pi = \int_{-\infty}^\infty (1+t^2)^{-1} dt$, it follows that
\[
\int_0^\infty|\langle|(x_0,y_0+t)|\rangle^{-2}dt
\le \pi \cdot \langle|x_0|\rangle^{-1}.
\]
Meanwhile, since $y_0 \ge 0$,
\[
\sup_{v \in \cR(x_0,y_0)}\langle|v|\rangle^{2-N} =
\langle|(x_0,y_0)|\rangle^{2-N}.
\]
This proves \eqref{eq:PLPS_3}. \qed


{\bf Proof of Lemma \ref{lemm:raybound}.}
By Lemma \ref{lemma:estimate_gamma}, Lemma \ref{lemma:PLPS}, and using the fact that $\norm{w}_\infty \le 1$,
\begin{eqnarray}\nonumber
|\langle w\cL ,  \gamma_{a,b,\theta}  \rangle |
& = & \Big| \int_{-\rho}^\rho w_2(y/\rho) \gamma_{a,b,\theta}(0,y) dy\Big|\\ \nonumber
& \le & \int_{-\rho}^\rho |\gamma_{a,b,\theta}(0,y)| dy\\ \nonumber
& \le & \int_{Seg(a,b,\theta)}   a^{-3/4} \langle |v| \rangle^{-N}dv
\end{eqnarray}
where we used an affine  transformation of variables to turn the anisotropic norm $| (0,y) |_{a,\theta}$
into the Euclidean norm $|v|$; the same transformation turns $\{0\} \times [-\rho,\rho]$
into  $Seg(a,b,\theta)$.
In the final expression, the integral is along a non-unit-speed curve traversing
$Seg(a,b,\theta)$, at speed $\sigma_1$.
Now let $Ray(a,b,\theta)$ denote the ray starting from $P_S$ and initially
traversing $Seg(a,b,\theta)$.  We continue
with
\begin{eqnarray}\nonumber
c_N \cdot  a^{-3/4}  \cdot  \int_{Seg(a,b,\theta)}  \langle |v| \rangle^{-N} dv  & \le & c_N \cdot  a^{-3/4}  \cdot   \int_{Ray(a,b,\theta)}  \langle |v| \rangle^{-N} dv  \\ \label{unitspeed}
& =& c_N \cdot  a^{-3/4}  \cdot \sigma_1^{-1}  \int_{\sigma_1 \cdot Ray(a,b,\theta)}  \langle |w| \rangle^{-N} dw  \\ \nonumber
& = & c_N \cdot a^{-3/4}   \cdot \sigma_1^{-1}  \int_{\sigma_1d_2}^\infty \langle |(d_1,t)|\rangle^{-N} dt \\ \nonumber
& = &  c_N \cdot a^{-3/4}   \cdot \sigma_1^{-1}  \cdot R_N(d_1,\sigma_1 d_2).
\end{eqnarray}
In (\ref{unitspeed}), the integral involves a unit-speed curve traversing $Ray(a,b,\theta)$,
which explains the appearance of the speed factor $\sigma_1$.  \qed



\subsubsection{Proof of  Lemma \ref{lemma:linecurvelet_estimate_other}}
\label{subsec:proof_of_linecurvelet_estimate_other}

By definition of the line singularity $w\cL$ and by \eqref{eq:actionofwL}, we can
rewrite $\ip{w\cL}{\gamma_{a,b,\theta}}$ in the following way:
\beq \label{eq:linecurvelet_estimate_other1}
2\pi\ip{w\cL}{\gamma_{a,b,\theta}} = \int (\hat{w} \star \hat{\gamma}_{a,b,\theta})(\xi_1,0) d\xi_1.
\eeq
Since
\[
(\hat{w} \star \hat{\gamma}_{a,b,\theta})(\xi)
= \int \rho \cdot \hat{w}_2(\rho(\xi_2-\eta_2)) \delta_0(\xi_1-\eta_1) \hat{\gamma}_{a,b,\theta}(\eta) d\eta,
\]
it follows that
\begin{eqnarray*}
(\hat{w} \star \hat{\gamma}_{a,b,\theta})(\xi_1,0)
& = & e^{i b_1 \xi_1} \int \rho \cdot \hat{w}_2(-\rho\eta_2) \hat{\gamma}_{a,0,\theta}(\xi_1,\eta_2)
e^{i b_2 \eta_2} d\eta_2.
\end{eqnarray*}
By \eqref{eq:linecurvelet_estimate_other1}, this implies
\[
2\pi\ip{w\cL}{\gamma_{a,b,\theta}}
= \int e^{i b_1 \xi_1} \Big[\int \rho \cdot \hat{w}_2(-\rho\eta_2) \hat{\gamma}_{a,0,\theta}(\xi_1,\eta_2)
e^{i b_2 \eta_2} d\eta_2 \Big] d\xi_1.
\]
Repeatedly applying integration by parts, and incorporating analyst's brackets $\langle| \, \cdot \, |\rangle$
as in the proof of Lemma \ref{lemma:pointwavelet}, we obtain
\beq \label{eq:linecurvelet_estimate_other2}
2\pi\absip{w\cL}{\gamma_{a,b,\theta}} \le \langle |b_1| \rangle^{-L} \cdot \langle |b_2| \rangle^{-M}
\cdot \norm{h_{L,M}}_{L^1(\bR)},
\eeq
where
\[
h_{L,M}(\xi_1) = \rho \cdot \int D^{L,M}\Big(\hat{w}_2(-\rho \eta_2) \hat{\gamma}_{a,0,\theta}(\xi_1,\eta_2) e^{i b_2 \eta_2}\Big) d\eta_2
\]
and for some `nice' $f \in L^2(\bR^2)$,
\[
D^{L,M}f(\eta_1,\eta_2) = \left(\frac{\partial}{\partial \eta_1}\right)^L\left(\frac{\partial}{\partial \eta_2}\right)^M
f(\eta_1,\eta_2).
\]

Next, we will estimate the term $|h_{L,M}(\xi_1)|$ from \eqref{eq:linecurvelet_estimate_other2},
and prove that
\beq \label{eq:linecurvelet_estimate_other3}
|h_{L,M}(\xi_1)| \le c \cdot a^{3/4} \cdot e^{-\rho \frac{|\sin \theta|}{2a}} \cdot
(a^{1/2} |\sin \theta| + a |\cos \theta|)^L \cdot (\rho + a^{1/2} |\cos \theta| + a |\sin \theta|)^M.
\eeq
Let $\Xi_{a,\theta}(\xi_1)$ denote the support of the function $\xi_1 \mapsto
D^{L,M}(\hat{w}(\rho \eta_2) \hat{\gamma}_{a,0,\theta}(\xi_1,\eta_2) e^{i b_2 \eta_2})$.
Then $h_{L,M}$ can be written as
\beq \label{eq:linecurvelet_estimate_other4}
h_{L,M}(\xi_1) = \rho \cdot \int_{\Xi_{a,\theta}(\xi_1)} D^{L,M}\Big(\hat{w}_2(-\rho \eta_2)
\hat{\gamma}_{a,0,\theta}(\xi_1,\eta_2) e^{i b_2 \eta_2}\Big) d\eta_2.
\eeq
We next rewrite the integrand as
\begin{eqnarray*}
\lefteqn{D^{L,M}\Big(\hat{w}_2(-\rho \eta_2) \hat{\gamma}_{a,0,\theta}(\xi_1,\eta_2) e^{i b_2 \eta_2}\Big)}\\
& = & \sum_{m=0}^M \left(\hspace*{-0.15cm}\begin{array}{c}M\\m\end{array}\hspace*{-0.15cm}\right)
\hat{w}^{(m)}_2(-\rho \eta_2) (-\rho)^m D^{L,M-m}(\hat{\gamma}_{a,0,\theta}(\xi_1,\eta_2)e^{i b_2 \eta_2}).
\end{eqnarray*}
This allows us to estimate $|h_{L,M}(\xi_1)|$ using \eqref{eq:linecurvelet_estimate_other4} and
\eqref{eq:exp_w2} by
\begin{eqnarray}\nonumber
|h_{L,M}(\xi_1)|
& \le & \sum_{m=0}^M  \left(\hspace*{-0.15cm}\begin{array}{c}M\\m\end{array}\hspace*{-0.15cm}\right)
\Big|\int_{\Xi_{a,\theta}(\xi_1)} \hat{w}^{(m)}_2(-\rho \eta_2)
(-\rho)^m D^{L,M-m}(\hat{\gamma}_{a,0,\theta}(\xi_1,\eta_2)e^{i b_2 \eta_2})d\eta_2\Big|\\ \nonumber
& \le & \sum_{m=0}^M  \left(\hspace*{-0.15cm}\begin{array}{c}M\\m\end{array}\hspace*{-0.15cm}\right)
\rho^{m+1} \cdot \norm{\hat{w}^{(m)}_2(\rho \, \cdot)}_{L^1[|\sin \theta|/(2a),\infty)}
N^{L,M-m}(a,\theta)\\ \label{eq:linecurvelet_estimate_other5}
& \le & c \cdot \sum_{m=0}^M  \left(\hspace*{-0.15cm}\begin{array}{c}M\\m\end{array}\hspace*{-0.15cm}\right)
\rho^{m+1} \cdot e^{-\rho \frac{|\sin \theta|}{2a}} \cdot N^{L,M-m}(a,\theta),
\end{eqnarray}
where
\[
N^{L,M-m}(a,\theta) = \norm{D^{L,M-m} \hat{\gamma}_{a,0,\theta}(\xi_1,\eta_2)}_{L^\infty(\Xi_{a,\theta}(\xi_1))}.
\]
Since, by simple decay estimates,
\[
|D_1^L \hat{\gamma}_{a,0,\theta}(\eta_1,\eta_2)| \le C_L \cdot a^{3/4} \cdot (a|\cos \theta|+a^{1/2}|\sin \theta|)^L
\]
and
\[
|D_2^M \hat{\gamma}_{a,0,\theta}(\eta_1,\eta_2)| \le C_M \cdot a^{3/4} \cdot (a^{1/2}|\cos \theta|+a|\sin \theta|)^M,
\]
the term $N^{L,M-m}(a,\theta)$ can be estimated by
\[
N^{L,M-m}(a,\theta) \le c_{L,M}\cdot a^{3/4} \cdot (a|\cos \theta|+a^{1/2}|\sin \theta|)^L(a^{1/2}|\cos \theta|+a|\sin \theta|)^M.
\]
Combining this finding with \eqref{eq:linecurvelet_estimate_other5}
proves \eqref{eq:linecurvelet_estimate_other3}.

Thus, in particular, by the support of
the function $h_{L,M}$, the $L^1$-norm of this function can be
estimated as
\begin{eqnarray*}
\lefteqn{\norm{h_{L,M}}_{L^1(\bR)}}\\
& \le & \hspace*{-0.2cm} c \cdot a^{-1} \cdot |\cos \theta| \cdot a^{3/4} \cdot e^{-\rho \frac{|\sin \theta|}{2a}} \cdot
(a^{1/2} |\sin \theta| + a |\cos \theta|)^L \cdot (\rho + a^{1/2} |\cos \theta| + a |\sin \theta|)^M.
\end{eqnarray*}
Combining this estimate with \eqref{eq:linecurvelet_estimate_other2} yields
\begin{eqnarray*}
\absip{w\cL}{\gamma_{a,b,\theta}}
& \le & c_{M,L} \cdot a^{-1/4} \cdot |\cos \theta| \cdot e^{-\rho \frac{|\sin \theta|}{2a}} \cdot \langle |b_1| \rangle^{-L} \cdot
(a^{1/2} |\sin \theta| + a |\cos \theta|)^L\\
& &  \cdot \langle |b_2| \rangle^{-M} \cdot (\rho + a^{1/2} |\cos \theta| + a |\sin \theta|)^M,
\end{eqnarray*}
as claimed. \qed

\subsection{Proofs of Results from Section \ref{sec-Sec2j}}
\label{subsec:proofs_5}

\subsubsection{Proof of Lemma \ref{lemm:t_estimate}}
\label{subsec:proofs_5_1}

\noindent
{\bf Proof.}
Below, various constants will appear, which for simplicity shall all be denoted by $c$.
Also we write $b_{\t{k}}$, $b_{k'}$, etc. in place of the full notation $b_{\t{j},\t{k},\t{\ell}} = R_{\theta_{j,\ell}} D_{2^{-j}} k$, etc.
Throughout the proof we associate $\t{\eta}$ with the triple $(\t{j},\t{k},\t{\ell})$ and similarly for  $\eta'$ and $(j',k',\ell')$.

Suppose first that
$j' \geq \t{j}$.
Then
\begin{eqnarray}  \nonumber
\lefteqn{\omega(\eta',\t{\eta})}\\ \nonumber
& = & 2^{|\t{j}-j'|}
\left(1+\min\{2^{\t{j}},2^{j'}\} \left[|\theta_{\t{j},\t{\ell}}-\theta_{j',\ell'}|^2
 + |b_{\t{k}}-b_{k'}|^2 + |\ip{e_{\eta'}}{b_{\t{k}}-b_{k'}}|\right]\right)\\ \nonumber
& \le & c \cdot \left( 2^{j'-\t{j}} + | \t{\ell}- 2^{(\t{j}-j')/2} \ell'|^2
+ 2^{\t{j}} \cdot |b_{\t{k}}-b_{k'}|^2 + 2^{\t{j}} |\ip{e_{\eta'}}{b_{\t{k}}-b_{k'}}| \right).
\end{eqnarray}
Similarly, if  $j' < \t{j}$:
\begin{eqnarray}\label{eq:estimateomega2}
\omega(\eta',\t{\eta})
& \le & c \cdot \left(2^{\t{j}-j'} + |2^{(j'-\t{j})/2} \t{\ell}- \ell'|^2
+ 2^{j'} \cdot |b_{\t{k}}-b_{k'}|^2 + 2^{j'} |\ip{e_{\eta'}}{b_{\t{k}}-b_{k'}}| \right).
\end{eqnarray}

Using the model (\ref{defCurvClust0}), we now seek to prove that, for sufficiently large $j$,
\beq \label{eq:t_estimate_2}
\#\{\eta' : |\t{M}_j^i(\eta',\eta)| \ge c_N \cdot n_j^{-N}\} = \#\{\eta' : \omega(\eta',(\t{j},\t{k},\t{\ell})) \le   n_j\} \ge n_j.
\eeq

We first study the case $j'  < \t{j}$. Let $\tilde{\omega}(\eta',\t{\eta})$ denote the RHS of \eqref{eq:estimateomega2}.
Note that, if we can prove \eqref{eq:t_estimate_2}
with $\tilde{\omega}(\eta',\t{\eta})$ in place of $\omega(\eta', \t{\eta})$, this immediately
implies the original claim of \eqref{eq:t_estimate_2}.
We repeatedly use the (trivial)
\begin{lemma} \label{lem:countLattice}
For each fixed $x \in \bR$, the set of $k \in \bZ$ satisfying $|k - x | \leq R$ has cardinality $\geq R-1 $.
\end{lemma}

Define $C_{j,j',\t{j}} = n_j/c - 2^{j'-\t{j}}$, where
$c$ is the constant in (\ref{eq:estimateomega2}).
The condition $\tilde{\omega} \leq n_j$
is equivalent to 
\begin{eqnarray}
\|2^{(j'-\t{j})/2} \t{\ell}- \ell'|^2
+ 2^{j'} \cdot |b_{\t{k}}-b_{k'}|^2 + 2^{j'} |\ip{e_{\eta'}}{b_{\t{k}}-b_{k'}}| & \le & C_{j,j',\t{j}}. \label{equivCond}
\end{eqnarray}
We next derive conditions making each of the three terms  smaller than  $C_{j,j',\t{j}}/3$,
thus implying (\ref{equivCond}).
Fix $j'$. By Lemma \ref{lem:countLattice}, there are at least
$C_{j,j',\t{j}}^{1/2}/\sqrt{3} -1$ integer values $\ell' \in \bZ$ obeying
\beq \label{condell}
|2^{(j'-\t{j})/2} \t{\ell}- \ell'|^2 \le \frac13 C_{j,j',\t{j}}.
\eeq
Now, define $\tilde{x} = \tilde{x}(\t{\eta},\eta') = R_{\theta_{j',\ell'}}^{-1} b_{\t{k}}$.
Since $R_{\theta_{j',\ell'}}$ is an isometry,
\[
   |b_{\t{k}} - b_{k}|^2 = | \t{x} -D_{2^{-j'}}k'|^2.
\]
At the same time,
\[
   |\ip{e_{\eta'}}{b_{\t{k}}-b_{k'}}| =  |\t{x}_1 - 2^{-j'} k_1' | .
\]
Now there are at least
$ C_{j,j',\t{j}}^{1/2}/\sqrt{6} - 1$ integers $k_2' \in\bZ$  obeying
\beq \label{condk2}
   | 2^{j'/2}\t{x}_2  - k'_2 | \leq   C_{j,j',\t{j}}^{1/2}/\sqrt{6} .
\eeq
For large $j'$,
there are at least  $ C_{j,j',\t{j}}/3 - 1$
integers $k_1' \in\bZ$  obeying {\it both}
\[
   | 2^{j'} \t{x}_1 - k'_1 | \leq   2^{j'/2} \cdot C_{j,j',\t{j}}^{1/2}/\sqrt{6}
\]
{\it and}
\beq \label{condk1b}
   | 2^{j'} \t{x}_1  - k'_1 | \leq   C_{j,j',\t{j}}/3.
\eeq
Every pair  $k'= (k'_1,k'_2) \in \bZ^2$
satisfying the conditions (\ref{condk2})-(\ref{condk1b}) simultaneously,
satisfies
\beq \label{condboth}
|b_{\t{k}} - b_{k}|^2  \leq  C_{j,j',\t{j}}/3 \quad \mbox{and} \quad  |\ip{e_{\eta'}}{b_{\t{k}}-b_{k'}}|  \leq  C_{j,j',\t{j}}/3 .
\eeq

Combining the above displays, we have  at least
$( C_{j,j',\t{j}}/3 - 1) \cdot  (C_{j,j',\t{j}}^{1/2}/\sqrt{6} -1)$  points $k'= (k'_1,k'_2) \in \bZ^2$
satisfying  (\ref{condboth}). Every pair $(k',\ell')$ satisfying (\ref{condboth}) and (\ref{condell})
satisfies  (\ref{equivCond}). So
focusing just on $j' = j$, we obtain a large number of pairs $(k',\ell')$  satisfying (\ref{equivCond}):
$ \geq c \cdot C_{j,j',\t{j}}^{2} \sim c \cdot n_j^{2} \gg n_j$ such pairs.
In particular, for all large $j$, inequality $\omega(\t{\eta},\eta') \leq n_j $ is satisfied by at
least $ n_j$ triples $\eta' = (j',k',\ell')$.

For the case $j' < \t{j}$, using \eqref{eq:estimateomega2}, we can similarly prove that $\omega \leq n_j$
is satisfied by at least $ n_j$ triples $(j',k',\ell')$.

Finally, we observe that $\omega(\t{\eta},\eta') \leq n_j$ can only hold if
\beq \label{eq:finalcount1}
|j'-\t{j}| \le c \cdot \log n_j, \quad
|2^{(j'-\t{j})/2} \t{\ell}- \ell'| \le c  \cdot n_j
\eeq
and
\beq \label{eq:finalcount2}
2^{j'} \cdot \max\{|b_{\t{k}}-b_{k'}|^2,|b_{\t{k}}-b_{k'}|\} \le c \cdot n_j.
\eeq
Concluding, each $\eta'$ contained in the sets in \eqref{eq:t_estimate_2} must satisfy
both \eqref{eq:finalcount1} and \eqref{eq:finalcount2}.
\qed

\subsubsection{Proof of Lemma \ref{lemm:lboundk}}
\label{subsec:proofs_Lemma_7_9}

For $\eta \in \t{\cS}_j$, we have $\theta = 0$
and $|b_2| < n_j/2^{j}$ while $|b_1| < 2 \rho /2^{j/2}$.
WLOG suppose that, for the patch $i$ in question, we
have that the $\eta' \in S_j^i$ of interest has
$b' = R_{\theta'} D_{2^{-j'}} k'$ for some $\theta$
obeying $|\theta| \leq c n_j/2^{j/2}$.  For such a pair
$(\eta,\eta')$, note that
\[
k = D_{2^{j}} R_0^{-1} b, \qquad k' = D_{2^{j'}} R_{\theta'}^{-1} b' ;
\]
it is also convenient to define $x' = D_{2^{j-j'}} k'$.
Then
\[
   k_2 - x'_2 = 2^{j/2} \left( (1-\cos(\theta')) b_2' + (b_2 - b'_2) - \sin(\theta') b'_1 \right).
\]
Now $|b'_1| < C$ for $\eta' \in \Fwd(\eta)$, and, from Lemma 7.8 we infer
\[
   |1 - \cos(\theta') | < c n_j^2/2^{j} , \qquad |b'_2 - b_2| < c n_j / 2^j , \qquad  |\sin(\theta')| < c n_j / 2^{j/2}  .
\]
Combining these,
\[
   |k_2 - x'_2 | < 2c n_j^2 ,
\]
so $ |x'_2| > |k_2| - 2c n_j^2 $, and
\[
   |k'_2| = 2^{(j'-j)/2 } |x'_2|  \geq 2^{- |j - j'|/2} \cdot \left( |k_2| - 2c n_j^2 \right) .
\]
Now $2^{|j - j'|} < cn_j$ for $\eta' \in \Fwd(\eta)$, so
$|k(\eta)| \leq |k_2(\eta)| + |k_1(\eta)| \leq |k_2(\eta)| + c n_j$, for $\eta \in \t{\cS}_j$.
So
\[
  |k(\eta')| \geq |k'_2| \geq  2^{- |j - j'|/2} \cdot \left( |k_2| - 2c n_j^2 \right) , \qquad j > j_0.
\]
\qed

\subsubsection{Proof of Lemma \ref{lemm:sumlemma}}
\label{subsec:proofs_Lemma_7_10}

Let $Q_\ell$ denote the square of sidelength $2^{\ell}+1$ centered at the origin.
Within each annulus  $A_\ell = Q_\ell - Q_{\ell-1}$  there are, for $\ell \geq  2$,
fewer than $2^{2\ell}$ points in the annulus and each one has $\ell^2$ norm at least $2^{\ell-2}$.
Partition the sum $\sum_{\bZ^2} =  \sum_{Q_1 + A_2 + A_3 + \dots }$.
Starting at $\ell=2$ each annulus contributes at most $2^{2\ell} \cdot \langle(a 2^{\ell-2} - b)_+\rangle^{-N}$
to the sum, and the inner square $Q_1$ contributes at most 9. Then
\begin{eqnarray*}
   \sum_{k \in \bZ^2} \langle (a |k| - b )_+ \rangle^{-N}  &\leq&  9 +
   \sum_{m=2}^\infty   2^{2\ell}   \cdot  \langle (a 2^{\ell-2} - b )_+ \rangle^{-N}  \\
      &=& 9 + 16 \cdot  \left( \sum_{m=0}^\infty  2^{2m}   \cdot  \langle (a 2^{m} - b )_+ \rangle^{-N} \right) .
\end{eqnarray*}

Let $m_0$ satisfy $2 b \geq a 2^{m_0} \geq b$. Now
\begin{eqnarray*}
T_1 := \sum_{m=0}^{m_0} 2^{2m} \langle
 |a 2^m - b |_+ \rangle^{-N} & \leq & \sum_{m=0}^{m_0} 2^{2m}  \leq  2^{2m_0+1} .
\end{eqnarray*}
Now $2^{2m_0} \leq (2b/a)^2$, so $2^{m_0} \leq 2b/a$.
So  $T_1 \leq 2 \cdot (2b/a)^2$.

On the other hand, $a 2^m - b \geq 0$ for $m \geq m_0$,
while $a 2^{m-m_0} \leq a2^m - b $. Then
\begin{eqnarray*}
T_2 := \sum_{m>m_0} 2^{2m} \langle | a 2^{m-m_0} | \rangle^{-N} & \leq & \sum_{m=m_0+1}^\infty 2^{2m} (a 2^{m-m_0})^{-N} \\
    &\leq& 2^{2(m_0+1)} (a 2^{-m_0})^N \sum_{h=0}^\infty 2^{2h} 2^{-hN} \\
    &=& 2^{m_0(2-N)}  \cdot a^{N} \cdot 2 \cdot (1 - 2^{-(N-2)})^{-1} \\
    &\leq&   (b/a)^{2-N} \cdot a^{N} \cdot 4  \\
      &\leq & 4 \cdot (2b/a)^2 \cdot b^{-N}
\end{eqnarray*}
Hence
\[
T_1 + T_2 \leq  (b/a)^2 \cdot \left(  8 + 16 b^{-N} \right).
\]
Combining these displays gives the lemma, with explicit constants. \qed

\end{document}